\documentclass[11pt]{article}

\usepackage{latexsym,amsmath,amsfonts,amssymb,amstext,theorem,euscript,array,enumerate,mathrsfs}
\usepackage{color}
\usepackage{comment}
\usepackage{hyperref}

\newtheorem{Theorem}{Theorem}[section]
\newtheorem{Definition}{Definition}[section]

\newtheorem{Lemma}{Lemma}[section]

\newtheorem{Remark}{Remark}[section]

\numberwithin{equation}{section}

\def\esssup_#1{\underset{#1}{\mathrm{ess\,sup\, }}}
\def\essinf_#1{\underset{#1}{\mathrm{ess\,inf\, }}}

\def \trans{^{\scriptscriptstyle{\intercal}}}

\def\sqr#1#2{{\vcenter{\vbox{\hrule height .#2pt \hbox{\vrule
 width .#2pt height#1pt \kern#1pt \vrule
width .#2pt} \hrule height .#2pt}}}}

\def\ds{\begin{displaystyle}}
\def\eds{\end{displaystyle}}
\def\dis{\displaystyle }
\def\<{\langle }
\def\>{\rangle }

\def \N{\mathbb{N}}
\def \R{\mathbb{R}}

\def \E{\mathbb{E}}
\def \F{\mathbb{F}}

\def \P{\mathbb{P}}
\def \Q{\mathbb{Q}}

\def \S{\mathbb{S}}

\def \Ac{{\cal A}}
\def \Bc{{\cal B}}
\def \Cc{{\cal C}}

\def \Ec{{\cal E}}
\def \Fc{{\cal F}}

\def \Hc{{\cal H}}

\def \Jc{{\cal J}}
\def \Kc{{\cal K}}
\def \Lc{{\cal L}}
\def \Pc{{\cal P}}
\def \Mc{{\cal M}}

\def \Rc{{\cal R}}
\def \Sc{{\cal S}}

\def \Vc{{\cal V}}
\def \Wc{{\cal W}}

\def \Vc{{\cal V}}

\def\calb{{\cal B}}

\def\calf{{\cal F}}

\def\calp{{\cal P}}

\def \eps{\varepsilon}

\def \ep{\hbox{ }\hfill$\Box$}

\def\upv{\text{\Large$\upsilon$}}

\def\Dt#1{\Frac{\partial #1}{\partial t}}

\def \Sum{\displaystyle\sum}

\def \Frac{\displaystyle\frac}

\def \Sup{\displaystyle\sup}

\def \b1{\bf{1}}

\def\reff#1{{\rm(\ref{#1})}}

\def\beqs{\begin{eqnarray*}}
\def\enqs{\end{eqnarray*}}
\def\beq{\begin{eqnarray}}
\def\enq{\end{eqnarray}}

\allowdisplaybreaks

\begin{document}

\title{Randomized filtering  and Bellman equation in Wasserstein space for
partial observation control problem}

\author{Elena BANDINI\thanks{Dipartimento di Economia e Finanza, Libera Universit\`a degli Studi Sociali ``Guido Carli'', Roma, Italy; e-mail: \texttt{ebandini@luiss.it}}
\and
Andrea COSSO\thanks{Politecnico di Milano, Dipartimento di Matematica, via Bonardi 9, 20133 Milano, Italy; e-mail: \texttt{andrea.cosso@polimi.it}}
\and
Marco FUHRMAN\thanks{Politecnico di Milano, Dipartimento di Matematica, via Bonardi 9, 20133 Milano, Italy; e-mail: \texttt{marco.fuhrman@polimi.it}}
\and
Huy\^{e}n PHAM\thanks{Laboratoire de Probabilit\'es et Mod\`eles Al\'eatoires, CNRS, UMR 7599, Universit\'e Paris Diderot, and CREST-ENSAE; e-mail: \texttt{pham@math.univ-paris-diderot.fr}}
}

\maketitle

\begin{abstract}
We study a stochastic optimal control problem for a partially observed diffusion. By using  the control randomization method in \cite{BCFP16b}, we
prove a corresponding 
randomized dynamic programming principle (DPP) for the value function, which is obtained from a flow property of an associated filter process. This DPP is the key step towards our main result: 
a characterization of the value function
of the partial observation control problem as the unique viscosity solution to the corresponding dynamic programming Hamilton-Jacobi-Bellman (HJB) equation. The latter
 is formulated as a new, fully non linear partial differential
 equation on the Wasserstein space of probability measures.
An important feature of our approach is that it does not require any non-degeneracy condition on the diffusion coefficient, and no condition is imposed to guarantee 
existence of a density for the filter process solution to the controlled Zakai equation, as usually done for the se\-parated problem.
%and  leading to a fully nonlinear equation with unbounded terms in the infinite dimensional space of density functions, see  \cite{Be}, \cite{par88}, \cite{lio89}, \cite{gozswi15}.
Finally, we give an explicit solution to our HJB equation in the case of  a partially observed non Gaussian linear quadratic model.
\end{abstract}

\vspace{5mm}


\vspace{5mm}


\maketitle

\date{}

\newpage

%\section{Introduction}

\section{Introduction and problem formulation}  \label{S:Formulation}

\setcounter{equation}{0} \setcounter{Assumption}{0}
\setcounter{Theorem}{0} \setcounter{Proposition}{0}
\setcounter{Corollary}{0} \setcounter{Lemma}{0}
\setcounter{Definition}{0} \setcounter{Remark}{0}

We formulate  the optimal control problem of partially observed diffusion process.
Let $(\bar\Omega,\bar\calf,\bar\P)$ be a complete probability space, on which a $(m+d)$-dimensional Brownian motion  $B$ is defined.
%and denote by $\bar\F$ $=$ $(\bar\Fc_t)_{t\geq 0}$ its natural filtration (namely, the right-continuous and $\bar\P$-complete filtration generated by $B$), %augmented with an independent $\sigma$-algebra $\bar\Fc^o$ $\subset$ $\bar\Fc$.
We will distinguish the $m$-dimensional and $d$-dimensional components of $B$  by writing $B=(V,W)$ where $W$ plays the role of the observation process. We suppose that there exist three probability spaces $(\Omega^o,\Fc^o,\P^o)$, $(\Omega^1,\Fc^1,\P^1)$, $(\Omega,\Fc,\P)$ such that $\bar\Omega=\Omega^o\times\Omega^1\times\Omega^2$, $\bar\Fc$ is the completion of $\Fc^o\otimes\Fc^1\otimes\Fc$ with respect to $\P^o\otimes\P^1\otimes\P$, and $\bar\P$ is the extension of $\P^o\otimes\P^1\otimes\P$ to $\bar\Fc$.
%Moreover, $\bar\Fc^o=\{E\times\Omega^1\times\Omega^2\colon E\in\Fc^o\}$ is the canonical extension of $\Fc^o$ to $\bar \Omega$.
In addition, $\Omega^o$ is a Polish space, $\Fc^o$ its Borel $\sigma$-algebra, $\P^o$ an atomless probability measure on $(\Omega^o,\Fc^o)$. The space $(\Omega^o,\Fc^o,\P^o)$ supports the initial condition $\xi$ of the state equation \eqref{dynX} introduced below, while $(\Omega^1,\Fc^1,\P^1)$ supports $V$, and $(\Omega,\Fc,\P)$ supports $W$. This is graphically described by the following figure:
\[
\bar\Omega \ := \ \underset{\!\!\!\underset{\xi}{\uparrow}}{\Omega^o} \times \underset{\!\!\!\underset{V}{\uparrow}}{\Omega^1} \times \underset{\underset{W}{\uparrow}}{\Omega}
\]
An element $\bar\omega$ $\in$ $\bar\Omega$ is then written as $\bar\omega$ $=$ $(\omega^o,\omega^1,\omega)$, and
we extend canonically $V,W$ on $\bar\Omega$ by setting $V(\omega^o,\omega^1,\omega)$ $=$ $V(\omega^1)$,
$W(\omega^o,\omega^1,\omega)$ $=$ $W(\omega)$, and similarly any random variable $\xi$ on $(\Omega^o,\Fc^o,\P^o)$ is extended to
$(\bar\Omega,\bar\Fc,\bar\P)$ by setting $\xi(\omega^o,\omega^1,\omega)$ $=$ $\xi(\omega^o)$.
We denote by $\bar\F$ $=$ $(\bar\Fc_t)_{t\geq 0}$ the natural filtration  of $B$
(namely, the right-continuous and $\bar\P$-complete filtration generated by $B$), augmented with the  independent $\sigma$-algebra $\Fc^o$. The partial observation character of our 
problem is modeled as follows:
we denote by  $\F^W=(\calf^W_t)_{t\ge 0}$   the $\P$-completion of the filtration generated by $W$
and we define as
an admissible control process  any
$\F^W$-progressive process $\alpha$ with values in some given Borel space $A$.  The set of admissible control processes is denoted by $\Ac^W$.  Given $t$
$\in$ $[0,T]$,  $\xi$ $\in$ $L^2(\Fc^o;\R^n)$,  the set of square-integrable $\Fc^o$-measurable random variables, valued in $\R^n$,
and $\alpha$ $\in$ $\Ac^W$, the controlled state equation is
\beq \label{dynX}
dX_s^{t,\xi,\alpha} &=& b(X_s^{t,\xi,\alpha},\alpha_s) ds + \sigma(X_s^{t,\xi,\alpha},\alpha_s) dB_s,  \;\;\;  t \leq s \leq T, \;\;\;  X_t^{t,\xi,\alpha} = \xi.
\enq
The coefficients  $b$ and $\sigma$ $=$ $(\sigma_{_V} \; \sigma_{_W})$  are deterministic measurable functions from $\R^n\times A$ into $\R^n$ and $\R^{n\times (m+d)}$, and assumed to satisfy the following standing assumptions.

\vspace{2mm}

\hspace{-7mm} {\bf (H1)}   There exists some positive constant $C_{1}$ such that for all $x,x'$ $\in$ $\R^n$, $a$ $\in$ $A$,
\beqs
|b(x,a)-b(x',a)| +  |\sigma(x,a)-\sigma(x',a)| & \leq & C_{1} |x-x'| \\
|b(0,a)| + |\sigma(0,a)| & \leq & C_{1}.
\enqs

Under {\bf (H1)}, it is shown by standard arguments that there exists a unique strong solution $X^{t,\xi,\alpha}$ to \reff{dynX}, which is $\bar\F$-adapted,  and satisfies
\beq \label{estimX}
\sup_{\alpha\in\Ac^W} \bar\E \Big[ \sup_{t\leq s\leq T} \big|X_s^{t,\xi,\alpha}\big|^2\Big] & \leq & C \big( 1 + \bar\E|\xi|^2 \big) \; < \; \infty,
\enq
for some positive constant $C$, where $\bar\E$ denotes the expectation with respect to $\bar\P$.

The aim is to maximize, over all admissible control processes $\alpha$ $\in$ $\Ac^W$,  the gain functional:
\beqs
\Jc(t,\xi,\alpha) & := & \bar\E \Big[ \int_t^T f(X_s^{t,\xi,\alpha},\alpha_s) ds + g(X_T^{t,\xi,\alpha}) \Big],
\enqs
where $f$ $:$ $\R^n\times A$ $\rightarrow$ $\R$, and $g$ $:$ $\R^n$ $\rightarrow$ $\R$ are continuous functions satisfying the quadratic growth condition:

\vspace{2mm}

\hspace{-7mm} {\bf (H2)} There exists some positive constant $C_{2}$ such that for all $(x,a)$ $\in$ $\R^n\times A$,
\beqs
|f(x,a)| + |g(x)| & \leq & C_{2}\big( 1 + |x|^2 \big).
\enqs

Under {\bf (H2)}, the gain functional is well-defined and finite, and we consider the value function of the partial observation control problem
\beq \label{defv}
\upv(t,\xi) & := & \sup_{\alpha\in\Ac^W} \Jc(t,\xi,\alpha), \;\;\;\;\;  t \in [0,T], \; \xi \in L^2(\Fc^o;\R^n),
\enq
which satisfies the quadratic growth estimate:
\beq \label{vquadra}
|\upv(t,\xi)| & \leq & C\big(1 + \bar\E|\xi|^2 \big),
\enq
for some positive constant $C$ depending on $C_1$ and $C_2$ in {\bf (H1)}-{\bf (H2)}.

This formulation of a partial observation control problem is general enough to include large classes of optimization models with latent unobservable factors, and classical partially observed control problems with observation process driven
by an additive noise of the controlled state signal process  can be recast  in our above form \reff{dynX}-\reff{defv}  by the usual  method of reference probability, which transforms  the observation process $W$ into  a Wiener process;
we refer to the introduction and Section 2.3  in our companion paper \cite{BCFP16b} for the details. As for \cite{BCFP16b}, we remark that the results of the present paper are still valid if we require  $b$ and $\sigma$ to be only locally Lipschitz and satisfy a linear growth condition in $x$, uniformly with respect to $a$. However, to alleviate the presentation, we assume  {\bf (H1)} throughout.

Optimal control of partially observed diffusions is a difficult problem, and a standard way of dealing with it is through
the so-called {\it separation} method. It consists in introducing the controlled filter, that is the conditional distribution defined for any test function
$\varphi$ $:$ $\R^n$ $\rightarrow$ $\R$ by
\beqs
\rho_s^{t,\xi,\alpha}[\varphi] &=& \bar\E\big[ \varphi(X_s^{t,\xi,\alpha}) | \Fc_s^W \big],
\enqs
%which is an infinite dimensional process,
and solution to the so-called controlled Zakai equation:
\beq \label{Zakaicontrol1}
d \rho_s^{t,\xi,\alpha} [\varphi] &=&  \rho_s^{t,\xi,\alpha}[L^{\alpha_s}\varphi]\,ds +
\rho_s^{t,\xi,\alpha} [M^{\alpha_s}\varphi]\,dW_s, \;\;\; t \leq s \leq T,
\enq
where $L^a\varphi$  $=$ $b(.,a)\cdot D_x\varphi  + \frac{1}{2}\text{tr}\big[\sigma\sigma\trans(.,a)D_x^2\varphi \big]$,
$M^a\varphi$ $=$ $D_x\varphi \trans \sigma_{_W} (.,a)$. The reward functional to be maximized is then rewritten from Bayes formula as
\beq \label{Jsepar}
\Jc(t,\xi,\alpha) & = & \bar\E \Big[ \int_t^T \rho_s^{t,\xi,\alpha} [ f(.,\alpha_s)] ds + \rho_T^{t,\xi,\alpha}[g]  \Big],
\enq
and the separated problem is a fully observable optimal control problem where one has to control the infinite dimensional filter process valued in the space of probability measures. This makes the
problem rather challenging, and the subject of intensive study. We give a short sketch of the possible approaches and refer the reader to the treatises
\cite{par88} and \cite{Be} for more complete results and references.  A first  method is the stochastic Pontryagin maximum principle, which provides necessary conditions for the optimality of a control process in terms of an adjoint equation
and can be used to solve successfully the problem in a number of cases: see \cite{Be} and \cite{Tang98}.
A second approach is the dynamic programming method applied to the separated control problem  \eqref{Zakaicontrol1}-\eqref{Jsepar},
often under the condition that  $\rho_s^{t,\xi,\alpha}(dx)$ admits a density $p_s^{t,\xi,\alpha}(x)$ with respect to the Lebesgue measure on $\R^n$,
and the controlled Zakai equation is then written as an equation for $p_s^{t,\xi,\alpha}$,
considered as a process with values in some Hilbert space (for instance the space $L^2(\R^n)$  or some weighted $L^2$-space),
the so-called controlled Duncan-Mortensen-Zakai (DMZ) stochastic partial  differential equation
\beq \label{DuncanMZcontrolintro}
dp_s^{t,\xi,\alpha} &=&  (L^{\alpha_s})^* p_s^{t,\xi,\alpha} \; ds  \; + \;  (M^{\alpha_s})^* p_s^{t,\xi,\alpha} \; dW_s, \;\;\; t \leq s \leq T,
%p_s^{t,\eta}  & = &  \eta +  \int_t^s(\call^{I_r})^*p_r^{t,\eta} \,dr + \int_t^s(\Mc^{I_r})^*p_r\,dW_r.
\enq
where $(L^a)^* p$ $=$ $- {\rm div}( b(.,a) p ) + \frac{1}{2} {\rm tr}\big[ D_x^2(\sigma\sigma\trans(.,a) p \big) \big]$,
$(M^a)^* p$  $=$ $- {\rm div}\big( \sigma_{_W}(.,a) p \big)$, are the adjoint operators of $L^a$ and $M^a$ (here ${\rm div}(.)$ denotes the divergence operator). This leads to a fully nonlinear  Bellman (HJB) equation with unbounded terms in the infinite dimensional space of density functions,  and has been studied by viscosity methods in \cite{lio88}, later significantly extended in \cite{gozzi_swiech}. However, existence of a density process solution to the DMZ equation \eqref{DuncanMZcontrolintro} is a difficult issue and needs some stringent conditions, see \cite{Be}, \cite{Rozovskii}, and Chapter 7 of \cite{BaCri}
for the uncontrolled case.  Moreover, a rigorous proof of the dynamic programming principle  for the value function to the separated control problem seems not available in the literature,
even though this is  the key step for proving the viscosity property to the associated HJB equation. Let us mention the so-called Nisio semigroup approach, strictly related to the dynamic programming principle,
and used for instance in  \cite{Be} for partial observation control problem but with drift control only, and which allows to characterize the value function via the non-linear semigroup as a mild solution to the dynamic programming equation.

We propose in this paper an alternative approach for tackling partial observation control problem, which allows us to relax the density and non-degeneracy assumptions mentioned above, and provides
an original partial differential equation (PDE) characterization of the value function.  We use a {\it randomization method} introduced in \cite{KP12} for classical Markovian models, and extended to the case of partial observation in \cite{BCFP16b}. We notice that the randomization method has already been successfully implemented for the study of Markovian models 
in \cite{KMPZ10}, \cite{bou09}, \cite{EKa}, \cite{ChoukrounCosso16}, \cite{CossoFuhrmanPham16}, \cite{CPH15}, \cite{Bandini15}, \cite{BandiniFuhrman15}, but in all these papers
only the full observation case was addressed and
the associated HJB equation was defined in a finite dimensional Euclidean space
(or it was a merely integral equation related to a controlled pure jump process),
whereas in the present paper we will be dealing with a HJB equation
on a genuinely infinite-dimensional space.
The basic idea is to consider a stepwise process $I$ associated to a Poisson random measure $\mu$, and replace the control process $\alpha$ by this exogenous process $I$ in the  dynamics of the process $X$, thus arriving at regime switching diffusion, called {\it randomized forward} process. We next consider an auxiliary optimization problem, called randomized or {\it dual} problem (in contrast with the primal problem defined in \reff{defv}),  which consists in optimizing among equivalent
changes of probability measures which only affect the intensity measure of $\mu$ but not the law of $W$. In the randomized problem, an admissible control is a bounded  positive   map  $\nu$ defined on $\Omega\times\R_+\times A$,
 which is predictable with respect to the filtration $\F^{W,\mu}$ generated  by $W$ and $\mu$. It turns out that the dual and primal problems coincide, and our first main result is to prove the so-called randomized dynamic programming principle (DPP)
 for the dual (hence for the primal) problem.  For this purpose, the key step is to show the flow property of the conditional distribution of the randomized forward process given the randomized observation filtration generated by $W$ and $\mu$.
Next, for exploiting the DPP, we use a notion of differentiability with respect to probability measures introduced by P.L. Lions in his lectures at the Coll\`ege de France \cite{lio12}, and detailed in the notes  \cite{car12}.
This notion of derivative is based on the lifting of functions defined on the Hilbert space of square integrable random variables distributed according to the ``lifted" probability measure.
By combining with a special It\^o's chain rule for flows of  conditional distributions in the spirit of \cite{cardel14b},  we derive the dynamic programming Hamilton-Jacobi-Bellman (HJB) equation  
for the  partial observation control problem.
This equation is a fully nonlinear second order partial differential equation (PDE)  in the infinite dimensional Wasserstein space  of probability measures,  and  in the case where the probability measures have a density, it reduces to the Bellman equation
in the space of density functions associated to the controlled DMZ 
stochastic PDE \reff{DuncanMZcontrolintro}, as derived in \cite{par88}, \cite{lio89}.
By adapting standard arguments to our context, we prove the viscosity property of the value function to the HJB  equation from the randomized dynamic programming principle. To complete our PDE characterization of the value function with a uniqueness result, it is convenient to work in the lifted Hilbert space of square integrable random variables instead of the Wasserstein metric space of probability measures, in order to rely on the general results for viscosity solutions  of second order Hamilton-Jacobi-Bellman equations in separable Hilbert spaces, see \cite{lio88}, \cite{lio89b}, \cite{gozswi15}.  We also state a verification theorem which is useful for getting an analytic feedback form of the optimal control when there is a smooth solution to the Bellman equation. Finally, we apply our results to a non Gaussian  linear-quadratic (LQ) partial observation control problem for which one can obtain explicit solutions.

The rest of the paper is organized as follows. Section \ref{preli} introduces notations and recalls some useful preliminaries about the Borel $\sigma$-algebra associated to the Wasserstein space of  probability measures.
We formulate in Section \ref{Randomized} the randomized partial observation control problem, and prove in Section \ref{secDPP} the randomized dynamic programming principle.
Section \ref{secHJB} is devoted to a verification theorem and the viscosity characterization of the value function to the HJB equation.
Finally,  Section \ref{secLQ} presents an application to a non Gaussian linear-quadratic  partial observation model with explicit solutions.

\section{Notations and preliminaries} \label{preli}

\setcounter{equation}{0} \setcounter{Assumption}{0}
\setcounter{Theorem}{0} \setcounter{Proposition}{0}
\setcounter{Corollary}{0} \setcounter{Lemma}{0}
\setcounter{Definition}{0} \setcounter{Remark}{0}

We introduce some notations used in the sequel of the paper.

\vspace{1mm}

\noindent {\bf Notations:} We denote by $x.y$ the scalar product of two Euclidian vectors $x$ and $y$,  and by $M\trans$ the transpose of a matrix or vector $M$.   We denote by  $\Pc_{_2}(\R^n)$ the space of probability measures $\pi$  on $\R^n$ with finite second order moment, namely
$\|\pi\|_{_2}^2$ $:=$ $\int_{\R^n} |x|^2 \pi(dx)$ $<$ $\infty$.  For any  $\pi$ $\in$ $\Pc_{_2}(\R^n)$,  we denote by $L_\pi^2(\R^q)$ the set of measurable functions $\varphi$ $:$ $\R^n$ $\rightarrow$ $\R^q$, which are square integrable with respect to $\pi$, by $L_{\pi\otimes\pi}^2(\R^q)$ the set of measurable functions $\psi$ $:$ $\R^n\times\R^n$ $\rightarrow$ $\R^q$,
which are square integrable with respect to the product measure  $\pi\otimes\pi$, and we set
\beqs
\pi[\phi] \; := \;   \int_{\R^n}\varphi(x)\,\pi(dx), & & \pi\otimes\pi[\psi] \;  := \;  \int_{\R^n\times\R^n} \psi(x,x') \pi(dx)\pi(dx').
\enqs
We also define  $L_\pi^\infty(\R^q)$ (resp.  $L_{\pi\otimes\pi}^\infty(\R^q)$) as the subset of elements $\varphi$ in $L_\pi^2(\R^q)$
(resp. $L_{\pi\otimes\pi}^2(\R^q)$), which are bounded $\pi$ (resp. $\pi\otimes\pi$) a.e., and $\|\varphi\|_\infty$ is their essential supremum.
For every random variable $\xi\colon(\Omega^o,\Fc^o)\rightarrow\R^n$, we denote by $\Lc(\xi)$ its probability law (also called distribution) under $\P^o$ (or equivalently $\bar\P$).

\begin{Remark}\label{R:OmegaPolish}
{\rm
We notice that  $\Fc^o$  is rich enough to carry $\R^n$-valued random variables with any arbitrary square-integrable distribution, namely $\Pc_{_2}(\R^n)$ $=$ $\{\Lc(\xi)\colon \xi \in L^2(\Fc^o;\R^n)\}$.
Indeed, when $\Omega^o=[0,1]$, $\Fc^o=\Bc([0,1])$ and $\P^o=\lambda$ the Lebesgue measure, then the result is well-known (see for instance Theorem 3.1.1 in \cite{Zabczyk96}). In the general case, by Corollary 7.16.1 in \cite{BertsekasShreve78} there exists a Borel isomorphism $\varphi\colon\Omega^o\rightarrow[0,1]$, which allows to define an atomless probability measure
$\P^o_{_\varphi}=\P^o\circ\varphi^{-1}$ on $([0,1],\Bc([0,1]))$. Let $F_{_\varphi}\colon[0,1]\rightarrow[0,1]$ be the cumulative distribution function corresponding to $\P^o_{_\varphi}$, namely $F_{_\varphi}(t)=\P^o_{_\varphi}([0,t])$, for all $t\in[0,1]$. Since $F_{_\varphi}$ is continuous, the function $G_{_\varphi}(s)=\inf\{t\in[0,1]\colon F_{_\varphi}(t)>s\}$, $s\in[0,1]$, satisfies $F_{_\varphi}(G_{_\varphi}(t))=G_{_\varphi}(F_{_\varphi}(t))=t$ (see for instance Lemma 1.37 in \cite{HeWangYan}). Notice that $\lambda=\P^o_{_\varphi}\circ G_{_\varphi}$ is the Lebesgue measure. Now, let $\pi\in\Pc_{_2}(\R^n)$ and consider a random variable $\xi\colon[0,1]\rightarrow\R^n$ with $\pi=\Lc(\xi)$.
Then, the random variable $\tilde \xi\colon\Omega^o\rightarrow\R^n$ given by $\tilde \xi=\xi\circ F_{_\varphi}\circ\varphi$ has distribution $\pi$. The claim follows from the arbitrariness of $\pi$.
\ep
}
\end{Remark}

Given a function $u$ defined on $\Pc_{_2}(\R^n)$, we call {\it lifted function} of $u$, the function $\upsilon$ defined on
$L^2(\Fc^o;\R^n)$ by $\upsilon(\xi)$ $=$ $u(\Lc(\xi))$.  Conversely, given a function $\upsilon$ defined on $L^2(\Fc^o;\R^n)$, we call {\it inverse-lifted function} of $\upsilon$ the function $u$ defined on $\Pc_{_2}(\R^n)$ by $u(\pi)$ $=$ $\upsilon(\xi)$ for $\pi$ $=$ $\Lc(\xi)$. Notice that $u$ exists if and only if $\upsilon(\xi)$ depends only on the distribution of $\xi$, for every $\xi\in L^2(\Fc^o;\R^n)$, and in this case there is a one-to-one correspondence between functions defined on $\Pc_{_2}(\R^n)$ and functions defined on   $L^2(\Fc^o;\R^n)$ via this lifting relation.
%and we shall sometimes identify $u$ and  $\upsilon$ when there is no ambiguity.
Recall that  $\Pc_{_2}(\R^n)$ is a metric space when endowed with the $2$-Wasserstein distance
$\Wc_2$ defined by
\begin{align*}
\Wc_2(\pi,\pi') \ &:= \ \inf\Big\{ \Big( \int_{\R^n\times\R^n} |x-x'|^2 \gamma(dx,dx') \Big)^{1\over 2}: \gamma \in \Pc_{_2}(\R^n\times\R^n) \mbox{ with marginals }  \pi \mbox{ and } \pi' \Big\}   \\
&= \ \inf\Big\{ \Big( \E|\xi-\xi'|^2 \Big)^{1\over 2}: \xi, \xi' \in  L^2(\Fc^o;\R^n) \mbox{ with }  \Lc(\xi) \; = \pi, \; \Lc(\xi') \; = \; \pi' \Big\}.
\end{align*}

We end this section providing two useful properties related to the space $\Pc_{_2}(\R^n)$, endowed with topology (and the Borel $\sigma$-algebra $\Bc(\Pc_{_2}(\R^n))$) induced by the Wasserstein metric $\Wc_2$. In the sequel we denote by
$\mathscr B_{_2}(\R^n)$ the set of Borel-measurable functions on $\R^n$ with quadratic growth.

\begin{Lemma}[Skorohod's representation theorem for $\Wc_2$-convergence]\label{L:W2-convergence}
Let $(\pi_m)_m$ be a sequence in $\Pc_{_2}(\R^n)$ such that $\Wc_2(\pi_m,\pi)\rightarrow0$, for some $\pi\in\Pc_{_2}(\R^n)$. Then, there exists a sequence of random variables $(\xi_m)_m\subset L^2(\Fc^o;\R^n)$, with $\Lc(\xi_m)=\pi_m$, converging $\bar\P$-a.s. and in $L^2(\Fc^o;\R^n)$ to some
$\xi\in L^2(\Fc^o;\R^n)$, with $\Lc(\xi)=\pi$.
\end{Lemma}
\textbf{Proof.}
By Theorem 6.9 in \cite{Villani}, we have that $\Wc_2(\pi_m,\pi)\rightarrow0$ is equivalent to the following convergence:
\[
\int_{\R^n} \varphi(x)\,\pi_m(dx) \ \overset{n\rightarrow\infty}{\longrightarrow} \ \int_{\R^n} \varphi(x)\,\pi(dx), \qquad \text{for every }\varphi\in\mathscr B_{_2}(\R^n).
\]
In particular, $\pi_m$ converges weakly to $\pi$, namely $\int \varphi\,d\pi_m \rightarrow \int \varphi\,d\pi$ for any bounded continuous $\varphi$. Therefore, by Skorohod's representation theorem (see for instance Theorem 2.2.2 in \cite{Borkar}) and in particular Remark \ref{R:OmegaPolish}, there exist random variables $(\xi_m)_m$ and $\xi$ in $L^2(\Fc^o;\R^n)$, with $\Lc(\xi_m)=\pi_m$ and $\Lc(\xi)=\pi$, such that $\xi_m$ converges $\bar\P$-a.s. to $\xi$. It remains to prove the convergence in $L^2(\Fc^o;\R^n)$. To this end, we notice that, by Proposition 7.1.5 in \cite{AmbrosioGS} the sequence $(|\xi_m|^2)_m$ is uniformly integrable. Then, it follows that $\xi_m\rightarrow\xi$ in $L^2(\Fc^o;\R^n)$.
\ep

\begin{Lemma}\label{L:P2-Meas}
Let $(E,\Ec)$ be a measurable space and consider a map $\Pi\colon E\rightarrow\Pc_{_2}(\R^n)$. For every $\varphi\in\mathscr B_{_2}(\R^n)$, we denote by $\Pi[\varphi]\colon E\rightarrow\R$ the following map:
\[
\Pi[\varphi](e) \ = \ \int_{\R^n}\varphi(x)\,\Pi(e,dx), \qquad \text{for every }e\in E.
\]
Then, $\Pi$ is measurable if and only if, for every continuous $\varphi\in\mathscr B_{_2}(\R^n)$, $\Pi[\varphi]$ is measurable.
\end{Lemma}
\textbf{Proof.}
For every continuous $\varphi\in\mathscr B_{_2}(\R^n)$, define the map $\Lambda_{_\varphi}\colon\Pc_{_2}(\R^n)\rightarrow\R$ as follows:
\[
\Lambda_{_\varphi}\pi \ = \ \int_{\R^n}\varphi(x)\,\pi(dx), \qquad \text{for every }\pi\in\Pc_{_2}(\R^n).
\]
By Theorem 6.9 in \cite{Villani}, we have that $\Wc_2(\pi_m,\pi)\rightarrow0$ if and only if, for every $\varphi\in\mathscr B_{_2}(\R^n)$, $\Lambda_{_\varphi}\pi_m\rightarrow\Lambda_{_\varphi}\pi$. Therefore, the $\sigma$-algebra $\sigma(\Lambda_{_\varphi},\,\varphi\in\mathscr B_{_2}(\R^n)\text{ continuous})$ coincides with $\Bc(\Pc_{_2}(\R^n))$.

Notice that, for every $\varphi\in\mathscr B_{_2}(\R^n)$, the map $\Pi[\varphi]$ is given by the composition of $\Lambda_{_\varphi}$ and $\Pi$, namely:
\[
\Pi[\varphi](e) \ = \ (\Lambda_{_\varphi}\circ\Pi)(e), \qquad \text{for every }e\in E.
\]
Then, recalling that $\sigma(\Lambda_{_\varphi},\,\varphi\in\mathscr B_{_2}(\R^n)\text{ continuous})$ is equal to $\Bc(\Pc_{_2}(\R^n))$, we deduce the equivalence between the measurability of $\Pi$ and the measurability of $\Lambda_{_\varphi}\circ\Pi$, for every continuous $\varphi\in\mathscr B_{_2}(\R^n)$.
\ep

\section{Randomized partial observation control problem} \label{Randomized}

\setcounter{equation}{0} \setcounter{Assumption}{0}
\setcounter{Theorem}{0} \setcounter{Proposition}{0}
\setcounter{Corollary}{0} \setcounter{Lemma}{0}
\setcounter{Definition}{0} \setcounter{Remark}{0}

%\subsection{Randomized formulation}

Up to a product extension, we can assume without loss of generality that the complete probability space $(\Omega,\Fc,\P)$ supports in addition to the Brownian motion $W$,
a Poisson random measure $\mu$ on $\R_+\times A$, independent of $W$.  The random measure $\mu$ has compensator $\lambda(da)dt$, for some finite positive measure
$\lambda$ on $(A,\Bc(A))$, with full topological support, and thus  $\mu$ is a sum  of Dirac measures of the form $\mu=\sum_{n\ge 1}\delta_{(T_n,A_n)}$,   where $(A_n)_{n\ge 1}$ is a sequence of $A$-valued random variables
 and  $(T_n)_{n\ge 1}$ is a strictly increasing sequence of random variables with values in $(0,\infty)$, and for any $C\in\calb(A)$ the process $\mu((0,t]\times C)-t\lambda(C)$, $t\ge 0$, is a $\P$-martingale.

For every $(t,a)\in[0,T]\times A$, we associate to  the random measure $\mu$ the $A$-valued pure jump process
\beq\label{defI}
I_s^{t,a} & =&   \sum_{n\ge 1} a\,1_{[t,T_n)}(s) + \sum_{\substack{n\ge 1 \\ t<T_n}} A_n\,1_{[T_n,T_{n+1})}(s), \qquad s\ge t.
\enq
Notice that, when $A$ is a subset of the linear space, formula \reff{defI} can be written as
\beqs
I_s^{t,a} & = & a + \int_t^s\int_A(a'-I_{r-}^{t,a})\,\mu(dr\, da'), \qquad s\ge t.
\enqs
From the definition of $I^{t,a}$ in \eqref{defI}, we see that the flow property holds:
\beq\label{flowI}
I_s^{t,a} &=& I_s^{\theta,I_\theta^{t,a}},
\enq
for all $\theta,s\in[t,T]$, with $\theta\leq s$. The probability space $(\bar\Omega,\bar\Fc,\bar\P)$  on which we shall define the randomized partial observation control problem, can be described graphically as follows:
%is defined is given by $(\Omega^\Rc,\Fc^\Rc,\P^\Rc)$, where $\Omega^\Rc=\Omega^o\times\Omega^1\times\Omega^2\times\Omega^3$, $\Fc^\Rc$ is the completion of $\Fc^o\otimes\Fc^1\otimes\Fc^2\otimes\Fc^3$ with respect to $\P^o\otimes%\P^1\otimes\P^2\otimes\P^3$, and $\P^\Rc$ is the extension of  $\P^o\otimes\P^1\otimes\P^2\otimes\P^3$ to $\Fc^\Rc$.
%We also denote by $\Fc^{o,\Rc}=\{E\times\Omega^1\times\Omega^2\times\Omega^3\colon E\in\Fc^o\}$
%the canonical extension of $\Fc^o$ to $\Omega^\Rc$.
%Recalling that $(\Omega^o,\Fc^o,\P^o)$ supports the initial condition $\xi$ of the state equation, $(\Omega^1,\Fc^1,\P^1)$ supports $V$, $(\Omega^2,\Fc^2,\P^2)$ supports $W$, and since $(\Omega^3,\Fc^3,\P^3)$ supports $\mu$, we can %describe these relations graphically as follows:
\[
\bar\Omega \ = \ \overbrace{\underset{\!\!\!\underset{\xi}{\uparrow}}{\Omega^o} \times \underset{\!\!\!\underset{V}{\uparrow}}{\Omega^1}}^{\bar\Omega^1 \; :=} \times \underset{\underset{W,\mu}{\uparrow}}{\Omega},
%\times \underset{\!\!\!\underset{\mu}{\uparrow}}{\Omega^3}}^{\Omega\,:=},
\]
and we extend canonically $V,W,\mu,\xi$ on $(\bar\Omega,\bar\Fc,\bar\P)$.
As reported in the graph, we have defined the space  $(\bar\Omega^1,\bar\Fc^1,\bar\P^1)$, where $\bar\Omega^1=\Omega^o\times\Omega^1$, $\bar\Fc^1$ is the completion of $\Fc^o\otimes\Fc^1$ with respect to
$\P^o\otimes\P^1$,  $\bar\P^1$ is the extension of $\P^o\otimes\P^1$ to $\bar\Fc^1$.
%and  $(\Omega,\Fc,\P)$, where $\Omega=\Omega^2\times\Omega^3$, $\Fc$ is the completion of $\Fc^2\otimes\Fc^3$ with respect to $\P^2\otimes\P^3$, and
%$\P$ is the extension of $\P^2\otimes\P^3$ to $\calf$.

For every $t\geq0$, we denote by $\F^{B,\mu,t}=(\Fc_s^{B,\mu,t})_{s\geq t}$ the $\P^1\otimes\P$-completion of the filtration on $\Omega^1\times\Omega$ generated by $(B_s-B_t)_{s\geq t}$ and $\mu\,1_{[t,\infty)}$. We then consider, for every $(t,x,a)\in[0,T]\times\R^n\times A$, the following regime switching diffusion equation:
\beq \label{dynXI_0}
dX_s &=& b(X_s,I_s^{t,a}) ds + \sigma(X_s,I_s^{t,a}) dB_s,  \;\;\;  t \leq s \leq T, \qquad  X_t = x.
\enq
We also require that $X_s=x$, for every $s\in[0,t)$.

\begin{Lemma}\label{L:MeasVersion}
Under Assumption {\bf (H1)}, there exists a random field $(X_s^{t,x,a})$, with $(t,s,x,a)\in[0,T]\times[0,T]\times\R^n\times A$, such that:
\begin{enumerate}
\item[\textup{(i)}] for every $(t,x,a)\in[0,T]\times\R^n\times A$, $(X_s^{t,x,a})_{s\in[t,T]}$ is the unique $($up to indistinguishability$)$ $\F^{B,\mu,t}$-adapted continuous process solution to \reff{dynXI_0}, and $X_s^{t,x,a}=x$, for every $s\in[0,t)$;
\item[\textup{(ii)}] the map $X_s^{t,x,a}\colon([0,T]\times[0,T]\times\Omega^1\times\Omega\times\R^n\times A,\Bc([0,T])\otimes\Bc([0,T])\otimes\Fc^1\otimes\Fc\otimes\Bc(\R^n)\otimes\Bc(A))\rightarrow(\R^n,\Bc(\R^n))$ is measurable.
\end{enumerate}
\end{Lemma}
\textbf{Proof.}
In order to construct the random field $(X_s^{t,x,a})$, it is enough to show that there exists a random field $(\mathscr X_s^{t,x,a})$, with $(t,s,x,a)\in[0,T]\times[0,T]\times\R^n\times A$, such that:
\begin{enumerate}
\item[\textup{(i)'}] for every $(t,x,a)\in[0,T]\times\R^n\times A$, $(\mathscr X_s^{t,x,a})_{s\in[t,T]}$ is the unique $($up to indistinguishability$)$ $\F^{B,\mu,t}$-adapted continuous process solution to \reff{dynXI_0}, and $\mathscr X_s^{t,x,a}=x$, for every $s\in[0,t)$;
\item[\textup{(ii)'}] the map $\mathscr X_s^{t,x,a}\colon([0,T]\times[0,T]\times\Omega^1\times\Omega\times\R^n\times A,\Bc([0,T])\otimes\Bc([0,T])\otimes\Fc_\infty^{B,\mu,t}\otimes\Bc(\R^n)\otimes\Bc(A))\rightarrow(\R^n,\Bc(\R^n))$ is measurable.
\end{enumerate}
As a matter of fact, suppose that we have already constructed the random field $(\mathscr X_s^{t,x,a})$. Then, by a monotone class argument, we see that there exists a $\P^1\otimes\P$-null set $N\in\Fc^1\otimes\Fc$ and a random field $(X_s^{t,x,a})$, with $(t,s,x,a)\in[0,T]\times[0,T]\times\R^n\times A$, which is measurable as a map from $([0,T]\times[0,T]\times\Omega^1\times\Omega\times\R^n\times A,\Bc([0,T])\otimes\Bc([0,T])\otimes\Fc^1\otimes\Fc\otimes\Bc(\R^n)\otimes\Bc(A))$ into $(\R^n,\Bc(\R^n))$, such that $(\mathscr X_s^{t,x,a})$ coincides with $(X_s^{t,x,a})$ on $[0,T]\times[0,T]\times((\Omega^1\times\Omega)\backslash N)\times\R^n\times A$. It is then easy to see that $(X_s^{t,x,a})$ satisfies items (i) and (ii) of the Lemma.

It remains to construct a random field $(\mathscr X_s^{t,x,a})$ satisfying (i)' and (ii)'. To this end, we exploit the Picard iterations method. More precisely, we define recursively a sequence of processes $(\mathscr X^{m,t,x,a})_m$, setting $\mathscr X^{0,t,x,a}\equiv0$ and then defining $\mathscr X^{m+1,t,x,a}$ from $\mathscr X^{m,t,x,a}$ as follows:
\beqs
d\mathscr X_s^{m+1,t,x,a} &=& b(\mathscr X_s^{m,t,x,a},I_s^{t,a}) ds + \sigma(\mathscr X_s^{m,t,x,a},I_s^{t,a}) dB_s,  \;\;\;  t \leq s \leq T, \;\;\;  \mathscr X_t^{m+1,t,x,a} = x.
\enqs
We also set $\mathscr X_s^{m+1,t,x,a}=x$, for every $s\in[0,t)$. It can be shown (proceeding for instance as in Lemma 14.20 of \cite{jacod_book}) that there exists a random field $(\mathscr X_s^{t,x,a})$ such that we have the following convergence in probability:
\begin{equation}\label{mathscrXm-->X}
\sup_{s\in[0,T]}|\mathscr X_s^{m,t,x,a}-\mathscr X_s^{t,x,a}|\ \overset{\P^1\otimes\P}{\longrightarrow} \ 0,  \qquad \text{as $m$ tends to infinity}.
\end{equation}
Then, for every $(t,x,a)\in[0,T]\times\R^n\times A$, the process $(\mathscr X_s^{t,x,a})_{s\in[t,T]}$ is $\F^{B,\mu,t}$-adapted and continuous, and solves equation \reff{dynXI_0}. Moreover, $\mathscr X_s^{t,x,a}=x$, for every $s\in[0,t)$. Therefore, by Theorem 14.23 in \cite{jacod_book}, $(\mathscr X_s^{t,x,a})_{s\in[t,T]}$ is the unique (up to indistinguishability) $\F^{B,\mu,t}$-adapted continuous process solution to \reff{dynXI_0}. This proves item (i)'.

On the other hand, we notice that, for every $m$, the map $\mathscr X_s^{m,t,x,a}\colon([0,T]\times[0,T]\times\Omega^1\times\Omega\times\R^n\times A,\Bc([0,T])\otimes\Bc([0,T])\otimes\Fc_\infty^{B,\mu,t}\otimes\Bc(\R^n)\otimes\Bc(A))\rightarrow(\R^n,\Bc(\R^n))$ is measurable. Then, by the convergence in probability \eqref{mathscrXm-->X}, we see (proceeding for instance as in the first point of Exercise IV.5.17 in \cite{revyor99}) that the map $\mathscr X_s^{t,x,a}\colon([0,T]\times[0,T]\times\Omega^1\times\Omega\times\R^n\times A,\Bc([0,T])\otimes\Bc([0,T])\otimes\Fc_\infty^{B,\mu,t}\otimes\Bc(\R^n)\otimes\Bc(A))\rightarrow(\R^n,\Bc(\R^n))$ can be taken to be measurable. This implies that item (ii)' holds true.
\ep

\vspace{3mm}

For every $t\in[0,T]$, let us denote by $\bar\F^{o, B, \mu,t}=(\bar\Fc_s^{0, B, \mu,t})_{s\in[t,T]}$ the $\bar\P$-completion of the filtration generated by
$\Fc^{o}$, $(B_s-B_t)_{s\geq t}$, and $\mu\,1_{(t,\infty)}$, and by $L^2(\bar\Fc_t^{o, B, \mu};\R^n)$ the set of square integrable
$\bar\Fc_t^{o,B,\mu}$ ($:=$ $\bar\Fc_t^{o,B,\mu,0}$)-measurable random variables on $(\bar\Omega,\bar\Fc,\bar\P)$.
By Lemma \ref{L:MeasVersion}, the pathwise uniqueness of a solution to \reff{dynXI}, and identity \eqref{flowI}, we deduce the following standard result.

\begin{Lemma}\label{L:xi}
Under Assumption {\bf (H1)}, for every $(t,\xi,a)\in[0,T]\times L^2(\bar\Fc_t^{o, B, \mu};\R^n)\times A$, the process $(X_s^{t,\xi,a})_{s\in[t,T]}$, defined as $X_s^{t,\xi,a}(\omega^o,\omega^1,\omega)=X_s^{t,\xi(\omega^o),a}(\omega^1,\omega)$, is the unique (up to indistinguishability) $\bar\F^{o, B, \mu,t}$-adapted continuous solution to the equation:
\beq\label{dynXI}
dX_s^{t,\xi,a} &=& b(X_s^{t,\xi,a},I_s^{t,a}) ds + \sigma(X_s^{t,\xi,a},I_s^{t,a}) dB_s,  \;\;\;  t \leq s \leq T, \;\;\;  X_t^{t,\xi,a} = \xi,
\enq
and satisfies the standard estimate
\beq \label{estimXrandom}
\bar\E \Big[ \sup_{t\leq s\leq T} \big|X_s^{t,\xi,a}\big|^2\Big] & \leq & C \big( 1 + \bar\E|\xi|^2 \big) \; < \; \infty,
\enq
for some positive constant $C$. In addition, we have the flow property:
\beq\label{flowXI}
X_s^{t,\xi,a} &=& X_s^{\theta,X_\theta^{t,\xi,a},I_\theta^{t,a}},
\enq
$\P^\Rc$-a.s., for all $\theta,s\in[t,T]$, with $\theta\leq s$.
\end{Lemma}

\begin{Remark}
{\rm  By misuse of notation, we have written here $X^{t,\xi,a}$ for the solution to the randomized SDE \reff{dynXI}, but we point out that this differs from the solution $X^{t,\xi,\alpha}$ to the controlled SDE \reff{dynX} when
$\alpha$ is constant equal to $a$.
}
\ep
\end{Remark}

We next consider on $(\Omega,\Fc,\P)$, for every $t\geq0$, the so-called {\it randomized observation filtration} $\F^{W,\mu,t}=(\Fc_s^{W,\mu,t})_{s\geq t}$, as the $\P$-completion of the filtration generated by $(W_s-W_t)_{s\geq t}$ and $\mu\,1_{[t,\infty)}$, and denote by $\calp(\F^{W,\mu,t})$ the corresponding predictable
$\sigma$-algebra on $\Omega\times[t,\infty)$.

We can then define, for every $t\in[0,T]$, the randomized optimal control problem as follows: the set $\Vc_{_{W,\mu,t}}$  of admissible controls  consists of all $\nu=\nu_s(\omega,a): \Omega\times[t,\infty)\times A\to (0,\infty)$,
which are $\calp(\F^{W,\mu,t})\otimes \calb(A)$-measurable and bounded.

\begin{Remark}\label{R:V_W,mu,t}
{\rm The randomized problem can be equivalently formulated (in the sense that the value function is the same) with $\bar\Vc_{_{W,\mu,t}}$ instead of $\Vc_{_{W,\mu,t}}$. More precisely, define on $(\bar\Omega,\bar\Fc,\bar\P)$,
for every $t\in[0,T]$, the filtration $\bar\F^{W,\mu,t}=(\bar\Fc_s^{W,\mu,t})_{s\geq t}$ as the $\bar\P$-completion of the filtration generated by $(W_s-W_t)_{s\geq t}$ and $\mu\,1_{[t,\infty)}$,
and let $\calp(\bar\F^{W,\mu,t})$ denote the corresponding predictable $\sigma$-algebra on $\bar\Omega\times[t,\infty)$.
Consider the set $\bar\Vc_{_{W,\mu,t}}$ of maps $\bar\nu=\bar\nu_s(\bar\omega,a): \bar\Omega\times [t,\infty)\times A\to (0,\infty)$,   which are $\calp(\bar\F^{W,\mu,t})\otimes \calb(A)$-measurable and bounded. Then, the value function of the randomized optimal control problem (defined below in \eqref{defvrandom}) remains the same if we replace $\Vc_{_{W,\mu,t}}$ by $\bar\Vc_{_{W,\mu,t}}$.
Indeed, for every $s\geq t$, $\bar\Fc_s^{W,\mu,t}$ is the $\bar\P$-completion of the canonical extension of $\Fc_s^{W,\mu,t}$ to $\bar\Omega$. As a consequence, using for instance Proposition 1.1 in \cite{jacod_book}, we see that given $\bar\nu\in\bar\Vc_{_{W,\mu,t}}$ there exist $\nu\in\Vc_{_{W,\mu,t}}$ and a $\bar\P$-null set $\bar N\in\bar\Fc$ such that $\bar\nu_s(\bar\omega,a)=\nu_s(\omega,a)$ for every $(\bar\omega,s,a)\in(\bar\Omega\backslash \bar N)\times[t,\infty)\times A$, where $\bar\omega=(\omega^o,\omega^1,\omega)$. Notice that in \cite{BCFP16b} the randomized problem was formulated using the class $\bar\Vc_{_{W,\mu,t}}$.
\ep}
\end{Remark}

For every $t\in[0,T]$, $\nu$ $\in$ $\Vc_{_{W,\mu,t}}$, consider the Dol\'eans exponential process
\beqs \label{densitykappa}
\kappa_s^{\nu,t} \ &=& \ \Ec_s\bigg(\int_t^\cdot\int_A (\nu_r(a) - 1)\,  \big(\mu(dr,da) - \lambda(da)dr\big) \bigg).
%\\ &=& \ \exp\left(\int_0^t\int_A (1 - \nu_s(a))\lambda(da)\,ds \right)\prod_{T_n\le t}\nu_{T_n}(A_n),\qquad t\ge 0,\label{doleans}
\enqs
We then define a new probability measure on $(\bar\Omega,\bar\Fc)$ as $d\bar\P^{\nu,t}$ $=$ $\kappa_T^{\nu,t}\,d\bar\P$. Since $\kappa^{\nu,t}$ is defined on $\Omega$, we also define a new probability measure on $(\Omega,\Fc)$ as $d\P^{\nu,t}$ $=$ $\kappa_T^{\nu,t}\,d\P$. From
Girsanov's theorem it follows that under $\P^{\nu,t}$ the compensator of $\mu$ on the set
$[t,T]\times A$ is the random measure $\nu_s^t(a)\lambda(da)ds$, while $B$ remains a Brownian motion under $\P^{\nu,t}$, and the law of a
$\Fc^{o}$-measurable random variable $\xi$ is the same under  $\bar\P$ and $\bar\P^{\nu,t}$. Therefore, under {\bf (H1)}, we obtain by standard arguments the following estimate:
\beq \label{EstimateX_nu}
    \sup_{\nu\in\Vc_{_{W,\mu,t}}}\,\bar\E^{\nu,t}\,\Big[\sup_{s\in [t,T]}|X_s^{t,\xi,a}|^2 \Big]  & \leq  &  C \big( 1 + \bar\E\big[|\xi|^2\big] \big) \; < \; \infty,
\enq
where $\bar\E^{\nu,t}$ (resp. $\E^{\nu,t}$) denotes the expectation with respect to $\bar\P^{\nu,t}$ (resp. $\P^{\nu,t}$). We finally introduce the gain functional of the randomized control problem
\beq \label{defJrandomized}
\Jc^\Rc(t,\xi,\nu,a) & := &  \bar\E^{\nu,t}\Big[ \int_t^Tf(X_s^{t,\xi,a},I_s^{t,a}) ds + g(X_T^{t,\xi,a}) \Big]
\enq
and the corresponding value function:
\beq \label{defvrandom}
\upv^\Rc(t,\xi,a) &:=&   \sup_{\nu\in\Vc_{_{W,\mu,t}}} \Jc^\Rc(t,\xi,\nu,a), \;\;\;\;\;  t \in [0,T], \; \xi \in L^2(\Fc^{o};\R^n),\;a\in A.
\enq

It is proved in \cite{BCFP16b} that the value functions of the primal partial observation control problem \reff{defv} and of the randomized control problem \reff{defvrandom} coincide:
\beq \label{relvvrandom}
\upv(t,\xi) &=& \upv^\Rc(t,\xi,a), \;\;\;\;\;  t \in [0,T], \; \xi \in L^2(\Fc^{o};\R^n),\;a\in A.
\enq
This key relation will be exploited in the sequel for deriving a so-called randomized dynamic programming principle and then a Bellman equation characterizing (in the viscosity sense) the value function to the partial observation control problem.

\vspace{2mm}

\begin{Remark} \label{remI}
{\rm  Relation \reff{relvvrandom} shows that the value function of the randomized control problem depends only on the objects $B$ $=$ $(V,W)$, $b$, $\sigma$, $f$, $g$ and $A$ of the primal  partial observation control problem, and thus not on the intensity measure $\lambda$ of $\mu$, nor on the initial point $a$ taken by the jump process $I^{t,a}$.
}
\ep
\end{Remark}

 \section{Randomized filtering and dynamic programming principle} \label{secDPP}

\setcounter{equation}{0} \setcounter{Assumption}{0}
\setcounter{Theorem}{0} \setcounter{Proposition}{0}
\setcounter{Corollary}{0} \setcounter{Lemma}{0}
\setcounter{Definition}{0} \setcounter{Remark}{0}

 \subsection{Randomized filter process  and Zakai equation}

We first consider the  filtering of the randomized process solution to \reff{dynXI} with respect to the  randomized observation filtration.
%which consists in the process of conditional distributions $\rho_s^{t,\xi}$ of $X_s^{t,\xi}$ given ${\cal F}_s^{W,\mu}$ for $t$ $\in$ $[0,T]$,  $\xi$ $\in$ $L^2(\Fc_0;\R^n)$.
More precisely, for every $(t,s,\omega,\pi,a)\in[0,T]\times[t,T]\times\Omega\times\Pc_{_2}(\R^n)\times A$, we define $\rho_s^{t,\pi,a}(\omega)\in\Pc_{_2}(\R^n)$ as follows:
\begin{equation}\label{rho}
\rho_s^{t,\pi,a}(\omega)[\varphi] \ := \ \E^1\bigg[\int_{\R^n} \varphi(X_s^{t,x,a}(\cdot,\omega)) \,\pi(dx)\bigg],
\end{equation}
for every $\varphi$ $\in$ $\mathscr B_{_2}(\R^n)$. We also set $\rho_s^{t,\pi,a}(\omega)=\pi$, for every $(t,s,\omega,\pi,a)\in[0,T]\times[0,t)\times\Omega\times\Pc_{_2}(\R^n)\times A$.

\begin{Remark}\label{R:Filter}
{\rm
Let $\xi\in L^2(\Fc^{o};\R^n)$ with law $\pi$ and consider the process $(X_s^{t,\xi,a})_{s\in[t,T]}$ solution to \eqref{dynXI}. Then, we see that
\begin{align*}
\rho_s^{t,\pi,a}(\omega)[\varphi] \ := \ \E^1\bigg[\int_{\R^n} \varphi(X_s^{t,x,a}(\cdot,\omega)) \,\pi(dx)\bigg] \ &= \ \bar\E^1\big[\varphi(X_s^{t,\xi(\cdot),a}(\cdot,\omega))\big] \\
&= \ \bar\E\big[\varphi(X_s^{t,\xi,a})\big| \bar\Fc_s^{W,\mu,t}\big](\omega),
\end{align*}
$\P(d\omega)$-a.s., where the $\sigma$-algebra $\bar\Fc_s^{W,\mu,t}$ was introduced in Remark \ref{R:V_W,mu,t}. In other words,
$\rho_s^{t,\pi,a}(\omega)$ is the law of  the random variable $X_s^{t,\xi(.),a}(.,\omega)$ on $(\bar\Omega^1,\bar\Fc^1,\bar\P^1)$, or equivalently the
realization at $\omega$ of the conditional law of $X_s^{t,\xi,a}$ given $\bar\Fc_s^{W,\mu,t}$.
Notice that $X_s^{t,\xi,a}$ is independent of the increments of $W_r-W_s$, and $\mu([s,r)\times B)$ for all $r\geq s$, $B$ $\in$ $\Bc(A)$. Thus, $\rho_s^{t,\pi,a}$ is equal to the conditional law of
$X_s^{t,\xi,a}$ given $\bar\Fc_\infty^{W,\mu}$, denoted by $\Lc(X_s^{t,\xi,a}|W,\mu)$.
\ep}
\end{Remark}

\begin{Lemma}
Under Assumption {\bf (H1)}, the $\Pc_{_2}(\R^n)$-valued random field $(\rho_s^{t,\pi,a})$, with $t\in[0,T]$, $s\in[0,T]$, $\pi\in\Pc_{_2}(\R^n)$, $a\in A$, is such that:
\begin{enumerate}
\item[\textup{(i)}] the map $\rho_s^{t,\pi,a}\colon([0,T]\times[0,T]\times\Omega\times\Pc_{_2}(\R^n)\times A,\Bc([0,T])\otimes\Bc([0,T])\otimes\Fc\otimes\Bc(\Pc_{_2}(\R^n))\otimes\Bc(A))\rightarrow(\Pc_{_2}(\R^n),\Bc(\Pc_{_2}(\R^n)))$ is measurable;
\item[\textup{(ii)}] for every $(t,\pi,a)\in[0,T]\times\Pc_{_2}(\R^n)\times A$, the stochastic process $(\rho_s^{t,\pi,a})_{s\in[t,T]}$ is $\F^{W,\mu,t}$-predictable, and satisfies the estimate
\beq \label{estimrhopi}
\E\Big[ \sup_{t\leq s\leq T} |\rho_s^{t,\pi,a}|^2 \Big] & \leq & C(1 + \|\pi\|_{_2}^2),
\enq
for some positive constant $C$.
\end{enumerate}
\end{Lemma}
\textbf{Proof.}
Concerning item (i), we begin noting that it is enough to prove that for every $\varphi\in\mathscr B_{_2}(\R^n)$ the map $\rho_s^{t,\pi,a}(\omega)[\varphi]\colon([0,T]\times[0,T]\times\Omega\times\Pc_{_2}(\R^n)\times A,\Bc([0,T])\otimes\Bc([0,T])\otimes\Fc\otimes\Bc(\Pc_{_2}(\R^n))\otimes\Bc(A))\rightarrow(\R,\Bc(\R))$ is measurable. In order to prove this latter property, fix $\varphi\in\mathscr B_{_2}(\R^n)$ and notice that, by Lemma \ref{L:MeasVersion} and Fubini's theorem, the map $\E^1[\varphi(X_s^{t,x,a}(\cdot,\omega))]\colon([0,T]\times[0,T]\times\Omega\times\R^n\times A,\Bc([0,T])\otimes\Bc([0,T])\otimes\Fc\otimes\Bc(\R^n)\otimes\Bc(A))\rightarrow(\R,\Bc(\R))$ is measurable. Then, by a monotone class argument (proving firstly the measurability of $\rho_s^{t,\pi,a}(\omega)[\varphi]$ when the map $\E^1[\varphi(X_s^{t,x,a}(\cdot,\omega))]$ can be expressed as a product $\varphi_1(t,s,\omega,a)\varphi_2(x)$, for some measurable and bounded functions $\varphi_1$ and $\varphi_2$), we see that $\rho_s^{t,\pi,a}(\omega)[\varphi]$ is measurable, so that the claim follows.

Let us now prove property (ii). Fix $(t,\pi,a)\in[0,T]\times\Pc_{_2}(\R^n)\times A$. Notice that, by Lemma \ref{L:P2-Meas}, item (ii) follows if we prove that for every continuous $\varphi\in\mathscr B_{_2}(\R^n)$ the process $(\rho_s^{t,\pi,a}[\varphi])_{s\in[t,T]}$ is $\F^{W,\mu,t}$-predictable. In order to prove this latter property for a fixed continuous $\varphi\in\mathscr B_{_2}(\R^n)$, it is enough to show that:
\begin{enumerate}
\item[(ii.a)] $\P(d\omega)$-a.s., the map $s\mapsto\rho_s^{t,\pi,a}(\omega)[\varphi]$ is continuous;
\item[(ii.b)] the process $(\rho_s^{t,\pi,a}[\varphi])_{s\in[t,T]}$ is $\F^{W,\mu,t}$-adapted.
\end{enumerate}

Regarding item (ii.a), we begin noting that, since, for every $x\in\R^n$, $(\P^1\otimes\P)$-a.e. path of the process $(X_s^{t,x,a})_{s\in[t,T]}$ is continuous, we deduce that the map $s\mapsto \varphi(X_s^{t,x,a}(\omega^1,\omega))$ is $(\P^1\otimes\P)(d\omega^1,d\omega)$-a.s. continuous. This implies that
\begin{equation}\label{Continuity_rho0}
\int_{\Omega^1\times\Omega} \lim_{\delta\downarrow0}\sup_{\substack{|s' - s''|<\delta\\ s',s''\in\Q}}\big|\varphi(X_{s'}^{t,x,a}(\omega^1,\omega)) - \varphi(X_{s''}^{t,x,a}(\omega^1,\omega))\big|\,(\P^1\otimes\P)(d\omega^1,d\omega) \ = \ 0,
\end{equation}
where we observe that the map
\[
(\omega^1,\omega,x) \ \longmapsto \ \lim_{\delta\downarrow0}\sup_{\substack{|s' - s''|<\delta\\ s',s''\in\Q}}\big|\varphi(X_{s'}^{t,x,a}(\omega^1,\omega)) - \varphi(X_{s''}^{t,x,a}(\omega^1,\omega))\big|
\]
is $\Fc^1\otimes\Fc\otimes\Bc(\R^n)$-measurable, as a consequence of Lemma \ref{L:MeasVersion}(ii). Then, from \eqref{Continuity_rho0} and by Fubini's theorem, we obtain
\begin{align*}
0 \ &= \ \int_{\R^n} \bigg(\int_{\Omega^1\times\Omega} \lim_{\delta\downarrow0}\sup_{\substack{|s' - s''|<\delta\\ s',s''\in\Q}}\big|\varphi(X_{s'}^{t,x,a}(\omega^1,\omega)) - \varphi(X_{s''}^{t,x,a}(\omega^1,\omega))\big|\,(\P^1\otimes\P)(d\omega^1,d\omega)\bigg)\pi(dx) \\
&= \ \int_\Omega\bigg(\int_{\Omega^1\times\R^n} \lim_{\delta\downarrow0}\sup_{\substack{|s' - s''|<\delta\\ s',s''\in\Q}}\big|\varphi(X_{s'}^{t,x,a}(\omega^1,\omega)) - \varphi(X_{s''}^{t,x,a}(\omega^1,\omega))\big|\,\P^1(d\omega^1)\,\pi(dx)\bigg)\P(d\omega).
\end{align*}
We conclude that there exists a $\P$-null set $N\in\Fc$ such that for every $\omega\notin N$ the map $s\mapsto\varphi(X_s^{t,x,a}(\omega^1,\omega))$ is $\P^1(d\omega^1)\pi(dx)$-a.e. continuous, so that
\begin{equation}\label{int_sup_s'_s''}
\int_{\Omega^1\times\R^n} \lim_{\delta\downarrow0}\sup_{\substack{|s' - s''|<\delta\\ s',s''\in\Q}}\big|\varphi(X_{s'}^{t,x,a}(\omega^1,\omega)) - \varphi(X_{s''}^{t,x,a}(\omega^1,\omega))\big|\,\P^1(d\omega^1)\,\pi(dx) \ = \ 0, \qquad \forall\,\omega\notin N.
\end{equation}
In particular, this implies that, for every $\omega\notin N$, we have: $\P^1(d\omega^1)\pi(dx)$-a.e.
\begin{align}\label{sup_s'_s''}
&\sup_{|s' - s''|<\delta}\big|\varphi(X_{s'}^{t,x,a}(\omega^1,\omega)) - \varphi(X_{s''}^{t,x,a}(\omega^1,\omega))\big| \notag \\
&= \ \sup_{\substack{|s' - s''|<\delta\\ s',s''\in\Q}}\big|\varphi(X_{s'}^{t,x,a}(\omega^1,\omega)) - \varphi(X_{s''}^{t,x,a}(\omega^1,\omega))\big|.
\end{align}
Now, we observe that, by formula \eqref{rho}, we have
\[
\big|\rho_{s'}^{t,\pi,a}(\omega) - \rho_{s''}^{t,\pi,a}(\omega)\big| \ \leq \ \int_{\Omega^1\times\R^n} \big|\varphi(X_{s'}^{t,x,a}(\omega^1,\omega)) - \varphi(X_{s''}^{t,x,a}(\omega^1,\omega))\big| \,\P^1(d\omega^1)\,\pi(dx).
\]
Fix $\delta>0$ and take the supremum over all pairs $s',s''$ such that $|s'-s''|<\delta$, then, for every $\omega\notin N$,
\begin{align}\label{Continuity_rho}
&\sup_{|s' - s''|<\delta}\big|\rho_{s'}^{t,\pi,a}(\omega) - \rho_{s''}^{t,\pi,a}(\omega)\big| \notag \\
&\leq \ \int_{\Omega^1\times\R^n} \sup_{|s' - s''|<\delta}\big|\varphi(X_{s'}^{t,x,a}(\omega^1,\omega)) - \varphi(X_{s''}^{t,x,a}(\omega^1,\omega))\big| \,\P^1(d\omega^1)\,\pi(dx) \notag \\
&= \ \int_{\Omega^1\times\R^n} \sup_{\substack{|s' - s''|<\delta\\ s',s''\in\Q}}\big|\varphi(X_{s'}^{t,x,a}(\omega^1,\omega)) - \varphi(X_{s''}^{t,x,a}(\omega^1,\omega))\big| \,\P^1(d\omega^1)\,\pi(dx),
\end{align}
where we used property \eqref{sup_s'_s''}. Taking the limit as $\delta\downarrow0$ in \eqref{Continuity_rho}, we obtain, for every $\omega\notin N$,
\begin{align}\label{Continuity_rho2}
&\lim_{\delta\downarrow0}\sup_{|s' - s''|<\delta}\big|\rho_{s'}^{t,\pi,a}(\omega) - \rho_{s''}^{t,\pi,a}(\omega)\big| \notag \\
&\leq \ \int_{\Omega^1\times\R^n} \lim_{\delta\downarrow0}\sup_{\substack{|s' - s''|<\delta\\ s',s''\in\Q}}\big|\varphi(X_{s'}^{t,x,a}(\omega^1,\omega)) - \varphi(X_{s''}^{t,x,a}(\omega^1,\omega))\big| \,\P^1(d\omega^1)\,\pi(dx).
\end{align}
From \eqref{int_sup_s'_s''}, we deduce that
\[
\lim_{\delta\downarrow0}\sup_{|s' - s''|<\delta}\big|\rho_{s'}^{t,\pi,a}(\omega) - \rho_{s''}^{t,\pi,a}(\omega)\big| \ = \ 0, \qquad \forall\,\omega\notin N.
\]
This proves that, for every $\omega\notin N$, the map $s\mapsto\rho_s^{t,\pi,a}(\omega)[\varphi]$ is continuous, and concludes the proof of item (ii.a).

Let us now prove item (ii.b). To this end, fix $\xi\in L^2(\Fc^{o};\R^n)$, with $\pi=\Lc(\xi)$. Recall from Lemma \ref{L:xi} that the process
$(X_s^{t,\xi,a})_{s\in[t,T]}$ is $\bar\F^{o,B,\mu,t}$-adapted. As a consequence, $(X_s^{t,\xi,a})_{s\in[t,T]}$ is adapted (and hence predictable, since it is continuous) with respect to the $\bar\P$-completion of the filtration
$(\Fc^o\otimes\Fc^1\otimes\Fc_s^{W,\mu,t})_{s\geq t}$. Then, by Proposition 1.1.(b) of \cite{jacod_book} it follows that there exists a $\bar\P$-null set $\bar N\subset\bar\Omega$ and a process $(\bar X_s^{t,\xi,a})_{s\in[t,T]}$ such that:
\begin{itemize}
\item $(\bar X_s^{t,\xi,a})_{s\in[t,T]}$ is $(\Fc^o\otimes\Fc^1\otimes\Fc_s^{W,\mu,t})_{s\geq t}$-adapted;
\item $(X_s^{t,\xi,a})_{s\in[t,T]}$ and $(\bar X_s^{t,\xi,a})_{s\in[t,T]}$ coincide on the complementary set of $\bar N$, namely they are $\bar\P$-indistinguishable.
\end{itemize}
In particular, there exists a $\P$-null set $N\subset\Omega$ such that, for every $\omega\notin N$, the processes $(X_s^{t,\xi,a}(\cdot,\omega))_{s\in[t,T]}$ and $(\bar X_s^{t,\xi,a}(\cdot,\omega))_{s\in[t,T]}$ are $\bar\P^1$-indistinguishable. Therefore, for every $s\in[t,T]$ and continuous $\varphi\in\mathscr B_{_2}(\R^n)$, we have
\beqs
\rho_s^{t,\pi,a}(\omega)[\varphi] \ = \ \E^1\Big[\int_{\R^n}\varphi(X_s^{t,x,a}(\cdot,\omega))\,\pi(dx)\Big]  &=&  \bar\E^1\big[\varphi(X_s^{t,\xi(\cdot),a}(\cdot,\omega))\big]  =  \bar\E^1\big[\varphi(\bar X_s^{t,\xi(\cdot),a}(\cdot,\omega))\big],
\enqs
where the last equality holds for every $\omega\notin N$. Now notice that, by Fubini's theorem, the map $\omega\mapsto\bar\E^1[\varphi(\bar X_s^{t,\xi(\cdot),a}(\cdot,\omega))]$ is $\Fc_s^{W,\mu,t}$-measurable. Since the random variable
$\omega\mapsto\rho_s^{t,\pi,a}(\omega)[\varphi]$ is equal $\P$-a.s. to $\omega\mapsto\bar\E^1[\varphi(\bar X_s^{t,\xi(\cdot),a}(\cdot,\omega))]$, $\rho_s^{t,\pi,a}[\varphi]$ is also $\Fc_s^{W,\mu,t}$-measurable. This implies that the stochastic process $(\rho_s^{t,\pi,a}[\varphi])_{s\in[t,T]}$ is $\F^{W,\mu,t}$-adapted, so that item (ii.b) holds.  Finally the estimate \reff{estimrhopi} follows from \reff{estimXrandom}, and this concludes the proof of property (ii).
\ep

\begin{Lemma} \label{lemflow}
Under Assumption {\bf (H1)}, for every $(t,\pi,a)\in[0,T]\times\Pc_{_2}(\R^n)\times A$, the flow property holds:
\begin{equation}\label{FlowProp_rho}
\rho_s^{t,\pi,a} \ = \ \rho_s^{\theta,\rho_\theta^{t,\pi,a},I_\theta^{t,a}},
\end{equation}
$\P$-a.s., for all $\theta,s\in[t,T]$, with $\theta\leq s$.
\end{Lemma}
\textbf{Proof.}
We recall from Theorem 2.18 in \cite{BaCri} that there exists a countable separating class $(\varphi_m)_m$ of continuous and bounded functions from $\R^n$ into $\R$, namely (see Definition 2.12 in \cite{BaCri}) given $\pi_1,\pi_2\in\Pc_{_2}(\R^n)$ the condition $\pi_1(\varphi_m)=\pi_2(\varphi_m)$, for every $m$, implies that $\pi_1=\pi_2$. As a consequence, the claim of the Lemma follows if we prove that
\begin{equation}\label{Flow_phi_m}
\rho_s^{t,\pi,a}(\omega)[\varphi_m] \ = \ \rho_s^{\theta,\rho_\theta^{t,\pi,a}(\omega),I_\theta^{t,a}(\omega)}(\omega)[\varphi_m], \qquad \P(d\omega)\text{-a.s., for every }m.
\end{equation}
Fix $\varphi\in(\varphi_m)_m$. Our aim is to prove \eqref{Flow_phi_m} for $\varphi$. From the definition of $\rho_s^{t,\pi,a}(\omega)[\varphi]$ in \eqref{rho} and Remark \ref{R:Filter}, we have
\[
\rho_s^{t,\pi,a}(\omega)[\varphi] \ = \ \E^1\bigg[\int_{\R^n}\varphi(X_s^{t,x,a}(\cdot,\omega))\,\pi(dx)\bigg] \ = \ \bar\E^1\big[ \varphi\big(X_s^{t,\xi(\cdot),a}(\cdot,\omega)\big) \big].
\]
Now, recall from \eqref{flowXI} that $X_s^{t,\xi,a}=X_s^{\theta,X_\theta^{t,\xi,a},I_\theta^{t,a}}$, $\P^\Rc$-a.s., therefore there exists a $\P$-null set $N\in\Fc$ such that $X_s^{t,\xi(\cdot),a}(\cdot,\omega)=X_s^{\theta,X_\theta^{t,\xi(\cdot),a}(\cdot,\omega),I_\theta^{t,a}(\omega)}(\cdot,\omega)$, $\P^o\otimes\P^1$-a.s., for every $\omega\notin N$. Hence, for every $\omega\notin N$, we have
\begin{align*}
\rho_s^{t,\pi,a}(\omega)[\varphi] \ &= \ \bar\E^1\Big[\varphi\Big(X_s^{\theta,X_\theta^{t,\xi(\cdot),a}(\cdot,\omega),I_\theta^{t,a}(\omega)}(\cdot,\omega)\Big)\Big] \\
&= \ \E^1\bigg[\int_{\R^n}\varphi\Big(X_s^{\theta,X_\theta^{t,x,a}(\cdot,\omega),I_\theta^{t,a}(\omega)}(\cdot,\omega)\Big)\,\pi(dx)\bigg] \ = \ \int_{\R^n} \E^1\big[\Phi(X_\theta^{t,x,a}(\cdot,\omega),\cdot)\big]\,\pi(dx),
\end{align*}
where
\begin{equation}\label{Phi}
\Phi(y,\omega^1) \ := \ \varphi\big(X_s^{\theta,y,I_\theta^{t,a}(\omega)}(\omega^1,\omega)\big),
\end{equation}
for every $(y,\omega^1)\in\R^n\times\Omega$ (we omit the dependence of $\Phi$ on $t,\theta,s,a,\omega$). Notice that, for every $(t,\theta,s,\omega,a)$, the $\sigma$-algebra on $\Omega^1$ given by $\sigma(X_\theta^{t,x,a}(\cdot,\omega),\,x\in\R^n)$ is independent of $\sigma(\Phi(y,\cdot),\,y\in\R^n)$. Therefore, by a monotone class argument (first taking $\Phi$ of the form $\Phi(y,\omega)=\Phi_1(y)\Phi_2(\omega)$, for some measurable and bounded functions $\Phi_1$ and $\Phi_2$), we see that
\begin{align*}
\E^1\big[\Phi(X_\theta^{t,x,a}(\cdot,\omega),\cdot)\big] \ &= \ \int_{\Omega^1} \Phi(X_\theta^{t,x,a}(\omega^1,\omega),\omega^1)\,\P^1(d\omega^1) \\
&= \ \int_{\Omega^1\times\tilde\Omega^1} \Phi(X_\theta^{t,x,a}(\tilde\omega^1,\omega),\omega^1)\,(\P^1\otimes\tilde\P^1)(d\omega^1,d\tilde\omega^1),
\end{align*}
where $(\tilde\Omega^1,\tilde\Fc^1,\tilde\P^1)$ is a copy of $(\Omega^1,\Fc^1,\P^1)$. In conclusion, we find, $\P(d\omega)$-a.s.,
\begin{equation}\label{rho_s}
\rho_s^{t,\pi,a}(\omega)[\varphi] \ = \ \int_{\R^n}\int_{\Omega^1\times\tilde\Omega^1} \Phi(X_\theta^{t,x,a}(\tilde\omega^1,\omega),\omega^1) (\P^1\otimes\tilde\P^1)(d\omega^1,d\tilde\omega^1)\,\pi(dx).
\end{equation}
On the other hand, from formula \eqref{rho} we have
\begin{equation}\label{rho_s^rho_theta}
\rho_s^{\theta,\rho_\theta^{t,\pi,a}(\omega),I_\theta^{t,a}(\omega)}(\omega)[\varphi] \ = \ \E^1\bigg[\int_{\R^n} \varphi\big(X_s^{\theta,y,I_\theta^{t,a}(\omega)}(\cdot,\omega)\big)\,\rho_\theta^{t,\pi,a}(\omega)(dy) \bigg].
\end{equation}
Notice that $\rho_\theta^{t,\pi,a}(\omega)$, which is defined as
\[
\rho_\theta^{t,\pi,a}(\omega)[\psi] \ := \ \E^1\bigg[\int_{\R^n} \psi(X_\theta^{t,x,a}(\cdot,\omega)) \,\pi(dx)\bigg], \qquad \text{for every }\psi\in\mathscr B_{_2}(\R^n),
\]
is the distribution of the random variable $X_s^{t,\xi(\cdot),a}(\cdot,\omega)\colon\Omega^1\rightarrow\R^n$. Therefore, \eqref{rho_s^rho_theta} can be rewritten as (recalling that $\pi$ is the distribution of $\xi$)
\begin{align*}
\rho_s^{\theta,\rho_\theta^{t,\pi,a}(\omega),I_\theta^{t,a}(\omega)}(\omega)[\varphi] \ &= \ \E^1\bigg[\int_{\R^n\times\tilde\Omega^1} \varphi(X_s^{\theta,X_\theta^{t,x,a}(\tilde\omega^1,\omega),I_\theta^{t,a}(\omega)}(\omega^1,\omega)) (\pi\otimes\tilde\P^1)(dx,d\tilde\omega^1)\bigg] \\
&= \ \int_{\R^n}\int_{\Omega^1\times\tilde\Omega^1} \varphi(X_s^{\theta,X_\theta^{t,x,a}(\tilde\omega^1,\omega),I_\theta^{t,a}(\omega)}(\omega^1,\omega)) (\P^1\otimes\tilde\P^1)(d\omega^1,d\tilde\omega^1)\,\pi(dx),
\end{align*}
which coincides with \eqref{rho_s}, using the definition of $\Phi$ in \eqref{Phi}. This concludes the proof.
\ep

\vspace{3mm}

For sake of completeness, we end this paragraph by showing how one can derive  a filtering equation, namely an evolution equation for the randomized filter process $\rho$, although it will not be used in our approach.
We denote by $L^a$ and $M^a$,  $a$ $\in$ $A$, the second-order and first-order differential operators from $\Cc_b^2(\R^n)$ ($\subset$ $\mathscr B_{_2}(\R^n)$),
the set of bounded and twice continuously differentiable functions on $\R^n$ with bounded derivatives, into  $\mathscr B_{_2}(\R^n)$:
\beqs
L^a\varphi  & = &  b(.,a)\cdot D_x\varphi  + \frac{1}{2}\text{tr}\big[\sigma\sigma\trans(.,a)D_x^2\varphi \big], \\
M^a\varphi & = & D_x\varphi \trans \sigma_{_W} (.,a).
\enqs
%for any $\varphi$ $\in$ $C_b^2(\R^n)$, the set of twice continuously differentiable functions on $\R^n$ with bounded derivatives.
%Here $.$ denotes the scalar product, and $\trans$ denotes the transpose of  a vector or matrix.
By It\^o's formula to \reff{dynXI}, given $\varphi\in \Cc_b^2(\R^n)$, we have
\beq
\varphi(X_s^{t,\xi,a}) & = & \varphi(\xi) + \int_t^s L^{I_r^{t,a}}\varphi(X_r^{t,\xi,a})\,dr + \int_t^s D_x\varphi(X_r^{t,\xi,a})\trans \sigma_{_V} (X_r^{t,\xi,a},I_r^{t,a})\,dV_r \nonumber \\
& & \hspace{2cm}   + \;  \int_t^s M^{I_r^{t,a}}\varphi(X_r^{t,\xi,a})\,dW_r.   \label{itoforfiltering}
\enq
Taking in \eqref{itoforfiltering} the conditional expectation with respect to the $\sigma$-algebra $\bar\Fc_\infty^{W,\mu,t}$, introduced in Remark \ref{R:V_W,mu,t}, we obtain by \eqref{rho} and Remark \ref{R:Filter}
\beq\label{ZakaiRandomized}
\rho_s^{t,\pi,a} [\varphi]&=& \pi[\varphi] + \int_t^s \rho_r^{t,\pi,a}[L^{I_r^{t,a}}\varphi]\,dr +  \int_t^s\rho_r^{t,\pi,a} [M^{I_r^{t,a}}\varphi]\,dW_r,
\enq
for all $s\in[t,T]$. Observe that equation \eqref{ZakaiRandomized} looks like a controlled Zakai equation (see  equation \reff{Zakaicontrol1})
where the control parameter has been replaced by the stochastic process $I^{t,a}$. We refer to \eqref{ZakaiRandomized} as \emph{randomized Zakai equation}. Following standard arguments in nonlinear filtering (see for instance Chapter 7 in \cite{BaCri}), we can recast equation \eqref{ZakaiRandomized} in a strong form under the condition that, for all $s\in[t,T]$, $\rho_s^{t,\pi,a}$ has a density $p_s^{t,\pi,a}$ $\in$  $L^2(\R^n)$ (the $L^2$ space on $(\R^n,\Bc(\R^n))$ with respect to the Lebesgue measure), namely
\beqs
\rho_s^{t,\pi,a}[\varphi] &=&   \int_{\R^n} \varphi(x)\,p_s^{t,\pi,a}(x)\,dx \;\, =: \;\, <\varphi,p_s^{t,\pi,a}>_{_{L^2(\R^n)}}, \qquad \text{for every }s\in[t,T],
\enqs
for all bounded measurable functions $\varphi$ on $\R^n$. In particular, $\pi$ has a density $\eta$ $:=$ $p_t^{t,\pi,a}$ $\in$ $L^2(\R^n)$, and we set $p_s^{t,\eta,a}$ $:=$ $p_s^{t,\pi,a}$.  If  the function $x$ $\mapsto$ $p_s^{t,\eta,a}(x)$ is smooth enough (for instance, $p_s^{t,\eta,a}$ $\in$
$\Cc_c^2(\R^n)$ $\subset$ $L^2(\R^n)$, the set of twice continuously differentiable functions on $\R^n$ with compact support), integrating by parts in \eqref{ZakaiRandomized} we obtain:
\beqs
\rho_s^{t,\pi,a}[\varphi] & = & < \varphi, p_s^{t,\eta,a} > _{_{L^2(\R^n)}}  \\
& = & < \varphi, \Big(\eta + \int_t^s(L^{I_r^{t,a}})^*p_r^{t,\eta,a} \,dr + \int_t^s(M^{I_r^{t,a}})^*p_r^{t,\eta,a}\,dW_r\Big) >_{_{L^2(\R^n)}},
\enqs
for all $\varphi$ $\in$ $\Cc_c^2(\R^n)$, where
\begin{equation}\label{Ladjoint}
\begin{array}{ccl}
(L^a)^* p &=& - {\rm div}( b(.,a) p ) + \frac{1}{2} {\rm tr}\big[ D_x^2(\sigma\sigma\trans(.,a) p \big) \big], \\
% - \sum_{i=1}^n \frac{\partial(b_i(t,x,a)\varphi(x))}{\partial x_i} + \frac{1}{2} \sum_{i,j=1}^n \frac{\partial^2((\sigma\sigma\trans)_{ij}(t,x,a)\varphi(x))}{\partial x_i\partial x_j}, \\
(M^a)^* p  &=& - {\rm div}\big( \sigma_{_W}(.,a) p \big),
%\ - \bigg(\sum_{i=1}^n \frac{\partial(\sigma_{ij}^{(2)}(t,x,a)\varphi(x))}{\partial x_i}\bigg)_{1\leq j\leq d},
\end{array}
\end{equation}
are the adjoint operators of $L^a$ and $M^a$ (here ${\rm div}(.)$ denotes the divergence operator). This leads to consider the stochastic partial differential equation (SPDE)  for the density
$p^{t,\eta,a}$ in the Hilbert space   $L^2(\R^n)$:
\begin{equation} \label{DuncanMZ}
dp_s^{t,\eta,a} \ = \ (L^{I_s^{t,a}})^* p_s^{t,\eta,a}\,ds + (M^{I_s^{t,a}})^* p_s^{t,\eta,a}\,dW_s, \;\;\; t \leq s \leq T, \qquad  p_t^{t,\eta,a} \ = \ \eta \; \in \; L^2(\R^n),
\end{equation}
which is refereed to as the {\it randomized Duncan-Mortensen-Zakai (DMZ)  SPDE} (compare with the controlled DMZ equation in  \reff{DuncanMZcontrolintro}).

\subsection{Randomized dynamic programming principle}

From Bayes formula, and the $(\P,\F^{W,\mu,t})$-martingale property of the density process $\kappa^{\nu,t}$ in \reff{densitykappa} of the probability measure
$\P^{\nu,t}$, $\nu$ $\in$ $\Vc_{_{W,\mu,t}}$, we can  express the randomized gain functional \reff{defJrandomized} in a separated form involving  the randomized filter process \reff{rho}, namely:
\beqs
J^\Rc(t,\xi,\nu,a) &=& \E^{\nu,t} \bigg[ \int_t^T \rho_s^{t,\pi,a}[f(.,I_s^{t,a})] ds + \rho_T^{t,\pi,a}[g] \bigg].
\enqs
where $\xi\in L^2(\Fc^{o};\R^n)$ has distribution $\pi$. We then see  that the randomized gain functional $\Jc^\Rc$ depends upon $\xi$ only through its distribution $\pi$. Therefore, we define the inverse-lifted version of $\Jc^\Rc$ by:
\beqs
J^\Rc(t,\pi,\nu,a) & := & \Jc^\Rc(t,\xi,\nu,a) \;\, = \;\, \E^{\nu,t} \bigg[ \int_t^T \rho_s^{t,\pi,a}[f(.,I_s^{t,a})] ds + \rho_T^{t,\pi,a}[g] \bigg],
\enqs
for all $(t,\pi)\in[0,T]\times\Pc_{_2}(\R^n)$ and every $\nu\in\Vc_{_{W,\mu,t}}$. Similarly, we define the inverse-lifted versions of $\upv$ and $\upv^\Rc$ respectively as $v$ and $v^\Rc$. They are given by (recalling the equivalence result \reff{relvvrandom})
\beq \label{vpi}
v(t,\pi) \; = \; v^\Rc(t,\pi,a) \; = \; \sup_{\nu\in\Vc_{_{W,\mu,t}}} J^\Rc(t,\pi,\nu,a), \qquad (t,\pi) \in [0,T]\times\Pc_{_2}(\R^n),\,a\in A.
\enq
Noting that $v^\Rc$ is constant with respect to $a\in A$, from now on we simply write $v^\Rc(t,\pi)$ instead of $v^\Rc(t,\pi,a)$, for any $a\in A$.
Notice also from \reff{vquadra} that we have the quadratic growth condition: there exists some positive constant $C$ such that
\beq \label{vquadrapi}
|v(t,\pi)| & \leq & C (1 + \|\pi\|_{_2}^2 ), \;\;\; \forall (t,\pi) \in [0,T]\times\Pc_{_2}(\R^n).
\enq

The purpose of this paragraph is to prove a dynamic programming principle for the value function $v^\Rc$ of the  randomized control problem, and consequently for the value function $v$ of the primal partial observation control problem. We shall rely on the dual representation of $v^\Rc(t,\pi)$ in terms of a constrained backward stochastic differential equation (BSDE),  obtained in \cite{BCFP16b}, and refereed to as the \emph{randomized equation}.

For every $(t,\pi,a)\in[0,T]\times\Pc_{_2}(\R^n)\times A$, the randomized equation is formulated in terms of the filter process $(\rho_s^{t,\pi,a})_{s\in[t,T]}$ as a backward stochastic differential equation with nonnegative jumps on the filtered probability space $(\Omega,\Fc,\F^{W,\mu,t},\P)$:
\begin{equation}\label{BSDEconstrained}
\begin{cases}
\vspace{2mm} \dis Y_s \ = \ \rho_T^{t,\pi,a}[g]  + \int_s^T\rho_r^{t,\pi,a}[f_r(\cdot,I_r^{t,a})]\,dr + K_T - K_s \\
\vspace{2mm} \dis \qquad\;\;\;\, - \int_s^T Z_r\,dW_r - \int_s^T\!\int_A U_r(a)\,\mu(dr\,da), \qquad \text{for all }s\in[t,T], \\
\dis U_s(a) \ \leq \ 0, \qquad \P(d\omega)\,ds\,\lambda(da)\text{-a.e. on }\Omega\times[t,T]\times A.
\end{cases}
\end{equation}
We are interested in the minimal solution to \eqref{BSDEconstrained}, whose definition is reported below.

\begin{Definition}\label{D:BSDE}
For every $(t,\pi,a)\in[0,T]\times\Pc_{_2}(\R^n)\times A$, a quadruple $(Y,Z,U,K)$ is a solution to equation \eqref{BSDEconstrained} if:
\begin{enumerate}
\item $Y$ $\in$  $\Sc^2(\F^{W,\mu,t})$, the class of real c\`adl\`ag  $\F^{W,\mu,t}$-adapted processes $(Y_s)_{s\in[t,T]}$ such that $\|Y\|_{\Sc^2(\F^{W,\mu,t})}^2$ $:=$ $\E[\sup_{t\leq s\leq T}|Y_s|^2]$ $<$ $\infty$;
\item $Z$ $\in$ $L_W^2(\F^{W,\mu,t})$, the class of $\R^d$-valued $\Pc(\F^{W,\mu,t})$-measurable processes $(Z_s)_{s\in[t,T]}$ such that $\|Z\|_{L_W^2(\F^{W,\mu,t})}^2$ $:=$ $\E\big[\int_t^T|Z_s|^2ds\big]$ $<$ $\infty$;
  \item $U$ $\in$ $L_{\tilde\mu}^2(\F^{W,\mu,t})$, the class of $\Pc(\F^{W,\mu,t})\otimes\Bc(A)$-measurable maps $U\colon\Omega\times[t,T]\times A\rightarrow\R$ such that $\|U\|^2_{L_{\tilde\mu}^2(\F^{W,\mu,t})}$ $:=$ $\E\big[\int_t^T\int_A|U_s(a)|^2\lambda(da)ds\big]$ $<$ $\infty$;
  \item $K$ $\in$ $\Kc^2(\F^{W,\mu,t})$, the subclass of $\Sc^2(\F^{W,\mu,t})$ consisting of $\Pc(\F^{W,\mu,t})$-measurable nondecreasing processes satisfying $K_t=0$.
\end{enumerate}
A solution $(Y,Z,U,K)$ is called minimal if, given any other solution $(Y',Z',U',K')$, we have $Y_s\leq Y'_s$, $\P$-a.s., for all $s\in[t,T]$.
\end{Definition}

It is proved in \cite{BCFP16b}, Theorem 5.1, that for every $(t,\pi,a)\in[0,T]\times\Pc_{_2}(\R^n)\times A$ there exists a unique minimal solution $(Y^{t,\pi,a},Z^{t,\pi,a},U^{t,\pi,a},K^{t,\pi,a})$ $\in$
$\Sc^2(\F^{W,\mu,t})\times L_W^2(\F^{W,\mu,t})\times L_{\tilde\mu}(\F^{W,\mu,t})\times \Kc^2(\F^{W,\mu,t})$ to \reff{BSDEconstrained}. Moreover, we have the following dual representation (see Theorem 5.1 in \cite{BCFP16b}):
\beq \label{dualY}
Y_t^{t,\pi,a} &=& \sup_{\nu\in\Vc_{_{W,\mu,t}}} \E^{\nu,t} \bigg[ \int_t^T \rho_s^{t,\pi,a}[f(\cdot,I_s^{t,a})]\,ds + \rho_T^{t,\pi,a}[g] \bigg].
\enq
This implies, using \reff{vpi}, that (recall that we write $v^\Rc(t,\pi)$ instead of $v^\Rc(t,\pi,a)$, for any $a\in A$)
\beq \label{relvY}
v(t,\pi) \;\, = \;\, v^\Rc(t,\pi) &=& Y_t^{t,\pi,a}, \qquad (t,\pi) \in [0,T]\times\Pc_{_2}(\R^n),\text{ for any }a\in A.
\enq
As shown in Theorem 5.2 of \cite{BCFP16b}, the following generalization of the dual representation \reff{dualY} holds:
\beq \label{dualYdyn}
Y_t^{t,\pi,a} &=& \sup_{\nu\in\Vc_{_{W,\mu,t}}} \sup_\theta \E^{\nu,t} \bigg[ \int_t^\theta \rho_r^{t,\pi,a}[f(\cdot,I_r^{t,a})]\,dr + Y_\theta^{t,\pi,a} \bigg] \notag \\
&=& \sup_{\nu\in\Vc_{_{W,\mu,t}}} \inf_\theta \E^{\nu,t} \bigg[ \int_t^\theta \rho_r^{t,\pi,a}[f(\cdot,I_r^{t,a})]\,dr + Y_\theta^{t,\pi,a} \bigg],
\enq
where the $\sup_\theta$ and $\inf_\theta$ are taken over the class of $\F^{W,\mu,t}$-stopping times $\theta$ valued in $[t,T]$,
and this can be viewed as a dynamic programming relation for the solution to the BSDE.  The importance of the flow property in Lemma \ref{lemflow} is now to state  in the next Lemma the identification between the value function and the
solution to the BSDE at any time, which in turn allows to obtain the randomized dynamic programming principle for the value function.

\begin{Lemma}
Under Assumption {\bf (H1)}, for every $(t,\pi,a)\in[0,T]\times\Pc_{_2}(\R^n)\times A$, we have the identity: $\P$-a.s.,
\begin{equation}\label{vYdyn}
v(s,\rho_s^{t,\pi,a}) \ = \ v^\Rc(s,\rho_s^{t,\pi,a}) \ = \ Y_s^{t,\pi,a}.
\end{equation}
\end{Lemma}
\textbf{Proof.}
We provide a sketch of the proof. We introduce, for every $m\in\N$, the penalized backward stochastic differential equation:
\begin{align*}
Y_s^{m,t,\pi,a} \ &= \ \rho_T^{t,\pi,a}[g]  + \int_s^T\rho_r^{t,\pi,a}[f_r(\cdot,I_r^{t,a})]\,dr + K_T^{m,t,\pi,a} - K_s^{m,t,\pi,a} \\
&\quad \ - \int_s^T Z_r^{m,t,\pi,a}\,dW_r - \int_s^T\int_A U_r^{m,t,\pi,a}(a)\,\mu(dr\,da),
\end{align*}
for every $s\in[t,T]$, where
\[
K_s^{m,t,\pi,a} \ = \ m \int_t^s \int_A \big(U_r^{m,t,\pi,a}(a)\big)^+ \,\lambda(da)\,dr.
\]
By Lemma 2.4 in \cite{tang_li} it follows that, for every $m\in\N$, there exists a unique solution $(Y^{m,t,\pi,a},$ $Z^{m,t,\pi,a},U^{m,t,\pi,a})$ $\in$ $\Sc^2(\F^{W,\mu,t})\times L_W^2(\F^{W,\mu,t})\times L_{\tilde\mu}^2(\F^{W,\mu,t})$ to the penalized equation. Since this latter equation is not a constrained equation, it is easy to prove, by the flow properties \eqref{flowI} and \eqref{FlowProp_rho}, that the following flow property holds:
\begin{equation}\label{FlowProp_Ym}
Y_s^{m,t,\pi,a} \ = \ Y_s^{m,\theta,\rho_\theta^{t,\pi,a},I_\theta^{t,a}},
\end{equation}
$\P$-a.s., for all $\theta,s\in[t,T]$, with $\theta\leq s$. Let us now define, for every $m\in\N$, the deterministic map:
\[
v_m(t,\pi,a) \ = \ Y_t^{m,t,\pi,a}, \qquad \text{for every }(t,\pi,a)\in[0,T]\times\Pc_{_2}(\R^n)\times A.
\]
Then, from \eqref{FlowProp_Ym}, we obtain:
\begin{equation}\label{vYdyn_m}
v_m(s,\rho_s^{t,\pi,a},I_s^{t,a}) \ = \ Y_s^{m,t,\pi,a}, \qquad \P\text{-a.s., for every }(t,\pi,a)\in[0,T]\times\Pc_{_2}(\R^n)\times A.
\end{equation}
We now notice that, proceeding as in Theorem 2.1 of \cite{KP12}, we see that $Y_s^{m,t,\pi,a}$ increasingly converges $\P$-a.s. to $Y_s^{t,\pi,a}$, as $m$ goes to infinity, for every $s\in[t,T]$. In particular, $v_m(t,\pi,a)=Y_t^{m,t,\pi,a}$ converges to $v(t,\pi)=v^\Rc(t,\pi)=Y_t^{t,\pi,a}$, for every $(t,\pi,a)\in[0,T]\times\Pc_{_2}(\R^n)\times A$. Then, letting $m\rightarrow\infty$ in \eqref{vYdyn_m}, we obtain \eqref{vYdyn}.
\ep

\vspace{3mm}

By \reff{dualYdyn} and formula \eqref{vYdyn},  we then  obtain the so-called {\it randomized dynamic programming principle} for the value functions $v$ $=$ $v^\Rc$:
\begin{align} \label{DPPrandom}
v(t,\pi) \ = \ \sup_{\nu\in\Vc_{_{W,\mu,t}}} \sup_\theta \E^{\nu,t} \bigg[ \int_t^\theta  \rho_r^{t,\pi,a}[f(\cdot,I_r^{t,a})] dr + v(\theta,\rho_\theta^{t,\pi,a}) \bigg] \notag \\
\ = \ \sup_{\nu\in\Vc_{_{W,\mu,t}}} \inf_\theta \E^{\nu,t} \bigg[ \int_t^\theta  \rho_r^{t,\pi,a}[f(\cdot,I_r^{t,a})] dr + v(\theta,\rho_\theta^{t,\pi,a}) \bigg],
\end{align}
for every $(t,\pi)\in[0,T]\times\Pc_{_2}(\R^n)$, $a\in A$, with  $\theta$ varying in the class of $\F^{W,\mu,t}$-stopping times valued in $[t,T]$.

\begin{Remark}
\label{dynpro}
{\rm   It would be possible to derive a standard DPP from the randomized version. However, there is no advantage in doing it, since the randomized DPP is as powerful as the standard one for our purpose
of characterizing the value function through a dynamic progra\-mming  equation in the Wasserstein space of probability measures.  This is the focus of the next section.
\ep
}
\end{Remark}

\section{Bellman equation and viscosity solutions  in the Wasserstein space} \label{secHJB}

\setcounter{equation}{0} \setcounter{Assumption}{0}
\setcounter{Theorem}{0} \setcounter{Proposition}{0}
\setcounter{Corollary}{0} \setcounter{Lemma}{0}
\setcounter{Definition}{0} \setcounter{Remark}{0}

The purpose of this section is to derive from the dynamic programming principle \reff{DPPrandom} a partial differential equation (PDE), called Hamilton-Jacobi-Bellman (HJB) equation,  in the Wasserstein space of probability measures for characterizing the value function $v(t,\pi)$.

\subsection{Notion of differentiability and It\^o's formula}

We shall rely on the notion of derivative with respect to a probability measure, as introduced by P.L. Lions in his course at Coll\`ege de France \cite{lio12}, and detailed   in the lecture notes  \cite{car12}.
% see also \cite{buetal14}, \cite{chacridel15}.

%This notion is based on the lifting of functions $u$ $:$ $\Pc_{_2}(\R^n)$ $\rightarrow$ $\R$ into functions $\upsilon$ defined on
%$L^2(\tilde\Fc;\R^n)$ ($=$ $L^2(\tilde\Omega,\tilde\Fc,\tilde\P;\R^n)$)  over some probability space $(\tilde\Omega,\tilde\Fc,\tilde\P)$  by  setting
%$\upsilon(X)$ $=$  $u(\tilde\P_{_X})$,  $\tilde\Omega$ being a Polish space and $\tilde\P$ an atomless probability measure.
A function $u$ $:$ $\Pc_{_2}(\R^n)$ $\rightarrow$ $\R$ is said to be differentiable  (resp. $\Cc^1$) on
$\Pc_{_2}(\R^n)$ if its lifted function $\upsilon$ is Fr\'echet differentiable (resp. Fr\'echet differentiable with continuous derivatives) on
$L^2(\Fc^o;\R^n)$.  In this case, the Fr\'echet derivative $[D\upsilon](\xi)$, viewed as an element $D\upsilon(\xi)$ of $L^2(\Fc^o;\R^n)$  by Riesz' theorem:
$[D\upsilon](\xi)(\zeta)$ $=$ $\E^o[D\upsilon(\xi).\zeta]$,  can be represented as
\beq \label{Uu1}
D\upsilon(\xi) &=& \partial_\pi u(\Lc(\xi))(\xi),
\enq
for some function  $\partial_\pi u(\Lc(\xi))$ $:$ $\R^n$ $\rightarrow$ $\R^n$,  which is  called derivative of $u$ at $\pi$ $=$ $\Lc(\xi)$.  Moreover,
$\partial_\pi u(\pi)$ $\in$ $L^2_\pi(\R^n)$ for  $\pi$ $\in$ $\Pc_{_2}(\R^n)$ $=$ $\{\Lc(\xi), \xi \in L^2(\Fc^o;\R^n)\}$.
Following \cite{chacridel15}, we say that $u$ is fully  $\Cc^2$ if it is $\Cc^1$, and one can find, for any $\pi$ $\in$ $\Pc_{_2}(\R^n)$, a continuous version of the mapping $x\in\R^n$ $\mapsto$ $\partial_\pi u(\pi)(x)$, such that the mapping
$(\pi,x)$ $\in$ $\Pc_{_2}(\R^n)\times\R^n$ $\mapsto$ $\partial_\pi u(\pi)(x)$  is continuous at any point $(\pi,x)$ such that $x$ $\in$ Supp$(\pi)$, and
\begin{itemize}
\item[(i)] for each fixed $\pi$ $\in$ $\Pc_{_2}(\R^n)$, the mapping $x$ $\in$ $\R^n$ $\mapsto$  $\partial_\pi u(\pi)(x)$ is differentiable in the standard sense, with a gradient denoted by  $\partial_x  \partial_\pi u(\pi)(x)$  $\in$ $\R^{n\times n}$, and such that the mapping  $(\pi,x)$ $\in$ $\Pc_{_2}(\R^n)\times\R^n$
$\mapsto$  $\partial_x  \partial_\pi u(\pi)(x)$ is continuous
\item[(ii)] for each fixed  $x$ $\in$ $\R^d$, the mapping $\mu$ $\in$ $\Pc_{_2}(\R^n)$ $\mapsto$  $\partial_\pi u(\pi)(x)$ is differentiable in the above lifted sense.  Its derivative, interpreted thus as a mapping $x'$ $\in$ $\R^d$ $\mapsto$ $\partial_\pi \big[ \partial_\pi u(\pi)(x)\big](x')$ $\in$ $\R^{n\times n}$ in
$L^2_\pi(\R^{n\times n})$, is denoted by $x'$ $\in$ $\R^n$ $\mapsto$ $\partial_\pi^2 u(\pi)(x,x')$, and such that the mapping $(\pi,x,x')$ $\in$
$\Pc_{_2}(\R^n)\times\R^n\times\R^n$ $\mapsto$ $\partial_\pi^2 u(\pi)(x,x')$ is continuous.
\end{itemize}
We say that $u$ $\in$ $\Cc^2_b(\Pc_{_2}(\R^n))$ if it is fully $\Cc^2$,  $\partial_x  \partial_\pi u(\pi)$ $\in$ $L_\pi^\infty(\R^{n\times n})$, $\partial_\pi^2 u(\pi)$ $\in$ $L_{\pi\otimes\pi}^\infty(\R^{n\times n})$ for any $\pi$
$\in$ $\Pc_{_2}(\R^n)$,  and for any compact set $\Kc$ of $\Pc_{_2}(\R^n)$, we have
%there exists some constant $C_u$ such that for all $\mu$ $\in$ $\Pc_{_2}(\R^d)$,
\beq \label{condpi}
 \sup_{ \pi \in \Kc } \Big[ \int_{\R^n} \big| \partial_\pi u(\pi)(x) |^2\pi(dx)  +
%\int_{\R^d}  \big| \partial_x \partial_\mu u(\mu)(x) |^2 \mu(dx)
\big \| \partial_x \partial_\pi u(\pi)\|_{_\infty} +  \big \| \partial_\pi^2 u(\pi)\|_{_\infty}
\Big]  & < & \infty.
%& \leq & C_u\big( 1 + \|\mu\|^2_{_2} \big).
\enq
%In this case, notice that $\partial_\mu u(\mu)$ (resp.  $\partial_x \partial_\mu u(\mu)$)  lies in $L_\mu^2(\R^d;\R^d)$ (resp. $L_\mu^2(\R^d;\R^{d\times d})$) the set functions from $\R^d$ into $\R^d$ (resp. $\R^{d\times d}$),
%and square-integrable with respect to $\mu$.
If $u$ lies in $\Cc^2_b(\Pc_{_2}(\R^n))$, then its  lifted function $\upsilon$ is twice continuously Fr\'echet differentiable on   $L^2(\Fc^o;\R^n)$.
In this case, the second Fr\'echet derivative $D^2 \upsilon(\xi)$  is identified indifferently by Riesz' theorem  as a  bilinear form on
$L^2(\Fc^o;\R^n)$ or  as a self-adjoint operator  (hence bounded) on  $L^2(\Fc^o;\R^n)$, denoted  by $D^2\upsilon(\xi)$ $\in$ $S(L^2(\Fc^o;\R^n))$,  and we have the relation (see Appendix A.2  in \cite{cardel14b}):
\begin{equation} \label{Uu2}
\left\{
\begin{array}{ccccl}
 D^2\upsilon(\xi)[\zeta, \zeta] & = & \E^o\Big[ D^2\upsilon(\xi)(\zeta).\zeta \Big]  &=& \E^o \Big[ \E^{'o} \big[ {\rm tr} \big( \partial_\pi^2 u(\Lc(\xi))(\xi,\xi') \zeta(\zeta')\trans \big) \big] \Big] \\
&&  & & \;\;\;\;\;  + \;\;  \E^o\Big[ {\rm tr} \big( \partial_x \partial_\pi u(\Lc(\xi))(\xi) \zeta\zeta\trans \big) \Big], \\
 D^2\upsilon(\xi)[\eta N, \eta N]  & = & \E^o \Big[ D^2\upsilon(\xi)(\eta N).\eta N \Big]  &=&
 \E^o\Big[ {\rm tr} \big( \partial_x \partial_\pi u(\Lc(\xi))(\xi) \eta\eta\trans \big) \Big],
 \end{array}
 \right.
\end{equation}
for any $\xi$ $\in$ $L^2(\Fc^o;\R^n)$, $\zeta$ $\in$ $L^2(\Fc^o;\R^{n})$,  $\eta$ $\in$ $L^2(\Fc^o;\R^{n\times q})$, and where $(\xi',\zeta')$ is a copy of
$(\xi,\zeta)$ on another Polish and atomless probability space $(\Omega^{'o},\Fc^{'o},\P^{'o})$, $N$ $\in$ $L^2(\Fc^o;\R^q)$ is  independent of
$(\xi,\eta)$ with zero mean  and unit variance.

%\marginpar{CHECK IF $C^2$ Frechet differentiable implies $C^2$ in the space of measures}

We next need an It\^o's formula along a flow of conditional measures proved in \cite{cardel14b} for processes with common noise, and that we adapt and formulate here in our randomization  context for  partial observation. Let  $(\bar\Omega,\bar\Fc,\bar\P)$ be a probability space of the form: $\bar\Omega$ $=$
$\bar\Omega^1\times\Omega$, $\bar\Fc$ $=$  $\bar\Fc^1\otimes\Fc$, $\bar\P$ $=$ $\bar\P^1\otimes\P$, where
$\bar\Omega^1$ $=$ $\Omega^o\times\Omega^1$, $\bar\Fc^1$ $=$ $\Fc^o\otimes\Fc^1$, $\bar\P^1$ $=$ $\P^o\otimes\P^1$, with
$(\Omega^o,\Fc^o,\P^o)$ a ``rich" probability space (see Remark \ref{R:OmegaPolish}), $(\Omega^1,\Fc^1,\P^1)$ supporting a $m$-dimensional
$(\P^1,\F^1=(\Fc_t^1)_t)$-Brownian motion $V$,  $(\Omega,\Fc,\P)$ supporting a $d$-dimensional $(\P,\F=(\Fc_t)_t)$-Brownian motion $W$.  For our purpose, $\F$ is the $\P$-completion of
either the filtration generated by $W$, or the filtration generated by $W$ and a Poisson random measure $\mu$. Denote by
$\bar\F$ $=$ $(\bar\Fc_t)_t$ the $\bar\P$-completion of the  filtration generated by $\Fc_t$, $\Fc_t^1$ and $\Fc^o$, $t$ $\geq$ $0$.
Let us consider an It\^o process in $\R^n$  of the form:
\beq \label{ItoX}
X_t &=& X_0 +  \int_0^t b_s ds +  \int_0^t \sigma_s^V dV_s +  \int_0^t \sigma_s^W dW_s, \;\;\; 0 \leq t \leq T,
\enq
where  $X_0$ $\in$ $L^2(\Fc^o,\R^n)$,  $b$, $\sigma$ $=$ $(\sigma^V \; \sigma^W)$ are $\bar\F$-progressively measurable processes,
and satisfying the integrability condition:
\beqs
\bar\E \big[ \int_0^T |b_t|^2 + |\sigma_t|^2 dt \big] & < & \infty.
\enqs
For $\omega$ $\in$ $\Omega$, denote by $\rho_t(\omega)$ the law of  the random variable $X_t(.,\omega)$ on
$(\bar\Omega^1,\bar\Fc^1,\bar\P^1)$. Thus $\rho_t$ is the conditional law of $X_t$ given $\Fc_\infty$, hence also equal to the conditional law of $X_t$ given $\Fc_t$ (recall that $\F$ is the natural filtration of $W$ or of $(W,\mu)$),
which then defines an $\F$-adapted process $(\rho_t)$ valued in $\Pc_{_2}(\R^n)$.
Let $u$ $\in$ $\Cc^2_b(\Pc_{_2}(\R^n))$. Then, for all $t$ $\in$ $[0,T]$, we have
\beq
u(\rho_t) &=& u(\rho_0)  + \int_0^t \bar\E^1 \Big[ \partial_\pi u(\rho_s)(X_s).b_s
+ \frac{1}{2}{\rm tr}\big(\partial_x \partial_\pi u(\rho_s)(X_s) \sigma_s\sigma_s\trans\big) \Big] \nonumber \\
& & \;\;\;\;\;\;\;  + \;  \bar\E^1  \Big[ \bar\E^{'1}
\big[ \frac{1}{2}{\rm tr}\big(\partial_\pi^2 u(\rho_s)(X_s,X_s') \sigma_s^W(\sigma_s^{',W})\trans\big) \big] \Big] ds  \nonumber \\
& & \;\;\; + \; \int_0^t \bar\E^1 \Big[ \partial_\pi  u(\rho_s)(X_s)\trans \sigma_s^W \Big] dW_s, \label{Ito}
\enq
where  $X'$ and $\sigma^{',W}$ are copies of  $X$ and $\sigma^W$ on another probability space
$(\bar\Omega'=\bar\Omega^{'1}\times\Omega,\bar\Fc^{'1}\otimes\Fc,\bar\P^{'1}\otimes\P)$,
with $(\bar\Omega^{'1},\bar\Fc^{'1},\bar\P^{'1})$ supporting $V'$ a copy of $V$.

Alternatively, we can formulate It\^o's formula for the lifted $\upsilon$ on $L^2(\Fc^o;\R^n)$. For this consider a copy $\tilde V$ of $V$ on
$(\Omega^o,\Fc^o,\P^o)$, and denote by $\tilde X_0$, $\tilde b$, $\tilde\sigma^V$, $\tilde\sigma^W$ copies of $X_0$, $b$, $\sigma^V$ on
$(\Omega^o\times\Omega,\Fc^o\otimes\Fc,\P^o\otimes\P)$, and consider the It\^o process
\beq\label{Ito_tildeX}
\tilde X_t &=& \tilde X_0 +  \int_0^t \tilde b_s ds +  \int_0^t \tilde\sigma_s^V d\tilde V_s +  \int_0^t \tilde\sigma_s^W dW_s, \;\;\; 0 \leq t \leq T,
\enq
which is then a copy on $(\Omega^o\times\Omega,\Fc^o\otimes\Fc,\P^o\otimes\P)$ of $X$ in \reff{ItoX}.
The process $\check X$ defined by $\check X_t(\omega)$ $=$ $\tilde X_t(.,\omega)$, $0\leq t\leq T$,  is $\F$-progressive, valued in
$L^2(\Fc^o;\R^n)$ with $\Lc(\check X_t(\omega))$ $=$ $\rho_t(\omega)$.
Similarly, the processes defined by $\check b_t(\omega)$ $=$ $\tilde b_t(.,\omega)$, $\check \sigma_t^V(\omega)$ $=$
$\tilde \sigma_t(.,\omega)$,  $\check \sigma^W_t(\omega)$ $=$ $\tilde \sigma^W_t(.,\omega)$, $0\leq t\leq T$, are $\F$-progressive
valued in $L^2(\Fc^o;\R^N)$,
$\P$-almost surely.  It\^o's formula \reff{Ito}  is then written for the lifted function $\upsilon$  $\in$ $\Cc^2(L^2(\Fc^o;\R^d))$ from \reff{Uu1}-\reff{Uu2}  as:
\beq
\upsilon(\check X_t) &=& \upsilon(\check X_0) +  \int_0^t \E^o\Big[ D\upsilon(\check X_s).\check b_s
+ \frac{1}{2} D^2\upsilon(\check X_s)(\check\sigma_s^V N).\check\sigma_s^V N +
 \frac{1}{2} \textup{tr}\big[D^2\upsilon(\check X_s)(\check\sigma_s^W).\check\sigma_s^W\big]   \Big] ds \nonumber \\
 & & \;\;\; + \; \int_0^t \E^o \big[ D\upsilon(\check X_s)\trans\check\sigma_s^W \big] dW_s, \;\;\;\;\;\;\;  0 \leq t \leq T,   \; \P\text{-a.s.} \label{Ito2}
\enq
for any   Gaussian random variable $N$ $\in$ $L^2(\Fc^o;\R^m)$ independent of $(\tilde V,\tilde X_0)$, with zero mean and unit variance. We do not report the proof of It\^o's formula \eqref{Ito2}, since it can be done proceeding along the same lines as in the proof of Propositions 6.3 and 6.5 in \cite{cardel14b}, see also Theorem 7.1 in \cite{BLPR14}, and \cite{chacridel15} for related results.

\subsection{Bellman equation and verification theorem}

 We have now the ingredients for deriving the Bellman equation for the (inverse-lifted) value function $v(t,\pi)$ in \reff{vpi}, and it turns out that it takes the following form:
 \begin{equation} \label{HJBpi}
 \left\{
 \begin{array}{rcl}
 \partial_t v +  \Sup_{a\in A} \Big\{ \pi\big[f(.,a) +  \Lc^a v(t,\pi)    \big]  + \pi\otimes\pi\big[ \Mc^a v(t,\pi) \big] \Big\} &=& 0, \;\;\; (t,\pi) \in
 [0,T)\times\Pc_{_2}(\R^n),  \\
 v(T,\pi) &=& \pi[g], \;\;\; \pi \in \Pc_{_2}(\R^n),
 \end{array}
 \right.
 \end{equation}
 where for  $\phi$ $\in$ $\Cc^2_b(\Pc_{_2}(\R^n))$, $a$ $\in$ $A$, and $\pi$ $\in$ $\Pc_{_2}(\R^n)$, $\Lc^a \phi(\pi)$ $\in$ $L_\pi^2(\R)$ is the function
 $\R^n$ $\rightarrow$ $\R$ defined by
\beq \label{defLc}
\Lc^a\phi(\pi)(x) &=& \partial_\pi \phi(\pi)(x).b(x,a) + \frac{1}{2}{\rm tr}\big(\partial_x\partial_\pi \phi(\pi)(x)\sigma\sigma\trans(x,a)\big)
\enq
and $\Mc^a\phi(\pi)$ $\in$ $L_{\pi\otimes\pi}^2(\R)$ is the function $\R^n\times\R^n$ $\rightarrow$ $\R$ defined by
\beq \label{defMc}
\Mc^a\phi(\pi)(x,x') &=& \frac{1}{2}{\rm tr}\big( \partial_\pi^2 \phi(\pi)(x,x')\sigma_{_W}(x,a)\sigma_{_W}\trans(x',a) \big).
\enq

The Bellman equation \reff{HJBpi} is validated by the following verification theorem, whose proof will rely on the above It\^o's formula.

\begin{Theorem} \label{theoverif}
(Verification Theorem)

\noindent  Let $w$ be a real-valued function in $\Cc_b^{1,2}([0,T]\times\Pc_{_2}(\R^n))$, i.e. $w$ is continuous on $[0,T]\times\Pc_{_2}(\R^n)$,
$w(t,.)$ $\in$ $\Cc_b^2(\Pc_{_2}(\R^n))$ for all $t$ $\in$ $[0,T]$, and $w(.,\pi)$ $\in$ $\Cc^1([0,T))$ for all $\pi$ $\in$ $\Pc_{_2}(\R^n)$, and satisfying
a quadratic growth condition as in \reff{vquadrapi},  together with a linear growth condition for its derivative:
\beq \label{Dwpilin}
|\partial_\pi w(t,\pi)(x)| & \leq & C(1 + |x| + \|\pi\|_{_2}), \;\;\; \forall (t,x,\pi) \in [0,T]\times\R^d\times\Pc_{_2}(\R^d),
\enq
for some positive constant $C$.
Suppose that $w$ is solution to the Bellman equation \reff{HJBpi}, and there exists for all
$(t,\pi)$ $\in$ $[0,T)\times\Pc_{_2}(\R^n)$ an element $\hat a(t,\pi)$ attaining the supremum in \reff{HJBpi} such that  the stochastic McKean-Vlasov SDE
\beqs
d\hat X_s &=&  b\big(\hat X_s, \hat a(s,\Lc(\hat X_s|W)) \big) ds +  \sigma\big(\hat X_s, \hat a(s,\Lc(\hat X_s|W)) \big)  dB_s,  \;\;\; t \leq s \leq T,
\enqs
admits a unique solution $\hat X^{t,\xi}$ given an initial condition $\xi$ $\in$ $L^2(\Fc^o;\R^n)$.
Then, $w$ $=$ $v$ and the feedback control in $\Ac^W$ defined by
\beq \label{defhata}
\alpha_s^* &=& \hat a(s,\Lc(\hat X_s^{t,\xi}|W)),
%\;  \mbox{ with } \;\;\;  \hat\rho_s^{t,\pi}(\varphi) \; = \;  \bar\E \big[ \varphi(\hat X_s^{t,\xi}) | \Fc_s^W \big],  \;\;\;  \pi = \Lc(\xi),
\; t \leq s \leq T,
\enq
%$\P^W_{_{\hat X_s^{t,\xi}}}$
is an optimal control for $v(t,\pi)$, i.e. $v(t,\pi)$ $=$ $J(t,\xi,\alpha^*)$ for $\pi$ $=$ $\Lc(\xi)$.
\end{Theorem}
{\bf Proof.}  Fix $(t,\pi = \Lc(\xi))$ $\in$ $[0,T)\times\Pc_{_2}(\R^n)$, and consider an arbitrary control $\alpha$ $\in$ $\Ac^W$ associated to
$X^{t,\xi,\alpha}$  solution to \reff{dynX}.  From \reff{estimX} under {\bf (H1)}, we have
\beqs
\bar\E \Big[ \int_t^T |b(X_s^{t,\xi,\alpha},\alpha_s)|^2 + |\sigma(X_s^{t,\xi,\alpha},\alpha_s)|^2 ds \Big] & < & \infty.
\enqs
Consider the controlled filter  process $\{\rho_s^{t,\pi,\alpha},t\leq s\leq T\}$ of the controlled process $X^{t,\xi,\alpha}$, which defines  an $\F^W$-optional process  valued in $\Pc_{_2}(\R^n)$, given  for any bounded measurable function $\varphi$ by
\beq \label{controlfilter}
\rho_s^{t,\pi,\alpha}[\varphi] &=& \bar\E \Big[ \varphi(X_s^{t,\xi,\alpha}) \big| \Fc_s^W \Big],  \;\;\; t \leq s \leq T, \; \pi = \Lc(\xi) \in \Pc_{_2}(\R^n), \; \alpha\in\Ac^W.
\enq
Notice that $X_s^{t,\xi,\alpha}$ is independent of the increments of $W_r-W_s$, for all $r\geq s$. Thus, $\rho_s^{t,\pi,\alpha}$ is equal to the conditional law of
$X_s^{t,\xi,\alpha}$ given $\Fc_\infty^W$, denoted by $\Lc(X_s^{t,\xi,\alpha}|W)$, and notice that  $\P(d\omega)$-a.s., $\rho_s^{t,\pi,\alpha}(\omega)$  is the law of $X_s^{t,\xi,\alpha}(.,\omega)$ on $(\bar\Omega^1,\bar\Fc^1,\bar\P^1)$.
From  \reff{estimX}, we have the estimate:
\beq \label{estimrhopi2}
\E\Big[ \sup_{t\leq s\leq T} |\rho_s^{t,\pi,\alpha}|^2 \Big] & \leq & C(1 + \|\pi\|_{_2}^2), \;\;\; \forall (t,\pi) \in [0,T]\times\Pc_{_2}(\R^n), \alpha \in \Ac^W,
\enq
for some positive constant $C$. One can apply It\^o's formula \reff{Ito} to $w(s,\rho_s^{t,\pi,\alpha})$  between $t$ and the $\F^W$-stopping time $\theta_T^n$ $=$
$\inf\{ s\geq t: \|\rho_s^{t,\pi,\alpha}\|_{_2} \geq n \}$,  and obtain:
\beqs
w(\theta_T^n,\rho_{\theta_T^n}^{t,\pi,\alpha}) &=& w(t,\pi) + \int_t^{\theta_T^n} \partial_t w(s,\rho_s^{t,\pi,\alpha})   \; + \;
\bar\E^1 \Big[ \partial_\pi w(s,\rho_s^{t,\pi,\alpha})(X_s).b(X_s,\alpha_s)  \\
& & \;\;\;\;\;\;\;\;  + \;  \frac{1}{2}{\rm tr}\big(\partial_x \partial_\pi w(s,\rho_s^{t,\pi,\alpha})(X_s) \sigma_s\sigma_s\trans(X_s,\alpha_s)\big) \Big]   \nonumber \\
& & \;\;\;\;\;\;\;  + \;  \bar\E^1 \Big[ \bar\E^{'1} \big[ \frac{1}{2}{\rm tr}\big(\partial_\pi^2 w(s,\rho_s^{t,\pi,\alpha})(X_s,X_s') \sigma_{_W}(X_s,\alpha_s)
\sigma_{_W}(X_s',\alpha_s)\trans\big) \big] \Big] ds  \nonumber \\
& & \;\;\; + \; \int_t^{\theta_T^n} \bar\E^1 \Big[ \partial_\pi  w(s,\rho_s^{t,\pi,\alpha})(X_s)\trans \sigma_{_W}(X_s,\alpha_s) \Big] dW_s,
\enqs
where we set for alleviating notations: $X_s$ $=$ $X_s^{t,\xi,\alpha}$, and $X_s'$ is a copy of $X_s$, $t\leq s \leq T$,
on another probability space  $(\bar\Omega^{'1}\times\Omega,\bar\Fc^{'1}\otimes\Fc,\bar\P^{'1}\otimes\P)$,
with  $(\bar\Omega^{'1},\bar\Fc^{'1},\bar\P^{'1})$ supporting $V'$ a copy of $V$,
(recall that $\alpha_s$ is $\Fc_s^{W}$-measurable and thus is a copy of itself).  By definition of  $\rho_s^{t,\pi,\alpha}$ and since  $\alpha_s$  is measurable with respect to
$\Fc_\infty^{W}$, the above relation can be written as
\begin{align}
w(\theta_T^n,\rho_{\theta_T^n}^{t,\pi,\alpha}) \ &= \ w(t,\pi) + \int_t^{\theta_T^n} \big\{ \partial_t w(s,\rho_s^{t,\pi,\alpha}) + \rho_s^{t,\pi,\alpha}\big[\Lc^{\alpha_s}w(s,\rho_s^{t,\pi,\alpha})\big]
\label{itow} \\
&\quad \ +  \rho_s^{t,\pi,\alpha} \otimes \rho_s^{t,\pi,\alpha} \big[ \Mc^{\alpha_s} w(s,\rho_s^{t,\pi,\alpha}) \big] \big\} ds   +
\int_t^{\theta_T^n} \rho_s^{t,\pi,\alpha} \big[  \partial_\pi  w(s,\rho_s^{t,\pi,\alpha})(.)\trans \sigma_{_W}(.,\alpha_s) \big] dW_s,  \nonumber
\end{align}
by definition of $\Lc^a$ and $\Mc^a$ in \reff{defLc}-\reff{defMc}.
Now, the integrand  of the stochastic integral with respect to $W$ in \reff{itow} satisfies:
\begin{align*}
\Big| \rho_s^{t,\pi,\alpha} \big[  \partial_\pi  w(s,\rho_s^{t,\pi,\alpha})(.)\trans \sigma_{_W}(.,\alpha_s) \big] \Big|^2 &\leq  \Big(  \int_{\R^n} \big| \partial_\pi w(s,\rho_s^{t,\pi,\alpha})(x)\trans\sigma_{_W}(x,\alpha_s) \big| \rho_s^{t,\pi,\alpha}(dx) \Big)^2 \\
& \leq \int_{\R^n} \big| \partial_\pi w(s,\rho_s^{t,\pi,\alpha})(x)\big|^2 \rho_s^{t,\pi,\alpha}(dx)   \int_{\R^n} \big| \sigma_{_W}(x,\alpha_s)  \big|^2 \rho_s^{t,\pi,\alpha}(dx)  \\
& \leq C(1 + n^2)^2
 <  \infty, \qquad  t \leq s \leq \theta_T^n,
\end{align*}
from Cauchy-Schwarz inequality, the linear growth condition of $\sigma_{_W}$ in {\bf (H1)}, the choice of  $\theta_T^n$, and condition \reff{Dwpilin}. Therefore, the stochastic integral in \reff{itow} vanishes in $\bar\P$-expectation, and
we get
\beq
\bar\E \big[ w(\theta_T^n,\rho_{\theta_T^n}^{t,\pi,\alpha}) \big] &=& w(t,\pi) \; + \; \bar\E \Big[ \int_t^{\theta_T^n}  \Dt{w}(s,\rho_s^{t,\pi,\alpha})  + \rho_s^{t,\pi,\alpha} \big[  \Lc^{\alpha_s} w(s,\rho_s^{t,\pi,\alpha}) \big]  \nonumber  \\
& & \hspace{3cm}  + \;   \rho_s^{t,\pi,\alpha}\otimes\rho_s^{t,\pi,\alpha}\big[ \Mc^{\alpha_s} w(s,\rho_s^{t,\pi,\alpha}) \big]  ds \Big]  \nonumber  \\
& \leq &  w(t,\pi) \; -  \; \bar\E \Big[ \int_t^{\theta_T^n}  \rho_s^{t,\pi,\alpha} [f(.,\alpha_s)]  ds \Big],  \label{Itoverif2}
\enq
since $w$ satisfies the Bellman equation \reff{HJBpi}.  By sending $n$ to infinity into \reff{Itoverif2}, and from the dominated convergence theorem (under the condition that $w$, $f$ satisfy a quadratic growth condition and recalling the estimation \reff{estimrhopi2}), we obtain:
\beqs
w(t,\pi) & \geq &   \bar\E\Big[ \int_t^T \rho_s^{t,\pi,\alpha}[f(.,\alpha_s)] ds  + \rho_T^{t,\pi,\alpha}[g]  \Big] \\
&=&  \bar\E \Big[ \int_t^T f(X_s^{t,\xi,\alpha},\alpha_s) ds + g(X_T^{t,\xi,\alpha}) \Big] \; = \; \Jc(t,\xi,\alpha)
\enqs
Since $\alpha$ is arbitrary in $\Ac$, this shows that $w(t,\pi)$ $\geq$ $\upv(\xi)$, hence $w$ $\geq$ $v$.

In the final step, let us apply the same It\^o's argument as in \reff{itow} with the control $\alpha^*$ in \reff{defhata}.
%and observe that $\P(d\omega)$-a.s., $\rho_s^{t,\pi,\alpha^*}(\omega)$  is the law of  $\hat X_s^{t,\xi}(.,\omega)$ on $(\bar\Omega^1,\bar\Fc^1,\bar\P^1)$.
Since $\hat a$ attains the supremum in the Bellman equation, we thus get after taking expectation:
\beqs
w(t,\pi) & = & \bar\E \Big[ \int_t^T   \rho_s^{t,\pi,\alpha^*} [ f(.,\alpha_s^*]   +  \rho_T^{t,\pi,\alpha^*}[g]  \Big] \; = \;
\bar\E \Big[ \int_t^T f(X_s^{t,\xi,\alpha^*},\alpha_s^*) ds + g(X_T^{t,\xi,\alpha^*}) \Big]
\enqs
which shows that  $w(t,\pi)$ $=$ $\Jc(t,\xi,\alpha^*)$ ($\leq$ $\upv(t,\xi)$ $=$ $v(t,\pi)$),
and therefore the required result: $v(t,\pi)$ $=$ $w(t,\pi)$ $=$ $\Jc(t,\xi,\alpha^*)$.
\ep

\vspace{2mm}

\begin{Remark}
{\rm (Bellman equation in the space of density functions)

\noindent Assume that the probability measure $\pi$ $\in$ $\Pc_{_2}(\R^n)$ of the random variable $\xi$ $\in$ $L^2(\Fc^o;\R^n)$ has a density $\eta$, i.e. $\pi(dx)$ $=$ $\eta(x)dx$,  with $\eta$ $\in$ $L^2(\R^n)$,
and  to alleviate notations, we set   $v(t,\eta)$ $=$ $v(t,\pi)$ for the value function  defined on $[0,T]\times L^2(\R^n)$ or on $[0,T]\times\Pc_{_2}(\R^n)$.
Suppose that $v(t,.)$ is Fr\'echet differentiable on $L^2(\R^n)$ with derivative $D_\eta v(t,\eta)$ $\in$ $L^2(\R^n)$ by Riesz theorem.   As shown in \cite{ben15}, this derivative is related to the derivative in the Wasserstein space by:
\beqs
\partial_\pi v(t,\pi)(x) \; = \;  D_x \big( D_\eta v(t,\eta)(x) \big), & \mbox{ and so }  & \partial_x \partial_\pi v(t,\pi)(x) \; = \;  D_x^2 \big( D_\eta v(t,\eta)(x) \big).
\enqs
Moreover, if $v(t,.)$ is twice Fr\'echet differentiable on $L^2(\R^n)$, its second derivative  is a self-adjoint operator on $L^2(\R^n)$, which can  be identified by Riesz theorem, as a symmetric function of two arguments, namely
$D_\eta^2 v(t,\eta)(x,x')$, and it is related to the second derivative in the Wasserstein space by:
\beqs
\partial_\pi^2 v(t,\pi)(x,x') &=& D_x D_{x'} \big( D^2_\eta v(t,\eta)(x,x') \big).
\enqs
Recalling the definition \reff{defLc} of the operator $\Lc^a$, we then have
\beqs
\pi\big[\Lc^a v(t,\pi) \big] &=& \int_{\R^n} \big[ D_x \big( D_\eta v(t,\eta)(x) \big).b(x,a) + \frac{1}{2}{\rm tr}\big(D_x^2 \big( D_\eta v(t,\eta)(x) \big)\sigma\sigma\trans(x,a)\big) \big] \eta(x) dx \\
&=&  \int_{\R^n}  D_\eta v(t,\eta)(x) (L^a)^* \eta(x) dx  \; = \; \langle D_\eta v(t,\eta),(L^a)^*\eta\rangle_{_{L^2(\R^n)}},
\enqs
by integration by parts and definition \reff{Ladjoint} of the adjoint operators $(L^a)^*$.  Similarly, recalling the definition \reff{defMc} of the operator $\Mc^a$, we have
\beqs
\pi\otimes\pi\big[\Mc^a v(t,\pi) \big] &=& \frac{1}{2} \int_{\R^n\times\R^n}  {\rm tr}\big[ D_x D_{x'} \big( D^2_\eta v(t,\eta)(x,x') \big) \sigma_{_W}(x,a)\sigma_{_W}\trans(x',a) \big] \eta(x) \eta(x') dx dx' \\
&=&   \frac{1}{2} \int_{\R^n\times\R^n}   D^2_\eta v(t,\eta)(x,x') (M^a)^*\eta(x)  (M^a)^*\eta(x') dx dx' \\
& = &  \frac{1}{2}\langle D_\eta^2 v(t,\eta)(M^a)^*\eta,(M^a)^*\eta\rangle_{_{L^2(\R^n)}},
\enqs
by integration by parts and definition \reff{Ladjoint} of the adjoint operators $(M^a)^*$.  Therefore, the Bellman equation \reff{HJBpi} is reduced to the Bellman equation in $[0,T]\times L^2(\R^n)$:
\begin{equation} \label{HJBden}
\left\{
\begin{array}{rcl}
\frac{\partial v(t,\eta)}{\partial t} + \Sup_{a\in A}
\Big[\langle D_\eta v(t,\eta),(L^a)^*\eta\rangle_{_{L^2(\R^n)}} & + & \frac{1}{2}\langle D_\eta^2 v(t,\eta)(M^a)^*\eta,(M^a)^*\eta\rangle_{_{L^2(\R^n)}} \\
& &  \hspace{1.5cm}   + \;\; \langle f(\cdot,a),\eta\rangle_{_{L^2(\R^n)}}  \Big]  \; = \;   0, \\
v(T,\eta)  &= & \langle g,\eta\rangle_{_{L^2(\R^n)}},
\end{array}
\right.
\end{equation}
which is the dynamic programming equation associated to the controlled DMZ equation \reff{DuncanMZcontrolintro}. A clear advantage of the Bellman approach in the Wasserstein space  compared to the classical Bellman approach in the space of density function, is that it is more general since we do not need to have a density for the controlled filter process, i.e. the  existence of a solution to the controlled SPDE \reff{DuncanMZcontrolintro}, which typically requires some uniform ellipticity condition on the diffusion matrix $\sigma$ (not needed here), see  \cite{par88} for related results.   Moreover, the viscosity property and the uniqueness result of the value function to \reff{HJBden} are  delicate issues as this Bellman equation in $L^2(\R^n)$ is fully nonlinear with unbounded first and second order terms, and needs to be studied in weighted  $L^2$ spaces, see \cite{lio89} and \cite{gozswi15} for some results.  Alternatively, and as investigated in the next section, the viscosity solutions approach for the Bellman equation \reff{HJBpi}, when lifted in the space $L^2(\Fc^o;\R^n)$, appears quite convenient and leads to rather general results under mild conditions, again without any uniform ellipticity assumption.
}
\ep
\end{Remark}

\subsection{Viscosity solutions}

In general, there are no smooth solutions to the HJB equation, and we shall thus consider in this section a notion of viscosity solutions for the Bellman equation \reff{HJBpi} in the Wasserstein space of probability measures $\Pc_{_2}(\R^n)$.  We adopt the approach in \cite{lio12}, and detailed in \cite{car12}, which consists, after the lifting identification between measures and random variables,  in working in the Hilbert space $L^2(\Fc^o;\R^n)$  instead of working in  the Wasserstein space $\Pc_{_2}(\R^n)$, in order to use the various   tools developed for viscosity solutions in separable Hilbert spaces, in particular in our context,  for second order Hamilton-Jacobi equations.
We shall assume that  the $\sigma$-algebra  $\Fc^o$  is  coun\-tably generated up to null sets, which ensures  that  the  Hilbert space $L^2(\Fc^o;\R^d)$ is separable, see \cite{doo94}, p. 92.  This is satisfied for example when $\Fc^o$ is the Borel $\sigma$-algebra of a canonical space $\Omega^o$ of continuous functions on $\R_+$,  in which case, $\Fc^o$ $=$ $\vee_{_{s\geq 0}} \Fc_s^{B^0}$, where $(\Fc_s^{B^0})$  is the canonical filtration on
$\Omega^o$, and it is then known that $\Fc^o$ is countably generated, see for instance Exercise 4.21 in Chapter 1 of \cite{revyor99}.

Let us then consider the ``lifted" Bellman equation in $(0,T]\times L^2(\Fc^o;\R^n)$:
\begin{equation} \label{HJBlift}
\left\{
\begin{array}{rcl}
- \Dt{\upv} + \Hc(\xi, D\upv(t,\xi),D^2 \upv(t,\xi)) &=& 0    \;\;\;\;\;\;\;  \mbox{ on } \; (0,T)\times L^2(\Fc^o;\R^n), \\
\upv(T,\xi) &=& \hat g(\xi) \; := \;  \E^o[g(\xi)],   \;\;\;  \xi  \in   L^2(\Fc^o;\R^n),
\end{array}
\right.
\end{equation}
where $\Hc$ $:$ $L^2(\Fc^o;\R^n)\times L^2(\Fc^o;\R^n)\times S(L^2(\Fc^o;\R^n))$ $\rightarrow$ $\R$ is defined by
\beq
\Hc(\xi,P,\Gamma) &=& - \sup_{a\in A} \Big\{ \E^o \big[ f(\xi,a)   +  P.b(\xi,a)
 \label{defHc} \\
& & \hspace{2cm}  + \;  \frac{1}{2} \Gamma (\sigma_{_V}(\xi,a)N).\sigma_{_V}(\xi,a)N  +   \frac{1}{2} \textup{tr}\big[\Gamma (\sigma_{_W}(\xi,a)).\sigma_{_W}(\xi,a) \big]  \big]  \Big\},  \nonumber
\enq
with $N$ $\in$ $L^2(\Fc^o,\R^m)$ of zero mean and unit variance, and independent of $\xi$.   Observe that when $u$ $\in$
$\Cc_b^2([0,T]\times\Pc_{_2}(\R^n)$ is a (classical) solution to the Bellman equation \reff{HJBpi}, then its lifted function $\upsilon$
is a classical solution to the lifted Bellman equation \reff{HJBlift}.

 We then  define viscosity solutions for the Bellman equation \reff{HJBpi} from viscosity solutions to  \reff{HJBlift}. As usual, we say  that a function $u$ (resp. $\upsilon$) is locally bounded in   $(0,T]\times \Pc_{_2}(\R^n)$ (resp.  on  $(0,T]\times L^2(\Fc^o;\R^n)$) if it is bounded on bounded subsets of $(0,T]\times \Pc_{_2}(\R^n)$ (resp. of  $(0,T]\times L^2(\Fc^o;\R^n)$), and we denote by
 $u^*$ (resp. $\upsilon^*$)  its upper   semicontinuous envelope, and by  $u_*$ (resp.  $\upsilon_*$) its lower semicontinuous envelope. More precisely, concerning for instance $\upsilon$ we have
\[
\upsilon_*(t,\xi) \ = \ \liminf_{\substack{(s,\eta)\rightarrow (t,\xi) \\ s < T}} \upsilon(s,\xi) \qquad \text{ and } \qquad \upsilon^*(t,\xi) \ = \ \limsup_{\substack{(s,\eta)\rightarrow (t,\xi) \\ s < T}} \upsilon(s,\xi)
\]
for all $(t,\xi) \in (0,T]\times L^2(\Fc^o;\R^n)$.

We denote by $\Cc_{b,loc}^{1,2}((0,T)\times L^2(\Fc^o;\R^n))$ the set of real-valued functions $\varphi$ on
 $(0,T)\times L^2(\Fc^o;\R^n)$ which are $\Cc^{1,2}((0,T)\times L^2(\Fc^o;\R^n))$ in the sense of Fr\'echet, and $\varphi$, together with its derivatives $\Dt{\varphi}$, $D\varphi$, $D^2\varphi$, is bounded on bounded subsets of $(0,T)\times L^2(\Fc^o;\R^n)$.

 \begin{Definition} \label{defvisco}
 We say that a locally bounded  function $u$ $:$ $[0,T]\times\Pc_{_2}(\R^n)$ $\rightarrow$ $\R$ is a viscosity (super, sub) solution to \reff{HJBpi} if its lifted function $\upsilon$ $:$ $[0,T]\times L^2(\Fc^o;\R^n)$ $\rightarrow$ $\R$
is a viscosity (super, sub) solution to the lifted Bellman equation \reff{HJBlift}, that is:

\noindent\textup{(i)} $\upsilon_*(T,.)$ $\geq$ $\hat g$, and for any test function $\varphi$ $\in$ $\Cc_{b,loc}^{1,2}((0,T)\times L^2(\Fc^o;\R^n))$  such that
$\upsilon_*-\varphi$ has a minimum at $(t_{_0},\xi_{_0})$ $\in$ $(0,T)\times L^2(\Fc^o;\R^n)$, one has
 \beqs
-  \Dt{\varphi}(t_{_0},\xi_{_0}) +   \Hc(\xi_{_0},D\varphi(t_{_0},\xi_{_0}),D^2\varphi(t_{_0},\xi_{_0})) & \geq & 0.
 \enqs
\textup{(ii)} $\upsilon^*(T,.)$ $\leq$ $\hat g$, and for any test function $\varphi$ $\in$ $\Cc_{b,loc}^{1,2}((0,T)\times L^2(\Fc^o;\R^n))$  such that
$\upsilon^*-\varphi$ has a maximum at $(t_{_0},\xi_{_0})$ $\in$ $(0,T)\times L^2(\Fc^o;\R^n)$, one has
 \beqs
-  \Dt{\varphi}(t_{_0},\xi_{_0}) +  \Hc(\xi_{_0},D\varphi(t_{_0},\xi_{_0}),D^2\varphi(t_{_0},\xi_{_0})) & \leq & 0.
 \enqs
 \end{Definition}

\begin{Remark}
{\rm
Notice that every test function in Definition \ref{defvisco} is also a test function in the sense of \cite{gozswi15} (see Definition 3.32 in \cite{gozswi15}). This will allow us to deduce a uniqueness result for viscosity solutions to the lifted Bellman equation \reff{HJBlift} from Theorem 3.50 in \cite{gozswi15}, see Theorem \ref{T:ExistUniq} below.
\ep}
\end{Remark}

The main result of this section is to prove the viscosity characterization  of the (inverse-lifted) value function $v$ in \reff{vpi} to the Bellman equation \reff{HJBlift}, hence equivalently the viscosity characterization of the value function $\upv$ in \reff{defv} to the  lifted Bellman equation \reff{HJBlift}.
In addition to the standing conditions {\bf (H1)}-{\bf (H2)}, we shall make a uniform continuity assumption on the  running gain function.

\vspace{2mm}

\hspace{-7mm} {\bf (H3)}  There exists   some modulus  $m_{f}$ (namely, $m_f\geq0$ and $m_{f}(\delta)$ $\rightarrow$ $0$ when $\delta$ goes to zero) such that for all  $x,x'$ $\in$ $\R^n$, $a$ $\in$ $A$,
\beqs
|f(x,a) - f(x',a)| & \leq & m_f\big(|x-x'|\big).
\enqs

\begin{Theorem}\label{T:ExistUniq}
The value function $\upv$ in \reff{defv} is the unique viscosity solution to the lifted Bellman equation \reff{HJBlift} satisfying a quadratic growth condition \reff{vquadra}. Moreover, $\upv$ is continuous on $(0,T)\times L^2(\Fc^o;\R^n)$.
\end{Theorem}
{\bf Proof.} \textbf{Step I.} {\it Viscosity property}. Let us first reformulate the randomized dynamic programming principle \reff{DPPrandom}  for the lifted value function $\upv$  on $(0,T]\times L^2(\Fc^o;\R^n)$.
For this, we take a copy $\tilde V$  of $V$ on  the probability space $(\Omega^o,\Fc^o,\P^o)$, and  given $(t,\xi,a)$ $\in$ $[0,T]\times L^2(\Fc^o;\R^n)\times A$, we consider on
$(\tilde\Omega=\Omega^o\times\Omega,\tilde\Fc=\Fc^o\otimes\Fc,\tilde\P =\P^o\otimes\P)$ the solution $\tilde X^{t,\xi,a}$, $t\leq s\leq T$, to the equation
\beqs
\tilde X_s^{t,\xi,a} &=& \xi + \int_t^s b(\tilde X_r^{t,\xi,a},I_r^{t,a}) dr + \int_t^s \sigma_{_V}(\tilde X_r^{t,\xi,a},I_r^{t,a}) d\tilde V_r + \int_t^s \sigma_{_W}(\tilde X_r^{t,\xi,a},I_r^{t,a}) dW_r, \; t \leq s\leq T.
\enqs
Thus,  $\tilde X^{t,\xi,a}$ is a copy of  $X^{t,\xi,a}$ on $(\tilde\Omega,\tilde\Fc,\tilde\P)$, and denoting by $\check X_s^{t,\xi,a}(\omega)$ $=$
$\tilde X_s^{t,\xi,a}(.,\omega)$, $t\leq s\leq T$, we see that the process $\{\check X_s^{t,\xi,a},t\leq s \leq T\}$ is $\F^{W,\mu,t}$-progressive, valued in $L^2(\Fc^o;\R^n)$, and
$\Lc(\check X_s^{t,\xi,a}(\omega))$ $=$ $\rho_s^{t,\pi,a}(\omega)$ for $\pi$ $=$ $\Lc(\xi)$.  Therefore,  the lifted value function on $[0,T]\times L^2(\Fc^o;\R^n)$  satisfies
$\upv(s,\check X_s^{t,\xi,a})$ $=$ $v(s,\rho_s^{t,\pi,a})$, $t\leq s \leq T$. %Let $\check f\colon L^2(\Fc^o;\R^n)\times A\rightarrow\R$ be given by $\check f(\xi,a):=f(\xi(\omega^o),a)$, for every $(\xi,a)\in L^2(\Fc^o;\R^n)\times A$. In a similar way we define $\check b$, $\check\sigma^V$, $\check\sigma^W$.
By noting that $\rho_s^{t,\pi,a}[f(.,\alpha_s)]$ $=$ $\E^o\big[ f(\tilde X_s^{t,\xi,a},\alpha_s)\big]$,
we  obtain from \reff{DPPrandom} the lifted DPP: for all $(t,\xi)$ $\in$ $[0,T]\times L^2(\Fc^o;\R^n)$,
\beq
\upv(t,\xi) &=& \sup_{\nu\in\Vc_{_{W,\mu,t}}} \sup_\theta \E^{\nu,t} \bigg[ \int_t^\theta \E^o\big[ f(\tilde X_s^{t,\xi,a},I^{t,a}_s)\big]  ds \; + \; \upv(\theta,\check X_\theta^{t,\xi,a}) \bigg] \notag \\
&=& \sup_{\nu\in\Vc_{_{W,\mu,t}}} \inf_\theta \E^{\nu,t} \bigg[ \int_t^\theta \E^o\big[ f(\tilde X_s^{t,\xi,a},I^{t,a}_s)\big]  ds \; + \; \upv(\theta,\check X_\theta^{t,\xi,a}) \bigg],
\label{DPPstronglift}
\enq
for any $a\in A$, with $\theta$ varying in the class of $\F^{W,\mu,t}$-stopping times valued in $[t,T]$.

\vspace{1mm}

\noindent (1) {\it  Supersolution property on $(0,T)\times L^2(\Fc^o;\R^n)$.} Let $(t_{_0},\xi_{_0})\in (0,T)\times L^2(\Fc^o;\R^n)$, and let $\varphi\in$ $\Cc_{b,loc}^{1,2}((0,T)\times L^2(\Fc^o;\R^n))$ be such that
\begin{equation}\label{min_subsol}	
0 \ = \ 	(\upv_*-\varphi)(t_{_0},\xi_{_0}) \ = \ \min_{(t,\xi) \in(0,T)\times L^2(\Fc^o;\R^n)}(\upv_*-\varphi)(t,\xi).
\end{equation}
By definition of $\upv_*(t_{_0},\xi_{_0})$, we have
$$
(t_{_m},\xi_{_m}) \ \overset{m\rightarrow\infty}{\longrightarrow} \ (t_{_0},\xi_{_0}) \qquad \textup{and}\qquad\upv(t_{_m},\xi_{_m}) \ \overset{m\rightarrow\infty}{\longrightarrow} \ \upv_*(t_{_0},\xi_{_0}),
$$
for some sequence $(t_{_m},\xi_{_m})_{_m}$ in $ (0,T)\times L^2(\Fc^o;\R^n)$. From \eqref{min_subsol} and the continuity of $\varphi$,
$$
\gamma_{_m} \ := \ \upv(t_{_m},\xi_{_m})- \varphi(t_{_m},\xi_{_m}) \ \overset{m\rightarrow\infty}{\longrightarrow} \ 0,
$$
Fix $a \in A$, and  let $\tau_{_m}\colon\Omega\rightarrow[t_{_m},T]\cup\{+\infty\}$ be the $\F^{W,\mu,t_{_m}}$-stopping time given by
\[
\tau_{_m}(\omega) \ := \ \inf\Big\{s \geq t_{_m} \colon \E^o\big[\big|\check X_s^{t_{_m},\xi_{_m},a}(\omega)-\xi_{_m}\big|^2\big]^{1/2} \ \geq \ \eta/2\Big\},
\]
with the convention $\inf\emptyset=+\infty$, where  $\eta$ is  some positive constant.
Let also $(h_{_m})_{_m}$ be a sequence of strictly positive real numbers such that
$$
h_{_m} \ \overset{m\rightarrow\infty}{\longrightarrow} \ 0 \qquad\textup{and}\qquad\frac{\gamma_{_m}}{h_{_m}} \ \overset{m\rightarrow\infty}{\longrightarrow} \ 0.
$$
By the dynamic programming principle \eqref{DPPstronglift}, we obtain
\begin{align*}
\upv(t_{_m},\xi_{_m}) \ \geq \ \E^{\nu,t_{_m}} \bigg[ \int_{t_{_m}}^{\theta_{_m}} \E^o\big[ f\big(\tilde X_s^{t_{_m},\xi_{_m},a},I^{t_{_m},a}_s\big)\big]  ds \; + \; \upv\big(\theta_{_m},\check X_{\theta_{_m}}^{t_{_m},\xi_{_m},a}\big) \bigg],
\end{align*}
for any $\nu\in\Vc_{_{W,\mu,t_{_m}}}$, where $\theta_{_m} := \tau_{_m} \wedge (t_{_m} + h_{_m}) \wedge T_{_{m,1}}$, and $T_{_{m,1}}$ denotes the first jump time of $I^{t_{_m},a}$. In particular, we take $\nu\equiv1$, so that $\E^{\nu,t_{_m}}$ is the expectation $\E$ with respect to $\P$. Equation \eqref{min_subsol}	implies that $\upv \ \geq \ \upv_* \ \geq \ \varphi$, thus
\begin{align}\label{ineq_1}
\varphi(t_{_m},\xi_{_m}) + \gamma_{_m} \ \geq \  \E \bigg[ \int_{t_{_m}}^{\theta_{_m}} \E^o\big[ f\big(\tilde X_s^{t_{_m},\xi_{_m},a},I^{t_{_m},a}_s\big)\big]  ds \; + \; \varphi\big(\theta_{_m},\check X_{\theta_{_m}}^{t_{_m},\xi_{_m},a}\big) \bigg].
\end{align}
On the other hand, applying It\^o's formula \eqref{Ito2} to $\varphi(\cdot,\check X_\cdot^{t_{_m},\xi_{_m},a})$ between $t_{_m}$ and $\theta_{_m}$, we obtain
\begin{align*}
\varphi(\theta_{_m},\check X_{\theta_{_m}}^{t_{_m},\xi_{_m},a}) \ &= \ \varphi(t_{_m}, \xi_{_m})+  \int_{t_{_m}}^{\theta_{_m}} 	\frac{\partial \varphi}{\partial t}(s,\check X_s^{t_{_m},\xi_{_m},a}) \, ds\\
&\quad \ + \int_{t_{_m}}^{\theta_{_m}} \E^o\big[ D\varphi(s,\check X_s^{t_{_m},\xi_{_m},a}).b_s(\tilde X_s^{t_{_m},\xi_{_m},a},I^{t_{_m},a}_s) \big] ds\\
&\quad \ + \int_{t_{_m}}^{\theta_{_m}} \E^o\Big[\frac{1}{2} D^2\varphi(s,\check X_s^{t_{_m},\xi_{_m},a})(\sigma_s^V(\tilde X_s^{t_{_m},\xi_{_m},a},I^{t_{_m},a}_s) N).\sigma_s^V(\tilde X_s^{t_{_m},\xi_{_m},a},I^{t_{_m},a}_s) N \Big] ds \\
&\quad \ + \int_{t_{_m}}^{\theta_{_m}} \E^o\Big[
\frac{1}{2} \textup{tr}\big[D^2\varphi(s,\check X_s^{t_{_m},\xi_{_m},a})(\sigma_s^W(\tilde X_s^{t_{_m},\xi_{_m},a},I^{t_{_m},a}_s)).\sigma_s^W(\tilde X_s^{t_{_m},\xi_{_m},a},I^{t_{_m},a}_s) \big]  \Big] ds \nonumber \\
&\quad \ + \; \int_{t_{_m}}^{\theta_{_m}} \E^o \big[ D\varphi(s,\check X_s^{t_{_m},\xi_{_m},a})\trans\sigma_s^W(\tilde X_s^{t_{_m},\xi_{_m},a},I^{t_{_m},a}_s) \big] dW_s,
\end{align*}
%where $N$ $\in$ $L^2(\Fc^o;\R^m)$ is independent of $(\tilde V,\xi_{_m})$  with zero mean and unit variance.
for any   Gaussian random variable $N$ $\in$ $L^2(\Fc^o;\R^m)$ independent of $(\tilde V,\xi_{_m})$,  %\textcolor{red}{($N$ should be $N_{_m}$?)}
with zero mean and unit variance.  Without loss of generality, we assume that $N$ is independent of $\xi_{_m}$ for every $m$.
 Therefore, taking expectation in \eqref{ineq_1}, and noting that the stochastic integral with respect to $W$ vanishes by the localization with the stopping time $\theta_{_m}$, we  get
\begin{align}\label{final_sub}
\frac{\gamma_{_m}}{h_{_m}} -	\E \bigg[\frac{1}{h_{_m}}\int_{t_{_m}}^{\theta_{_m}} \frac{\partial \varphi}{\partial t}(s,\check X_s^{t_{_m},\xi_{_m},a})\,ds\bigg] & \notag \\
+ \, \E \bigg[\frac{1}{h_{_m}} \int_{t_{_m}}^{\theta_{_m}} \E^o\bigg[  - D\varphi(s,\check X_s^{t_{_m},\xi_{_m},a}).b_s(\tilde X_s^{t_{_m},\xi_{_m},a},I^{t_{_m},a}_s) & \notag \\
- \, \frac{1}{2} D^2\varphi(s,\check X_s^{t_{_m},\xi_{_m},a})(\sigma_s^V(\tilde X_s^{t_{_m},\xi_{_m},a},I^{t_{_m},a}_s) N).\sigma_s^V(\tilde X_s^{t_{_m},\xi_{_m},a},I^{t_{_m},a}_s) N & \notag \\
- \, \frac{1}{2} \textup{tr}\big[D^2\varphi(s,\check X_s^{t_{_m},\xi_{_m},a})(\sigma_s^W(\tilde X_s^{t_{_m},\xi_{_m},a},I^{t_{_m},a}_s)).\sigma_s^W(\tilde X_s^{t_{_m},\xi_{_m},a},I^{t_{_m},a}_s) \big] & \notag \\
- \, f(\tilde X_s^{t_{_m},\xi_{_m},a},I^{t_{_m},a}_s)\bigg] ds \bigg] \ &\geq \ 0.
\end{align}
Now we notice that, $\P$-a.s., $I_r^{t_{_m},a} = a$ for $r \in [t_{_m},\theta_{_m}]$. Moreover, $\E\big[\E^o\big[\big|\check X_s^{t_{_m},\xi_{_m},a}-\xi_{_m}\big|^2\big]\big]\rightarrow0$, therefore, up to a subsequence, $\E^o\big[\big|\check X_s^{t_{_m},\xi_{_m},a}(\omega)-\xi_{_m}\big|^2\big]\rightarrow0$, $\P(d\omega)$-a.s., so that, $\P(d\omega)$-a.s., $\theta_{_m}(\omega)= t_{_m} + h_{_m}$ when $m$ is sufficiently large $(m\geq M(\omega))$. Thus, by the mean value theorem, the random variable inside the expectation $\E$ in \eqref{final_sub} converges $\P$-a.s. to
\begin{align*}
- \, \frac{\partial \varphi}{\partial t}(t_{_0},\xi_{_0}) -\E^o\Big[ D\varphi(t_{_0},\xi_{_0}).b_{t_{_0}}(\xi_{_0},a)   +  \frac{1}{2} D^2\varphi(t_{_0},\xi_{_0})(\sigma_{t_{_0}}^V(\xi_{_0},a) N).\sigma_{t_{_0}}^V(\xi_{_0},a) N \\
+ \, \frac{1}{2} \textup{tr}\big[D^2\varphi(t_{_0},\xi_{_0})(\sigma_{t_{_0}}^W(\xi_{_0},a)).\sigma_{t_{_0}}^W(\xi_{_0},a)\big] + f(\xi_{_0},a)\Big],
\end{align*}
when $m$ goes to infinity. Then, by the Lebesgue dominated convergence theorem, we obtain from \eqref{final_sub}
\begin{align*}
- \, \frac{\partial \varphi}{\partial t}(t_{_0},\xi_{_0}) -\E^o\Big[ D\varphi(t_{_0},\xi_{_0}).b_{t_{_0}}(\xi_{_0},a)   +  \frac{1}{2} D^2\varphi(t_{_0},\xi_{_0})(\sigma_{t_{_0}}^V(\xi_{_0},a) N).\sigma_{t_{_0}}^V(\xi_{_0},a) N & \\
+ \, \frac{1}{2} \textup{tr}\big[D^2\varphi(t_{_0},\xi_{_0})(\sigma_{t_{_0}}^W(\xi_{_0},a)).\sigma_{t_{_0}}^W(\xi_{_0},a)\big] + f(\xi_{_0},a)\Big] \ &\geq \ 0.
\end{align*}
The claim follows from the arbitrariness of $a \in A$.

\vspace{1mm}

\noindent (2) {\it Supersolution property: terminal condition.} Fix $\xi_{_0}\in L^2(\Fc^o;\R^n)$ and let us prove that $\upsilon_*(T,\xi_{_0})$ $\geq$ $\hat g(\xi_{_0})$. From the definition of $\upv_*(T,\xi_{_0})$, we have
$$
(t_{_m},\xi_{_m}) \ \overset{m\rightarrow\infty}{\longrightarrow} \ (T,\xi_{_0}) \qquad \textup{and}\qquad\upv(t_{_m},\xi_{_m}) \ \overset{m\rightarrow\infty}{\longrightarrow} \ \upv_*(T,\xi_{_0}),
$$
for some sequence $(t_{_m},\xi_{_m})_{_m}$ in $ (0,T)\times L^2(\Fc^o;\R^n)$. Fix $a\in A$. By the dynamic programming principle \eqref{DPPstronglift},
\begin{align}
\upv(t_{_m},\xi_{_m}) \ &\geq \ \sup_{\nu\in\Vc_{_{W,\mu,t_{_m}}}} \E^{\nu,t_{_m}} \bigg[ \int_{t_{_m}}^T \E^o\big[ f(\tilde X_s^{t_{_m},\xi_{_m},a},I^{t_{_m},a}_s)\big]  ds \; + \; \E^o\big[g(\tilde X_T^{t_{_m},\xi_{_m},a})\big] \bigg] \notag \\
&\geq \ \E \bigg[ \int_{t_{_m}}^T \E^o\big[ f\big(\tilde X_s^{t_{_m},\xi_{_m},a},I^{t_{_m},a}_s\big)\big]  ds \; + \; \E^o\big[g\big(\tilde X_T^{t_{_m},\xi_{_m},a}\big)\big] \bigg], \label{Proof_terminal_condition}
\end{align}
where the inequality follows taking $\nu\equiv1$. Now, from the quadratic growth condition {\bf (H2)} and estimate \eqref{estimXrandom}, we obtain
\[
\E \bigg[ \int_{t_{_m}}^T \E^o\big[ f\big(\tilde X_s^{t_{_m},\xi_{_m},a},I^{t_{_m},a}_s\big)\big]  ds\bigg] \ \overset{m\rightarrow\infty}{\longrightarrow} \ 0.
\]
On the other hand, from the continuity of $g$ and the fact that $\E\big[\E^o\big[\big|\tilde X_T^{t_{_m},\xi_{_m},a}-\xi_{_0}\big|^2\big]\big]\rightarrow0$, it follows that $g(\tilde X_T^{t_{_m},\xi_{_m},a})$ converges in $\P\otimes\P^o$-probability to $g(\xi_{_0})$. By the quadratic growth condition {\bf (H2)} and estimate \eqref{estimXrandom}, we obtain from \eqref{Proof_terminal_condition}, using the Lebesgue dominated convergence theorem,
\[
\upv_*(T,\xi_{_0}) \ = \ \lim_{m\rightarrow\infty} \upv(t_{_m},\xi_{_m}) \ \geq \ \lim_{m\rightarrow\infty} \E\big[\E^o\big[g\big(\tilde X_T^{t_{_m},\xi_{_m},a}\big)\big]\big] \ = \ \lim_{m\rightarrow\infty} \E \big[ \hat g\big(\check X_T^{t_{_m},\xi_{_m},a}\big) \big] \ = \ \hat g(\xi_{_0}).
\]

\vspace{1mm}

\noindent (3) {\it Subsolution property on $(0,T)\times L^2(\Fc^o;\R^n)$.} Let $(t_{_0},\xi_{_0})\in (0,T)\times L^2(\Fc^o;\R^n)$ and $\varphi\in$ $\Cc_{b,loc}^{1,2}((0,T)\times L^2(\Fc^o;\R^n))$ be such that
	\begin{equation}\label{max_subsol}	
	0 \ = \ 	(\upv^*-\varphi)(t_{_0},\xi_{_0}) \ = \ \max_{(t,\xi) \in(0,T)\times L^2(\Fc^o;\R^n)}(\upv^*-\varphi)(t,\xi).
	\end{equation}
Without loss of generality, we can assume that $(t_{_0},\xi_{_0})$ is a strict maximum of $\upv^*-\varphi$. Then, for any $\eta>0$, we can define $\beta(\eta)$ as follows:
$$
\max_{\partial B((t_{_0},\xi_{_0});\eta)}(\upv^*-\varphi) \ =: \ -\beta(\eta) \ < \ 0,
$$
where
\[
\partial B((t_{_0},\xi_{_0});\eta) \ = \ \Big\{(t,\xi)\in(0,T)\times L^2(\Fc^o;\R^n)\colon\max\big(|t-t_{_0}|,\E^o\big[|\xi-\xi_{_0}|^2\big]^{1/2}\big) \ = \ \eta\Big\}.
\]
We will proceed by contradiction, assuming that
\begin{equation}\label{rcl}
- \Dt{\varphi}(t_{_0},\xi_{_0}) + \Hc(\xi_{_0}, D\varphi(t_{_0},\xi_{_0}),D^2 \varphi(t_{_0},\xi_{_0})) \ > \ 0.
\end{equation}
Then by the continuity of the function   $\Hc$ on the interior of its domain, there exists $\eta >0$, with corresponding $\beta(\eta)>0$, and $\varepsilon \in (0,\,\beta(\eta)]$ such that
\begin{equation}\label{rcl1}
- \Dt{\varphi}(t,\xi) + \Hc(\xi, D\varphi(t,\xi),D^2 \varphi(t,\xi)) \ \geq \ \varepsilon,
\end{equation}
for all $(t,\xi) \in B((t_{_0}, \xi_{_0});\eta)= \{(t,\xi) \in (0,T)\times L^2(\Fc^o;\R^n)\colon|t_0-t| < \eta,\, \E^o\big[|\xi-\xi_{_0}|^2\big]^{1/2}< \eta\}$.

By definition of $\upv^*(t_{_0},\xi_{_0})$, there exists a sequence $(t_{_m},\xi_{_m})_{_m}$ taking values in $B((t_{_0}, \xi_{_0});\eta)$ such that
$$
(t_{_m},\xi_{_m}) \ \overset{m\rightarrow\infty}{\longrightarrow} \ (t_{_0},\xi_{_0}) \qquad \textup{and } \qquad \upv(t_{_m},\xi_{_m}) \ \overset{m\rightarrow\infty}{\longrightarrow} \ \upv^*(t_{_0},\xi_{_0}).
$$
By the continuity of $\varphi$ and using \eqref{max_subsol}, we also obtain
$$
\gamma_{_m} \ := \ \upv(t_{_m},\xi_{_m})- \varphi(t_{_m},\xi_{_m}) \ \overset{m\rightarrow\infty}{\longrightarrow} \ 0.
$$
Fix $a \in A$. Using the dynamic programming principle \eqref{DPPstronglift}, we find
\begin{align*}
\upv(t_{_m},\xi_{_m}) \ \leq \ \sup_{\nu\in\Vc_{_{W,\mu,t_{_m}}}} \E^{\nu,t_{_m}} \bigg[ \int_{t_{_m}}^{\theta_{_m}} \E^o\big[ f\big(\tilde X_s^{t_{_m},\xi_{_m},a},I^{t_{_m},a}_s\big)\big]  ds \; + \; \upv\big(\theta_{_m},\check X_{\theta_{_m}}^{t_{_m},\xi_{_m},a}\big) \bigg],
\end{align*}
with $\theta_{_m} := \tau_{_m}   \wedge T$, where $\tau_{_m}:=\inf\{s \geq t_{_m}\colon(s,\check X_s^{t_{_m},\xi_{_m},a}) \notin B((t_{_0},\xi_{_0});\eta/2)\}$. In particular, recalling \eqref{max_subsol}, there exists  $\nu_{_m} \in \Vc_{_{W,\mu,t_{_m}}}$ such that
\begin{align*}
&\varphi(t_{_m},\xi_{_m}) + \gamma_{_m} \\
&\leq \ \E^{\nu_{_m},t_{_m}} \bigg[ \int_{t_{_m}}^{\theta_{_m}} \E^o\big[ f\big(\tilde X_s^{t_{_m},\xi_{_m},a},I^{t_{_m},a}_s\big)\big]  ds \; + \; \varphi\big(\theta_{_m},\check X_{\theta_{_m}}^{t_{_m},\xi_{_m},a}\big) -\beta(\eta)\,1_{\{\tau_{_m} < T\}} \bigg] + \frac{\varepsilon}{2}(T - t_{_m}).
\end{align*}
Applying It\^o's formula \eqref{Ito2} to $\varphi(\cdot,\check X_\cdot^{t_{_m},\xi_{_m},a})$ between $t_{_m}$ and $\theta_{_m}$, we get
\begin{align}\label{gammaeq}
\gamma_{_m} - \frac{\varepsilon}{2}(T - t_{_m})-	\E^{\nu_{_m},t_{_m}} \bigg[ \int_{t_{_m}}^{\theta_{_m}}\frac{\partial \varphi}{\partial t}(s,\check X_s^{t_{_m},\xi_{_m},a})\,ds\bigg] + \beta(\eta)\,\E^{\nu_{_m},t_{_m}}[1_{\{\tau_{_m} < T\}}\,(\theta_{_m} -t_{_m})] & \nonumber\\
-\E^{\nu_{_m},t_{_m}} \bigg[ \int_{t_{_m}}^{\theta_{_m}} \E^o\bigg[ f(\tilde X_s^{t_{_m},\xi_{_m},a},I^{t_{_m},a}_s)   +D\varphi(s,\check X_s^{t_{_m},\xi_{_m},a}).b_s(\tilde X_s^{t_{_m},\xi_{_m},a},I^{t_{_m},a}_s) & \nonumber \\
+\frac{1}{2} D^2\varphi(s,\check X_s^{t_{_m},\xi_{_m},a})(\sigma_s^V(\tilde X_s^{t_{_m},\xi_{_m},a},I^{t_{_m},a}_s) N).\sigma_s^V(\tilde X_s^{t_{_m},\xi_{_m},a},I^{t_{_m},a}_s) N & \nonumber \\
+\frac{1}{2} \textup{tr}\big[D^2\varphi(s,\check X_s^{t_{_m},\xi_{_m},a})(\sigma_s^W(\tilde X_s^{t_{_m},\xi_{_m},a},I^{t_{_m},a}_s)).\sigma_s^W(\tilde X_s^{t_{_m},\xi_{_m},a},I^{t_{_m},a}_s)  \big] \bigg] ds \bigg] & \nonumber \\
-\E^{\nu_{_m},t_{_m}} \bigg[\int_{t_{_m}}^{\theta_{_m}} \E^o \bigg[ D\varphi(s,\check X_s^{t_{_m},\xi_{_m},a})\trans\sigma_s^W(\tilde X_s^{t_{_m},\xi_{_m},a},I^{t_{_m},a}_s) \bigg] dW_s\bigg] \
\leq \ 0.&
\end{align}
Denote
\begin{align*}
\Lc^b(\xi,P,\Gamma) \ = \ \E^o \big[ f(\xi,b)   +  P.b(\xi,b) + \tfrac{1}{2} \Gamma (\sigma_{_V}(\xi,b)N).\sigma_{_V}(\xi,b)N  +   \tfrac{1}{2} \textup{tr}[\Gamma (\sigma_{_W}(\xi,b)).\sigma_{_W}(\xi,b) ]\big],
\end{align*}
for every $(\xi,P,\Gamma,b)\in L^2(\Fc^o;\R^n)\times L^2(\Fc^o;\R^n)\times S(L^2(\Fc^o;\R^n))\times A$,
for any   Gaussian random variable $N$ $\in$ $L^2(\Fc^o;\R^m)$ independent of $(\tilde V,\xi)$,
with zero mean and unit variance.
%with $N$ $\in$ $L^2(\Fc^o,\R^n)$ of zero mean and unit variance, and independent of $\xi$.
 Notice that
\[
\Hc(\xi,P,\Gamma) \ = \ - \, \sup_{b\in A}\Lc^b(\xi,P,\Gamma).
\]
%Then
%\begin{align*}
%&\Lc^{I_s^{t_{_m},a}}\big(s,D\varphi(s,\check X_s^{t_{_m},\xi_{_m},a}),D^2\varphi(s,\check X_s^{t_{_m},\xi_{_m},a})\big)\\
%&= \ \E^o \big[ f(\check X_s^{t_{_m},\xi_{_m},a},I_s^{t_{_m},a})   +  D\varphi(\check X_s^{t_{_m},\xi_{_m},a},I_s^{t_{_m},a}) .\check b_s(\check X_s^{t_{_m},\xi_{_m},a}, I_s^{t_{_m},a})   \\
%&\quad \ + \tfrac{1}{2} D^2\varphi(\check X_s^{t_{_m},\xi_{_m},a},I_s^{t_{_m},a})  (\sigma_s^{V}(\check X_s^{t_{_m},\xi_{_m},a},I_s^{t_{_m},a})N).\sigma_s^{V}(\check X_s^{t_{_m},\xi_{_m},a},I_s^{t_{_m},a})N \\
%&\quad \ +   \tfrac{1}{2} D^2\varphi(\check X_s^{t_{_m},\xi_{_m},a},I_s^{t_{_m},a})  (\sigma_s^{_W}(\check X_s^{t_{_m},\xi_{_m},a},I_s^{t_{_m},a})).\sigma_s^{_W}(\check X_s^{t_{_m},\xi_{_m},a},I_s^{t_{_m},a}) \big].
%\end{align*}
%where $N$ belongs to $L^2(\Fc^o,\R^n)$, has zero mean and unit variance, and is independent of $(\tilde V,\textcolor{red}{\xi_{_m}})$.
Therefore, inequality \eqref{gammaeq} reads
\begin{align}\label{lastineq}
&	\E^{\nu_{_m},t_{_m}} \bigg[ \int_{t_{_m}}^{\theta_{_m}}\left(-\frac{\partial \varphi}{\partial t}(s,\check X_s^{t_{_m},\xi_{_m},a})- \Lc^{I^{t_{_m},a}_s}\big(s,D\varphi(s,\check X_s^{t_{_m},\xi_{_m},a}),D^2\varphi(s,\check X_s^{t_{_m},\xi_{_m},a})\big)\right) \, ds \bigg] \nonumber\\
&+\gamma_{_m} - \frac{\varepsilon}{2}(T - t_{_m}) + \beta(\eta)\,\E^{\nu_{_m},t_{_m}}\big[1_{\{\tau_{_m} < T\}}\,(\theta_{_m} -t_{_m})\big] \ \leq \ 0.
\end{align}
Notice that
\begin{align*}
&-\frac{\partial \varphi}{\partial t}(s,\check X_s^{t_{_m},\xi_{_m},a}) - \Lc^{I^{t_{_m},a}_s}\big(s,D\varphi(s,\check X_s^{t_{_m},\xi_{_m},a}),D^2\varphi(s,\check X_s^{t_{_m},\xi_{_m},a})\big)\\
&\geq \ -\frac{\partial \varphi}{\partial t}(s,\check X_s^{t_{_m},\xi_{_m},a})-\sup_{b \in A} \Lc^{b}\big(s,D\varphi(s,\check X_s^{t_{_m},\xi_{_m},a}),D^2\varphi(s,\check X_s^{t_{_m},\xi_{_m},a})\big)\\
&= \ -\frac{\partial \varphi}{\partial t}(s,\check X_s^{t_{_m},\xi_{_m},a})+\Hc\big(s,D\varphi(s,\check X_s^{t_{_m},\xi_{_m},a}),D^2\varphi(s,\check X_s^{t_{_m},\xi_{_m},a})\big) \ \geq \ \varepsilon,
\end{align*}
where the last inequality follows from \eqref{rcl1}.  Therefore \eqref{lastineq} yields
\begin{align*}
\gamma_{_m} - \frac{\varepsilon}{2}(T - t_{_m}) \ &\leq \ - \varepsilon\,\E^{\nu_{_m},t_{_m}}[\theta_{_m}-t_{_m}]- \beta(\eta)\,\E^{\nu_{_m},t_{_m}}\big[1_{\{\tau_{_m} < T\}}\,(\theta_{_m} -t_{_m})\big] \\
&\leq \ - \varepsilon\,\E^{\nu_{_m},t_{_m}}[\theta_{_m}-t_{_m}]- \varepsilon\,\E^{\nu_{_m},t_{_m}}\big[1_{\{\tau_{_m} < T\}}\,(\theta_{_m} -t_{_m})\big] \\
&= \ - \varepsilon\,\E^{\nu_{_m},t_{_m}}\big[1_{\{\theta_{_m}=T\}}\,(\theta_{_m} -t_{_m})\big] \ = \ - \varepsilon\,(T - t_{_m}).
\end{align*}
Sending $m\rightarrow\infty$ we get the contradiction $-\frac{\varepsilon}{2}(T-t_{_0})\leq -\varepsilon(T-t_{_0})$.

\vspace{1mm}

\noindent (4) {\it Subsolution property: terminal condition.} Fix $\xi_{_0}\in L^2(\Fc^o;\R^n)$ and let us prove that $\upsilon^*(T,\xi_{_0})$ $\leq$ $\hat g(\xi_{_0})$. From the definition of $\upv^*(T,\xi_{_0})$, we have
$$
(t_{_m},\xi_{_m}) \ \overset{m\rightarrow\infty}{\longrightarrow} \ (T,\xi_{_0}) \qquad \textup{and}\qquad\upv(t_{_m},\xi_{_m}) \ \overset{m\rightarrow\infty}{\longrightarrow} \ \upv^*(T,\xi_{_0}),
$$
for some sequence $(t_{_m},\xi_{_m})_{_m}$ in $ (0,T)\times L^2(\Fc^o;\R^n)$. Fix $a\in A$. By the dynamic programming principle \eqref{DPPstronglift},
\begin{align*}
\upv(t_{_m},\xi_{_m}) \ \leq \ \sup_{\nu\in\Vc_{_{W,\mu,t_{_m}}}} \E^{\nu,t_{_m}} \bigg[ \int_{t_{_m}}^T \E^o\big[ f(\tilde X_s^{t_{_m},\xi_{_m},a},I^{t_{_m},a}_s)\big]  ds \; + \; \E^o\big[g(\tilde X_T^{t_{_m},\xi_{_m},a})\big] \bigg].
\end{align*}
In particular, given $\eps>0$, there exists $\nu_{_m} \in \Vc_{_{W,\mu,t_{_m}}}$ such that
\begin{align}
\upv(t_{_m},\xi_{_m}) \ &\leq \ \E^{\nu_{_m},t_{_m}} \bigg[ \int_{t_{_m}}^T \E^o\big[ f\big(\tilde X_s^{t_{_m},\xi_{_m},a},I^{t_{_m},a}_s\big)\big]  ds \; + \; \E^o\big[g\big(\tilde X_T^{t_{_m},\xi_{_m},a}\big)\big] \bigg] + \eps. \label{Proof_terminal_condition2}
\end{align}
From the quadratic growth condition {\bf (H2)} and estimate \eqref{EstimateX_nu}, we obtain
\[
\E^{\nu_{_m},t_{_m}} \bigg[ \int_{t_{_m}}^T \E^o\big[ f\big(\tilde X_s^{t_{_m},\xi_{_m},a},I^{t_{_m},a}_s\big)\big]  ds\bigg] \ = \ \int_{t_{_m}}^T \bar\E^{\nu_{_m},t_{_m}}\big[ f\big(\tilde X_s^{t_{_m},\xi_{_m},a},I^{t_{_m},a}_s\big)\big]  ds \ \overset{m\rightarrow\infty}{\longrightarrow} \ 0.
\]
On the other hand, by standard estimate on $\tilde X^{t,\xi,a}$, we notice  that
\beqs
\sup_{\nu\in\Vc_{_{W,\mu,t_{_m}}}}\E^{\nu,t_{_m}}\big[\E^o\big[\big|\tilde X_T^{t_{_m},\xi_{_m},a}-\xi_{_0}\big|^2\big]\big]
& \overset{m\rightarrow\infty}{\longrightarrow}  & 0,
\enqs
which implies from the quadratic growth condition {\bf (H2)} and the continuity of $g$ that
\[
\sup_{\nu\in\Vc_{_{W,\mu,t_{_m}}}}\E^{\nu,t_{_m}}\big[\E^o\big[\big|g\big(\tilde X_T^{t_{_m},\xi_{_m},a}\big)-g(\xi_{_0})\big|\big]\big] \ \overset{m\rightarrow\infty}{\longrightarrow} \ 0.
\]
Hence, letting $m\rightarrow\infty$ in \eqref{Proof_terminal_condition2},
\begin{align*}
\upv^*(T,\xi_{_0}) \ \leq \ \E^o\big[g(\xi_{_0})\big] + \eps \ = \ \hat g(\xi_{_0}) + \eps.
\end{align*}
The claim follows from the arbitrariness of $\eps$.

\vspace{1mm}

\noindent \textbf{Step II.} {\it  Uniqueness and continuity of $\upv$.} In view of our definition of viscosity solution, we have to show a comparison principle for viscosity solutions to the lifted Bellman equation \reff{HJBlift}.
We use the comparison principle proved in Theorem 3.50 in \cite{gozswi15} and check that the hypotheses of this theorem are satisfied in our context for the lifted Hamiltonian $\Hc$ defined in \reff{defHc}.
Notice that the Bellman equation \reff{HJBlift} is a bounded equation in the terminology  of  \cite{gozswi15} (see their section  3.3.1) meaning that there is no linear dissipative operator on $L^2(\Fc^o;\R^n)$ in the equation.
Therefore, the notion of $B$-continuity reduces to the standard notion of continuity in $L^2(\Fc^o;\R^n)$ since one can take for $B$ the identity operator.
Their Hypothesis 3.44 follows from the uniform continuity  of $b$, $\sigma$, and $f$ in {\bf (H1)} and {\bf (H3)}.  Hypothesis 3.45  is immediately satisfied since there is no discount factor in our equation, i.e.
$\Hc$ does not depend on $\upv$ but only on its derivatives.  The monotonicity condition in $\Gamma$ $\in$ $S(L^2(\Fc^o;\R^n))$ of $\Hc$ in Hypothesis 3.46  is clearly satisfied.  Hypothesis 3.47  holds directly  when dealing with  bounded equations.   Hypothesis 3.48  is obtained from the Lipschitz condition of $b,\sigma$ in {\bf (H1)}, and the uniform continuity condition on $f$ in {\bf (H3)}, while Hypothesis 3.49  follows from the growth condition of
$\sigma$ in  {\bf (H1)}.

We can then apply Theorem 3.50  in  \cite{gozswi15} and deduce that given any other viscosity solution $\text{\Large$u$}$ to the lifted Bellman equation \eqref{HJBlift} we have:
\[
\upv_* \ \geq \ \text{\Large$u$}^* \qquad \text{ and } \qquad \text{\Large$u$}_* \ \geq \ \upv^*,
\]
on $(0,T]\times L^2(\Fc^o;\R^n)$. Since $\upv^*\geq\upv_*$ and $\text{\Large$u$}^*\geq\text{\Large$u$}_*$, we get $\upv_*=\text{\Large$u$}^*=\text{\Large$u$}_*=\upv^*$ on $(0,T]\times L^2(\Fc^o;\R^n)$. We conclude that $\upv=\upv^*=\upv_*$ and $\text{\Large$u$}=\text{\Large$u$}^*=\text{\Large$u$}_*$ are continuous on $(0,T)\times L^2(\Fc^o;\R^n)$ and $\upv=\text{\Large$u$}$ on $(0,T)\times L^2(\Fc^o;\R^n)$.
\ep

\section{Application: a partially observed non Gaussian LQ model} \label{secLQ}

\setcounter{equation}{0} \setcounter{Assumption}{0}
\setcounter{Theorem}{0} \setcounter{Proposition}{0}
\setcounter{Corollary}{0} \setcounter{Lemma}{0}
\setcounter{Definition}{0} \setcounter{Remark}{0}

We consider a model where the observation process  $(O_t)_t$ valued in $\R^d$ is governed by a diffusion process driven by the $d$-dimensional Brownian motion $W$:
\beqs
dO_t &=& \kappa(O_t) dt + \gamma(O_t) dW_t, \;\;\; t \geq 0, \;  O_0  \in \R^d.
\enqs
Under suitable standard assumptions on $\kappa$ and the invertible matrix coefficient $\gamma$, there exists a unique strong solution to the above SDE for $O$, and the observation filtration $\F^O$  generated by $O$ is equal to
$\F^W$, the natural filtration of $W$.  We then consider a  model \reff{dynX} for the controlled signal  process $X$ $=$ $X^\alpha$ governed by the $m$-dimensional (unobserved) Brownian motion $V$ $=$ $(V^1,\ldots,V^m)$, and
the $d$-dimensional Brownian motion $W$ $=$ $(W^{1},\ldots,W^{d})$, with linear dynamics:
\beq
dX_t &=& (b_0 + BX_t + C\alpha_t) dt + \sum_{i=1}^m (\gamma_{_V}^i +  D_{_V}^i X_t + F_{_V}^i \alpha_t) dV_t^i  \nonumber \\
& & \;\;\;+ \;     \sum_{j=1}^d (\gamma_{_W}^j +  D_{_W}^j X_t + F_{_W}^j \alpha_t) dW_t^j, \;\;\; 0 \leq t \leq T, \; X_0 = x_0 \in \R^n,  \label{dynXlin}
\enq
and control $\alpha$ $\in$ $\Ac^W$,  valued in $A$ $=$ $\R^q$, and progressively measurable with respect to $\F^O$ $=$ $\F^W$.
%coefficients $b(x,a)$, $\sigma_{_V}(x,a)$ $=$ $(\sigma_{_V}^1,\ldots,\sigma_{_V}^m)(x,a)$,
%$\sigma_{_W}(x,a)$ $=$ $(\sigma_{_W}^1,\ldots,\sigma_{_W}^d)(x,a)$  from $\R^n\times A$ (with $A$ $=$ $\R^q$) into $\R^n$ (resp. $\R^{n\times m}$,
%and $\R^{n\times d}$) given by
%\begin{equation} \label{linX}
%\left\{
%\begin{array}{ccl}
%b(x,a)  & =  &  b_0 +   Bx + Ca, \\
%\sigma_{_V}^i(x,a) &=& \gamma_{_V}^i +    D_{_V}^i x + F_{_V}^i a, \;\;\; i =1,\ldots,m \\
%\sigma_{_W}^j(x,a) &=&  \gamma_{_W}^i + D_{_W}^j x + F_{_W}^j a, \;\;\; j =1,\ldots,d.
%\end{array}
%\right.
%\end{equation}
Here, $b_0$, $\gamma_{_V}^i$, $\gamma_{_W}^j$ are constant vectors in $\R^n$, $B$, $D_{_V}^i$, $D_{_W}^j$, are constant matrices in $\R^{n\times n}$,  $C$, $F_{_V}^i$, $F_{_W}^j$ are constant matrices in $\R^{n\times q}$.
Notice that due to the affine terms in the diffusion coefficients, the pair of processes signal/observation $(X,W)$ is not Gaussian (even in the absence of control), hence we do not have a finite dimensional representation of the filter process (law of $X$ given $W$) as in the Kalman-Bucy filter.
The cost functional to be minimize over $\alpha$ $\in$ $\Ac^W$ is
\beq \label{costquadra}
J(\alpha) &=& \E \Big[ \int_0^T  \big(X_t\trans Q X_t + \alpha_t \trans N \alpha_t \big)dt + X_T\trans P X_T \Big]
\;\;\; \rightarrow  \;  \upv_0 = \inf_{\alpha\in\Ac^W} J(\alpha),
\enq
%The cost functions $f(x,a)$ and $g(x)$ from $\R^n\times\R^q$ (resp. $\R^n$) into  $\R$ are in the quadratic form
%\beq \label{quacost}
%f(x,a) \; = \;  x\trans Q x + a\trans N a, \;\;\; & & g(x) \; = \; x\trans P x,
%\enq
where $Q,P$ are  constant valued in $\S^n$, the set of symmetric matrices in $\R^{n\times n}$, and $N$ is constant valued in $\S^q$.
We shall make the following assumptions  on the  coefficients of the model:

\vspace{2mm}

\noindent {\bf (C1)} \hspace{1cm} $Q$,  $P$,  $N$ are nonnegative a.s.

\vspace{1mm}

\noindent {\bf (C2)}  \hspace{1cm} One of the two following conditions holds:
\begin{itemize}
\item[(i)] $N$ is uniformly positive definite  i.e. $N$ $\geq$ $\delta I_q$, for some $\delta$ $>$  $0$,
\item[(ii)] $P$ or $Q$ is uniformly positive definite, and $F_{_V}^i$ is uniformly nondegenerate, i.e. $|F_{_V}^i|$ $\geq$ $\delta_i$,
for some $i$ $=$ $1,\ldots,m$, and $\delta_i$ $>$ $0$.
\end{itemize}

We cannot use the classical results in \cite{Be} for the linear quadratic (LQ) Gaussian partial observation  framework, and instead we show how one can derive explicit solutions by solving the Bellman equation \reff{HJBpi} in our
non Gaussian LQ model \reff{dynXlin}-\reff{costquadra}.  Let us introduce some notations. For any $\pi$ $\in$ $\Pc_{_2}(\R^n)$, we denote by
\beqs
\bar\pi &:=& \int_{\R^n} x \pi(dx),
\enqs
and define the functions on $\Pc_{_2}(\R^n)\times\S^n$ and $\Pc_{_2}(\R^n)\times\R^n$ by:
\beqs
{\rm Var}(\pi,k) & := & \int_{\R^n} (x-\bar\pi)\trans k(x-\pi) \pi(dx), \;\;\; \pi \in \Pc_{_2}(\R^n), \; k \in \S^n \\
v_2(\pi,\ell) & := & \bar\pi\trans \ell \bar\pi, \;\;\;  \pi \in \Pc_{_2}(\R^n), \; \ell \in \S^n \\
v_1(\pi,y) &:=& y\trans\bar\pi, \;\;\; \pi \in \Pc_{_2}(\R^n), \; y \in \R^n.
\enqs

We look for a value function solution to the Bellman equation \reff{HJBpi} in the form
\beq \label{wquadra1}
w(t,\pi) &=& {\rm Var}(\pi,K(t)) + v_2(\pi,\Lambda(t))  +  v_1(\pi,Y(t)) + \chi(t),
\enq
for some  functions $K$, $\Lambda$ $\in$ $C^1([0,T];\S^n)$, $Y$ $\in$ $C^1([0,T];\R^n)$, and $\chi$ $\in$ $C^1([0,T];\R)$ to be determined. One easily checks   that
$w$ lies in $\Cc_b^{1,2}([0,T]\times\Pc_{_2}(\R^n))$ with
\beqs
\partial_t w(t, \pi)  &=& {\rm Var}(\mu,K'(t)) + v_2(\pi,\Lambda'(t))  + v_1(\pi,Y'(t))  + \chi'(t), \\
\partial_\pi w(t, \pi)(x)&=& 2 K(t)(x-\bar\pi) + 2 \Lambda(t)\bar\pi  + Y(t), \\
\partial_x\partial_\pi w(t, \pi)(x)&=& 2 K(t),\\
\partial_\mu^2 w(t, \pi)(x,x') &=& 2(\Lambda(t)-K(t)).
\enqs
Together with the quadratic expression  of the running and terminal costs $f,g$  in \reff{costquadra},
and the linear form of the drift coefficient $b$, and diffusion coefficients
$\sigma_{_V}$, $\sigma_{_W}$ in \reff{dynXlin},
we then see after some tedious but direct calculations that $w$ satisfies the Bellman equation \reff{HJBpi} iff
\beq
{\rm Var}(\pi,K(T)) + v_2(\pi,\Lambda(T))  + v_1(\pi,Y(T)) + \chi(T)  &=&  {\rm Var}(\pi,P) + v_2(\pi,P), \label{HJBLQT1}
%+ (p_1 + \bar p_1).\bar\mu,
\enq
holds for all  $\pi$ $\in$ $\Pc_{_2}(\R^n)$, and
\beq
& & {\rm Var} \big(\pi,K'(t)+ Q + \sum_{i=1}^m (D_{_V}^i)\trans K(t)D_{_V}^i  + \sum_{j=1}^d (D_{_W}^j)\trans K(t) D_{_W}^j  + K(t) B + B\trans K(t) \big)   \nonumber \\
& & \;\;+\; v_2\big(\pi,\Lambda'(t) +Q + \sum_{i=1}^m (D_{_V}^i) \trans K(t)D_{_V}^i  +  \sum_{j=1}^d (D_{_W}^j) \trans \Lambda(t) D_{_W}^j   +   \Lambda(t)B + B\trans \Lambda(t) \big)  \nonumber \\
& & \;  +  \;  v_1\big(\pi,Y'(t)  + B\trans Y(t)    + \sum_{i=1}^m 2 (D_{_V}^i)\trans K(t) \gamma_{_V}^i  + 2 \sum_{j=1}^d  (D_{_W}^j) \trans\Lambda(t) \gamma_{_W}^j    + 2 \Lambda(t) b_0  \big)  \nonumber\\
& &  \;\;\;+ \;  \chi'(t)   \; + \;  Y(t)\trans b_0 + \sum_{i=1}^m (\gamma_{_V}^i) \trans K(t)\gamma_{_V}^i + \sum_{j=1}^d (\gamma_{_W}^j) \trans \Lambda(t) \gamma_{_W}^j +  \inf_{a \in  \R^q} G_t^{\pi}(a)  \nonumber \\
& =&  0, \label{HJBLQt1}
%& &+(\gamma'(t)+ q_1(t) +\bar q_1(t)+ \gamma(t)(B(t)+ \bar B(t)))\bar\mu +\chi'(t) =0
\enq
holds for all  $t$ $\in$ $[0,T)$, $\pi$ $\in$ $\Pc_{_2}(\R^n)$, where the function $G_t^\pi$ $:$ $\R^q$ $\rightarrow$ $\R$ is defined by
\beqs
G_t^\pi(a)&=&   a\trans \Gamma(K(t),\Lambda(t))a + \big[ 2 U(K(t),\Lambda(t))\trans\bar\pi + R(K(t),\Lambda(t),Y(t))]\trans a   \label{defGk}
 \enqs
 with
\begin{equation} \label{defUVSZ}
\left\{
\begin{array}{rcl}
\Gamma(k,\ell) &=& N +  \Sum_{i=1}^m (F_{_V}^i)\trans k F_{_V}^i +  \Sum_{j=1}^d  (F_{_W}^j)\trans \ell  F_{_W}^j \\
U(k,\ell) &=& \Sum_{i=1}^m (D_{_V}^i)\trans k F_{_V}^i + \Sum_{j=1}^d (D_{_W}^j)\trans \ell F_{_W}^j + \ell C \\
R(k,\ell,y) &=&   \Sum_{i=1}^m 2 (F_{_V}^i)\trans k \gamma_{_V}^i +   \Sum_{j=1}^d 2 (F_{_W}^j)\trans \ell \gamma_{_W}^j  +  C\trans y,
\end{array}
\right.
\end{equation}
for $(k,\ell)$ $\in$ $\S^n$, $y$ $\in$ $\R^n$.
Then, under the condition that the symmetric matrice $\Gamma(k,\ell)$ is positive definite in $\S^q$ (this will be checked later on), we get after square completion:
\beqs
G_t^\pi(a) &=&  (a-\hat a(\pi,K(t),\Lambda(t),Y(t))\trans\Gamma(K(t),\Lambda(t)) (a- \hat a(\pi,K(t),\Lambda(t),Y(t)) \\
& & \;\; -  \;  v_2\big(\pi,U(K(t),\Lambda(t))\Gamma^{-1}(K(t),\Lambda(t))U(K(t),\Lambda(t))\big) \\
& & \;\; - \; v_1\big(\pi, U(K(t),\Lambda(t))\Gamma^{-1}(K(t),\Lambda(t))R(K(t),\Lambda(t),Y(t)) \big) \\
& & \;\; - \;  \frac{1}{4} R\trans(K(t),\Lambda(t),Y(t))\Gamma^{-1}(K(t),\Lambda(t))R(K(t),\Lambda(t),Y(t))
\enqs
with
\beq \label{aopti}
\hat a(\pi,k,\ell,y)  &=&  -\Gamma^{-1}(k,\ell)\big[U\trans(k,\ell)\bar\pi + \frac{1}{2} R(k,\ell,y) \big]
\enq
This means that $G_t^\pi$ attains its infimum  at   $\hat a(\pi,K(t),\Lambda(t),Y(t))$, and plugging the above expression of
$G_t^\pi(\hat a(\pi,K(t),\Lambda(t),Y(t)))$ in \reff{HJBLQt1}, we observe that  the relation  \reff{HJBLQT1}-\reff{HJBLQt1}, hence the Bellman equation \reff{HJBpi},
is satisfied by identifying the terms in ${\rm Var}(\pi,.)$, $v_2(\pi,.)$,  $v_1(\pi,.)$,  which leads to the system of ordinary differential
equations (ODEs) for  $(K,\Lambda,Y,\chi)$:
\begin{equation} \label{linK}
\left\{
\begin{array}{rcl}
K'(t)+ Q + \Sum_{i=1}^m (D_{_V}^i)\trans K(t)D_{_V}^i  + \Sum_{j=1}^d (D_{_W}^j)\trans K(t) D_{_W}^j  + K(t) B + B\trans K(t) &=& 0 \\
K(T) &=& P,
\end{array}
\right.
\end{equation}
\begin{equation}\label{Riccatilambda}
\left\{
\begin{array}{rcl}
\Lambda'(t) +Q + \Sum_{i=1}^m (D_{_V}^i) \trans K(t)D_{_V}^i  +  \Sum_{j=1}^d (D_{_W}^j) \trans \Lambda(t) D_{_W}^j
+  \Lambda(t)B + B\trans \Lambda(t) & & \\
\;\;\; - \;  U(K(t),\Lambda(t))\Gamma^{-1}(K(t),\Lambda(t))U(K(t),\Lambda(t)) &=& 0 \\
\Lambda(T) \; &=& \;  P,
\end{array}
\right.
\end{equation}
\begin{equation} \label{linY}
\left\{
\begin{array}{rcl}
Y'(t)  + B\trans Y(t)    + \Sum_{i=1}^m 2 (D_{_V}^i)\trans K(t) \gamma_{_V}^i  + 2 \Sum_{j=1}^d  (D_{_W}^j) \trans\Lambda(t) \gamma_{_W}^j    + 2 \Lambda(t) b_0 \\
\;\;\; - \; U(K(t),\Lambda(t))\Gamma^{-1}(K(t),\Lambda(t))R(K(t),\Lambda(t),Y(t)) &=& 0,  \\
Y(T) \;& =& \;  0
\end{array}
\right.
\end{equation}
\begin{equation} \label{linchi}
\left \{
\begin{array}{rcl}
 \chi'(t)   \; + \;  Y(t)\trans b_0 + \Sum_{i=1}^m (\gamma_{_V}^i) \trans K(t)\gamma_{_V}^i + \Sum_{j=1}^d (\gamma_{_W}^j) \trans \Lambda(t) \gamma_{_W}^j & & \\
\;\;\; - \;  \frac{1}{4} R\trans(K(t),\Lambda(t),Y(t))\Gamma^{-1}(K(t),\Lambda(t))R(K(t),\Lambda(t),Y(t)) &=& 0, \\
\chi(T) \; &=& \;  0.
\end{array}
\right.
\end{equation}

Let us now show, under conditions {\bf (C1)} and {\bf (C2)}, the existence (and uniqueness) of a solution to the above system of ODEs.
Notice that  this system is decoupled:
\begin{itemize}
\item [(i)] One first considers  the linear ODE \reff{linK} for $K$, which is symmetric (in the sense that the ODE for $K\trans$ is the same).  Since the terminal condition
$P$ is nonnegative as well as the running term $Q$ under condition {\bf (C1)},  it is clear that there exists a unique solution $K$ $\in$ $C^1([0,T],\S^n)$ to \reff{linK}, and
this solution is nonnegative.
\item[(ii)] Given $K$, we next consider the ODE for $\Lambda$ with generator $\ell$ $\in$ $\S^n$ $\mapsto$
$Q + \sum_{i=1}^m (D_{_V}^i) \trans K(t)D_{_V}^i  +  \sum_{j=1}^d (D_{_W}^j) \trans \ell D_{_W}^j  +  \ell B + B\trans \ell$ $-$
$U_t(K_t,\ell,z)\Gamma_t^{-1}(K_t,\ell)U_t\trans(K_t,\ell,z)$  $\in$  $\S^n$, and terminal condition $P$.  This is a  Riccati  equation, and it is well-known (see e.g.
\cite{bis76}) that it is  associa\-ted to a stochastic standard LQ control problem  with controlled linear dynamics:
\beqs
d\tilde X_t &=& (B\tilde X_t +  C\alpha_t  ) dt + \sum_{j=1}^d (D_{_W}^j \tilde X_t + F_{_W}^i \alpha_t ) dW_t^j,
\enqs
and quadratic cost functional
\beqs
\tilde J^K(\alpha) &=& \E \Big[ \int_0^T \big( \tilde X_t\trans Q^K(t) \tilde X_t + \alpha_t\trans N^K(t)\alpha_t
+ 2\tilde X_t\trans M^K(t)\alpha_t \big) dt \; + \;  \tilde X_T\trans P\tilde X_T \Big],
\enqs
where $Q^K(t)$ $=$ $Q+ \Sum_{i=1}^m (D_{_V}^i)\trans K(t) D_{_V}$, $N^K(t)$ $=$ $N+\Sum_{i=1}^m F_{_V}\trans K(t) F_{_V}^i$,
$M^K(t)$ $=$ $\Sum_{i=1}^m (D_{_V})\trans K(t) F_{_V}$.
Under the condition that $N^K$ is positive definite, we can rewrite this cost functional after square completion as
\beqs
\tilde J^K(\alpha) &=& \E \Big[ \int_0^T \big( \tilde X_t\tilde Q^K(t) \tilde X_t + \tilde\alpha_t\trans N^K(t)\tilde\alpha_t \big) dt \;
+ \;  \tilde X_T\trans P \tilde X_T \Big],
\enqs
with $\tilde Q^K(t)$ $=$ $Q^K(t)-M^K(t)(N^K(t))^{-1}(M^K(t))\trans$, $\tilde\alpha_t$ $=$
$\alpha_t + (N^K(t))^{-1}(M^K(t))\trans\tilde X_t$. By noting that  $\tilde Q^K(t)$ $\geq$ $Q$, it follows that the symmetric matrices
$\tilde Q^K$ and $P$ are nonnegative under condition {\bf (C1)}, and assuming furthermore that $N^K$ is uniformly positive definite, it is known
from \cite{won70} that there exists a unique  solution $\Lambda$ valued in $\S^n$  to this Riccati equation, with $\Lambda$ being nonnegative.
This implies in particular that $\Gamma^{-1}(K(t),\Lambda(t))$ is well-defined.  Since $K$ is nonnegative under {\bf (C1)}, notice that the uniform positivity condition on $N^K$ is  satisfied under {\bf (C2)}: this is clear when $N$ is uniformly positive definite (as usually assumed in LQ problem), and holds also true
when  one of the $F_{_V}^i$ is  uniformly nondegenerate,  and $K$ is uniformly positive definite,  which occurs when $P$ or $Q$ is uniformly positive definite  from standard comparison principle for the linear ODE for $K$.
\item[(iii)]  Given $(K,\Lambda)$,  it is clear that there exists a unique solution $Y$ valued in $\R^n$ to the linear ODE \reff{linY}.
\item[(iv)] Finally, given $(K,\Lambda,Y)$, the equation for $\chi$ is explicitly solved as
\beqs
\chi_t &=& \int_t^T \big[Y(s)\trans b_0 + \Sum_{i=1}^m (\gamma_{_V}^i) \trans K(s)\gamma_{_V}^i + \Sum_{j=1}^d (\gamma_{_W}^j) \trans \Lambda(s) \gamma_{_W}^j   \\
& & \hspace{0.7cm} -   \frac{1}{4} R\trans(K(s),\Lambda(s),Y(s))\Gamma^{-1}(K(s),\Lambda(s))R(K(s),\Lambda(s),Y(s)) \big] ds, \;   \; 0 \leq t\leq T.
\enqs
\end{itemize}

By misuse of notation, we set: $\hat a(t,\bar\pi)$ $=$ $\hat a(\pi,K(t),\Lambda(t),Y(t))$, which is linear in $\bar\pi$,  and we see that there exists a unique solution $\hat X$
to the linear  stochastic McKean-Vlasov equation:
\beqs
d\hat X_t &=& (b_0 + B\hat X_t + \hat a(t,\E[\hat X_t|W]) dt + \sum_{i=1}^m (\gamma_{_V}^i +  D_{_V}^i \hat X_t + F_{_V}^i \hat a(t,\E[\hat X_t|W]) dV_t^i  \nonumber \\
& & \;\;\;+ \;     \sum_{j=1}^d (\gamma_{_W}^j +  D_{_W}^j \hat X_t + F_{_W}^j \hat a(t,\E[\hat X_t|W]) dW_t^j, \;\;\; 0 \leq t \leq T, \; \hat X_0 = x_0 \in \R^n.  \label{dynXopt}
\enqs
From the verification theorem \ref{theoverif}, we conclude that an optimal control for \reff{costquadra} is given  by
\beqs
\alpha_t^* &=& \hat a(t,\E[\hat X_t|W]), \;\;\; 0 \leq t \leq T, \\
&=&   -\Gamma^{-1}(K(t),\Lambda(t))\big[U\trans(K(t),\Lambda(t)) \E[\hat X_t|W] + \frac{1}{2} R(K(t),\Lambda(t),Y(t)) \big],
\enqs
while the optimal cost (at time $0$) is given by:
\beqs
\upv_0 &=& w(0,\delta_{x_0}) \; = \;  x_0\trans \Lambda(0)x_0 + Y(0)\trans x_0 + \chi(0).
\enqs

\vspace{5mm}

\small
\bibliographystyle{plain}
\bibliography{biblio7}

\end{document}